\documentclass[12pt]{amsart}
\usepackage{amsmath,amsfonts,amssymb,amscd,amsthm,latexsym}
\usepackage{pstricks,pst-node,pst-plot}
\usepackage{graphicx}

\newcommand{\cal}{\mathcal}

\begin{document}

\sloppy

\vskip 1cm

\centerline{\Large\bf Theoretic--model Properties of Regular
Polygons}
\bigskip

\centerline{\large Mikhalev A.V., Ovchinnikova E.V., Palyutin E.A.,
Stepanova A.A.} {\footnotetext[1]{This research was supported by
RFBR (grant 17-01-00531)}
\bigskip


\vskip 1cm {\bf\centerline {Abstract}}

 This work is dedicated to the results were got in the model
theory of the regular polygons. We give the characterization of the
monoids with axiomatizable and model complete class of regular
polygons. We describe the monoids with complete class of regular
polygons which satisfy the additional conditions. We study the
monoids, whose regular core is presented as a union of the finite
number of principal right ideals, all regular polygons over which
have the stable  and superstable theory. We prove the stability of
the axiomatizable model complete class of regular polygons and we
also describe the monoids with the superstable and $\omega$--stable
class of regular polygons when this class is axiomatizable and model
complete.

\sloppy  \vskip 1cm

This work is dedicated to the results were got in the model theory
of the regular polygons.

Under the left polygon ${}_SA$ over monoid $S$ or simply polygon
we understand a set $A$ upon which $S$ acts on left with the
identity of $S$ acting as the identity map on $A$. As it can be
seen from this definition we can look at the polygon over monoid
as the generalization of the notion of the module over ring. So
many notions and problems came to the model theory of the polygons
from the model theory  of the modules. In particular, the notion
of the regular polygon did. In  the module theory there are
several different definitions of the regular polygon. In the
polygon theory we are using the analogue of the regular of
Zelmanowitz module [Zel], introducing by Tran [Tra].

One of the standard problems in the model theory of the polygons
is a problem of monoids description, over which some class of
polygons possesses property $P$, where  $P$ can be the
axiomatizability, completeness, model completeness and others. In
the given work these questions are considered for regular polygon
class.

We tried to make the interpretation as  closed as possible, giving
all necessary definitions and statements, which became classical
already.

In the first two paragraphs we give information from the polygon
theory which will be necessary in future.

In the third paragraph we state the information from the model
theory.

In the fourth paragraph the characterization of the monoids with
axiomatizable class of regular polygons is given  (Theorem 4.1).
As a consequence the axiomatizability of the class of regular
polygons over group is received  (Corollary 5.1).

In the fifth paragraph the monoids with axiomatizable model
complete class of regular polygons are described (Theorem 5.1). In
this paragraph the model completeness of the class of regular
polygons over infinite group is proved  (Corollary 5.2).

In the sixth paragraph the questions of the completeness of the
class of regular polygons is regarded. Here we give the
characterization of monoids with complete class of regular
polygons which satisfy the additional conditions, namely, the
monoids over which the class of regular polygons satisfies the
condition of formula definability of isomorphic orbits (Theorem
6.1) and the linearly ordered monoids of depth 2 (Theorem 6.2).
Note that the question of the description of the monoids with
complete class of regular polygons is still open.

In the seventh and eighth paragraphs we study the monoids, whose
regular core is presented as a union of the finite number of
principal right ideals, all regular polygons over which have the
stable (Theorem 7.1) and superstable (Theorem 8.1) theory.

In the ninth paragraph we prove the stability of the axiomatizable
model complete class of regular polygons and we also describe the
monoids with the superstable and $\omega$--stable class of regular
polygons when this class is axiomatizable and model complete
(Theorem 9.1).

\vskip 1cm {\bf\centerline {\S~1. Polygons}}

\medskip

Throughout, by $S$ we denote a monoid, and by 1 unity in $S$. An
algebraic system $\langle A;s\rangle_{s\in S}$ of the language
$L_S=\{s\mid s\in S\}$ is a (left) $S$--polygon (or polygon over
$S$, or polygon) if $s_1(s_2a)=(s_1s_2)a$ and $1a=a$ for all
$s_1,s_2\in S$ and $a\in A$. A polygon $\langle A;s\rangle_{s\in
S}$ is denoted by ${}_SA$. All the polygons treated in the
article, unless specified otherwise, are left $S$--polygons.
Similarly the notion of the right $S$-polygon is defined. We
denote the class of all polygons by $S-Act$.

The subsystem ${}_SB$ of the polygon ${}_SA$ is called the
subpolygon of the polygon ${}_SA$. The polygon ${}_SA$ is called
finite generated if there exists $a_1,\ldots,a_n\in A$ such that
${}_SA=\bigcup\limits_{i=1}^n{}_SSa_i$. The polygon ${}_SA$ is
called cyclic if ${}_SA$ is one-generated polygon that is
${}_SA={}_SSa$ for some $a\in S$. The coproduct of the polygons
$_SA_i$ is a disjunctive union of this polygons. The coproduct of
the polygons $_SA_i$ we denote  by $\coprod\limits_{i\in
I}{}_SA_i$.

The following proposition we will often use without references.

\smallskip {\bf Proposition 1.1.} For any $a,e,f\in S,\;e^2=e,\;f^2=f$ we
have:

1) $aS\subseteq eS\Longleftrightarrow ea=a$;

2) $Sa\subseteq Se\Longleftrightarrow ae=a$.\hfill$\Box$

\smallskip {\bf Proposition 1.2.} Let $e,f\in S,\;e^2=e,\;f^2=f$,
$Se\subseteq Sf,\;fS\subseteq eS$. Then $e=f$.

 {\bf Proof.} Suppose $Se\subseteq Sf,\;fS\subseteq eS$, where $e,f\in S$ are
 the idempotents. Then from the first including by Proposition 1.1(2) it
 follows $e=ef$, from second one by Proposition 1.1(1) it
 follows $f=ef$. Consequently, $e=f$.\hfill$\Box$

\smallskip {\bf Proposition 1.3.} Let $T$ be a semigroup, $e\in T$ is idempotent.
Then the left ideal $Te$ is minimum by inclusion among the left
ideals of the semigroup $T$, generated by the idempotens if and
only if the right ideal $eT$ is minimum by inclusion among the
right ideals of the semigroup $T$, generated by the idempotens.

{\bf Proof. Necessity.} Suppose the left ideal $Te$ is minimum
among the left ideals of the semigroup $T$, generated by the
idempotens, $gT\subseteq eT,\;g^2=g$. Then $g=eg$. In view of
minimality of the ideal $Te$ the equality $Tge=Te$ is hold.
Consequently, $e=kge$ for some $k\in T$ and
$e=kge=kgge=kgege=ege=ge$. Thus $e\in gT$ and $gT=eT$.

{\bf Sufficiency} is proved similarly.\hfill $\Box$

A semigroup $T$ is said to be a rectangular band of  the groups
 if $T=\bigcup\{T_{ij}\mid i\in I(T),j\in J(T)\}$ is a
 decomposition of the semigroup $T$ on groups $T_{ij}$ and herewith
$T_{ij}\cdot T_{kl}\subseteq T_{il}$.

{\bf Remark 1.1.} If a semigroup $\langle T,\ast\rangle$ is a
rectangular band of groups $T_{ij}$ with units $e_{ij},$ $i\in
I(T),$ $j\in J(T),$ then for any $i,k\in I(T),$ $j,l\in J(T)$ the
following conditions are true:

1) $e_{ij}\cdot e_{kj}=e_{ij}$; $e_{ij}\cdot e_{il}=e_{il}$;

2) $e_{ij}\cdot T_{kl}=T_{ij}\cdot e_{kl}=T_{ij}\cdot
T_{kl}=T_{il}$;

3) $\langle T_{ij},\ast\rangle\cong\langle T_{kl},\ast\rangle$;

4) $Te_{ij}=\bigcup\limits_{p\in I(T)} T_{pi}$;
 $e_{ij}T=\bigcup\limits_{p\in J(T)} T_{ip}$;

5) for any $a\in T$ the condition $Ta=Te_{ij}$ ($aT=e_{ij}T$) is
equivalent to $a\in T_{pj}$ ( $a\in T_{ip}$) for some  $p\in I(T)$
($p\in J(T)$);

6) for any $a\in T$ the set $Ta$ ($aT$) forms a minimal left
(right) ideal.\hfill $\Box$

The union of all minimal left ideal of a semigroup $T$ is called
its kernel and will be denoted by $K(T)$.

  {\bf Proposition 1.4}
[Su]. If the kernel $K(S)$ contains an idempotent, then $K(S)$ is
a rectangular band of groups. \hfill $\Box$

The set of all idempotents of a semigroup $T$ will be denoted by
$E(T)$.

{\bf Proposition 1.5.} If $Te$ is a minimal left ideal of a
semigroup $T$ and $e\in E(T)$, then the set
$G_e\rightleftharpoons\{a\mid ea=ae=a\}$ forms a subgroup of the
semigroup $T$.

Proof. If $a,b\in G_e$, then $e(ab)=(ea)b=ab$ and
$(ab)e=a(be)=ab$. Consequently, the set $G_e$ is closed under a
semigroup operation. Obviously, the idempotent $e$ is an unit in
$G_e$.

Let $a$ be an arbitrary element of $G_e$. By $ae=a$ and the
minimality of left ideal $Se$ we get $Sa=Se$. Then there exists an
element $b\in S$ such that $ba=e$. Since $ebe\in G_e$ and
$(ebe)a=eb(ea)=eba=e$, then it is not difficult check that element
$a\cdot ebe=g$ is an idempotent. Since $g\in Se$, then $eg=ge=g$,
that is $g\in G_e$. In view of minimality $Se$ we have $Se=Sg$,
consequently, $eg=e$. Since $ea=a$ and $g=a\cdot ebe$, then
$eg=g$. Thereby, $g=e$, $a^{-1}=b$ and $G_e$ forms a subgroup in
the semigroup $S$. \hfill $\Box$

\vskip 1cm

\medskip {\bf\centerline {\S~2. Regular Polygons}}

\medskip

Let $_SA$ be a polygon. We call $a\in A$ an act--regular element
if there exists a homomorphism $\varphi:{}_SSa\longrightarrow
{}_SS$ such that $\varphi(a)a=a$, and $_SA$ is called a regular
polygon if every $a\in A$ is an act--regular element.

\smallskip {\bf Proposition 2.1} [KKM]. Let $_SA$ be a polygon.
The following conditions for  $a\in A$ are equivalent:

1) the element $a$ is an act--regular element;

2) there exists an idempotent $e\in S$ and an isomorphism
$\psi:{}_SSa\longrightarrow {}_SSe$ such that $\psi(a)=e$;

3) there exists an idempotent $e\in S$ such that ${}_SSa\cong
{}_SSe$, i.e. the polygon ${}_SSa$ is projective.

{\bf Proof.} 1)$\Rightarrow$2). Suppose $a\in A$ is act--regular
element, $\varphi:{}_SSa\longrightarrow {}_SS$ is a homomorphism
such that $\varphi(a)a=a$. Let $e=\varphi(a)$. Then
$e=\varphi(a)=\varphi(\varphi(a)a)=\varphi(a)\varphi(a)=e^2$ and
$ea=\varphi(a)a=a$. Furthermore an equality $sa=ta$ implies
$se=s\varphi(a)=\varphi(sa)=\varphi(ta)=t\varphi(a)=te$ for any
$s,t\in S$. Then the mapping $\psi:{}_SSa\longrightarrow {}_SSe$
such that $\psi(sa)=se$ for any $s\in S$, is an polygon
isomorphism.

The implication 2)$\Rightarrow$3) is obviously.

3)$\Rightarrow$1). Let $e$ be an idempotent of monoid $S$,
$\psi:{}_SSa\longrightarrow {}_SSe$ be an isomorphism,
$\psi(a)=ue,\;\psi(va)=e$. Then
$eva=va,\;e=\psi(va)=v\psi(a)=vue,\;\psi(ueva)=\psi(uva)=u\psi(va)=ve=\psi(a)$,
i.e. $e=vue$ and $ueva=a$. Suppose $f=uev\in
S,\;\varphi:Sa\longrightarrow S$ is the mapping such that
$\varphi(sa)=sf$ for any $s\in S$. The equality
$f^2=ue(vue)v=uev=f$ implies that $f$ is an idempotent. If $sa=ta$
for some $s,t$ then $s\psi(a)=t\psi(a),\;sue=tue$ and $sf=tf$,
consequently, $\varphi$ is the homomorphism. Since
$\varphi(a)a=ueva=a$, then $a$ is an act--regular
element.\hfill$\Box$

\smallskip {\bf Corollary 2.1.} The following conditions for a polygon ${}_SA$ are equivalent:

1) ${}_SA$ is a regular polygon;

2) for any $a\in A$ there exist an idempotent $e\in S$ and an
isomorphism $\psi:{}_SSa\longrightarrow {}_SSe$ such that
$\psi(a)=e$;

3) for any $a\in A$ there exists an idempotent $e\in S$ such that
${}_SSa\cong {}_SSe$.\hfill$\Box$

Recall that an element $a$ of semigroup $T$ is (von Neumann)
regular if $a=aba$ for some $b\in T$. The semigroup is called (von
Neumann) regular  if all its elements are regular.

\smallskip {\bf Proposition 2.2} [KKM]. If  $S$ is a (von Neumann) regular monoid then
 ${}_SS$ is a regular polygon but the converse is not true.

{\bf Proof.} Let $S$ be a (von Neumann) regular monoid, $a\in S$.
Then $a=aba$ for a some $b\in S$. Denote the element $ba$ of
monoid $S$ by $e$. Clear that $e$ is an idempotent. Furthermore,
$Sa=Saba=Sae\subseteq Se,\;Se=Sba\subseteq Sa$, that is $Sa=Se$.
Consequently, ${}_SS$ is a regular polygon.

On the other hand, let $S$ be a right cancellative monoid which is
not a group. Then for any $a\in S$ a polygon ${}_SSa$ is
isomorphous to a polygon ${}_SS$. Consequently, ${}_SS$ is a
regular polygon but $S$ is not a (von Neumann) regular
monoid.\hfill$\Box$

Let ${}_SA$ be a polygon which has a regular subpolygon. Note that
the union of all regular subpolygons of the polygon ${}_SA$ is
also a regular subpolygon. This subpolygon is called the regular
core of the polygon ${}_SA$ and denote by $R({}_SA)$. Instead
$R({}_SS)$ we will write ${}_SR$. A subsemigroup $R$ of monoid $S$
is called the the regular core of the monoid $S$. Hereinafter, we
assume that $R\neq \emptyset$. Denote the class all regular
polygons by ${{\mathfrak R}}$.

The elements $x,y$ of the polygon ${}_SA$ are called connected
(denoted by $x\sim y$) if there exist $n\in\omega$,
$a_0,\ldots,a_n\in A$, $s_1,\ldots,s_n\in S$ such that $x=a_0$,
$y=a_n$, and $a_i=s_ia_{i-1}$ or $a_{i-1}=s_ia_i$. The polygon
$_SA$ is called connected if  we have $x\sim y$ for any
$x,y\in{}_SA$. It is easy to check that $\sim$ is a congruence
relation on the polygon $_SA$. Let $B\subseteq A$. The elements
$x,y\in A\setminus B$ are called connected out of $B$, if there
exist $n\in\omega$, $a_0,\ldots,a_n\in A\setminus B$,
$s_1,\ldots,s_n\in S$ such that $x=a_0$, $y=a_n$, and
$a_i=s_ia_{i-1}$ or $a_{i-1}=s_ia_i$.

Let $_SA$ be a polygon. We will denote by Con$({}_SA)$ the lattice
of congruences of the polygon $_SA$, by $1_{{}_SA}$ and
$0_{{}_SA}$ unit and zero in the lattice Con$({}_SA)$ accordingly.
The congruence $\theta\in{\rm Con}({}_SA)$ is called an amalgam
congruence if $\theta\;\cap\sim\;=0_{{}_SA}$. Notice that amalgam
congruences identify elements, not connected with each other.

For an arbitrary class of polygons $K$ we will denote by

${\bf H}_A(K)$ the class of all polygons, isomorphic to
factor-polygons of polygons from $K$ by amalgam congruences;

${\bf S}(K)$ the class of all polygons, isomorphic to subpolygons
from $K$;

${\bf D}(K)$ the class of all polygons, isomorphic to coproducts
of polygons from $K$.

{\bf Proposition 2.3.} [Ov1] For any monoid $S$ the class of all
regular $S$-polygons $_S\Re$ coincides with the following classes:
${\bf H}_A{\bf DS}({}_SR)$, ${\bf SH}_A{\bf D}({}_SR)$, ${\bf
H}_A{\bf SD}({}_SR)$.\hfill $\Box$

A monoid $S$ is called regularly linearly ordered if for any $a\in
R$ the set $\{Sb\mid Sb\subseteq Sa\}$ is linearly ordered  by
inclusion.

{\bf Proposition 2.4.} 1) If $r\in R$, $e\in S$, $e^2=e$ and
$rS=eS$ then $e\in R$ and $rR=eR$.

2) If $e,f\in R$, $e^2=e$ and $f^2=f$ then the equality $eS=fS$ is
equivalent to the equality $eR=fR$.

{\bf Proof}. Suppose $r\in R,\;e\in S,\;e^2=e$ and $rS=eS$. From
the last equality we get $Sr\cong Se$. Since $r\in R$ then $e\in
R$. Since $e\in rS$ then $e=rt$, where $t\in S$. Since $e\in R$
then $te\in R$. Consequently, $e=rt=rte\in rR$, i.e. $eR\subseteq
rR$. In view of the inclusion $rS\subseteq eS$ we have $r=er\in eR
$ that is $rR\subseteq eR$. Thus, $rR=eR$ and 1) is proved.

Suppose $e,f\in R$, $e^2=e$ and $f^2=f$. If $eR=fR$ then $e=ee\in
eR=fR\subseteq fS$ and $eS\subseteq fS$; similarly, $fS\subseteq
eS$ that is $eS=fS$. If $eS=fS$ then $e=fe\in fR$ and $eR\subseteq
fR$; similarly, $fR\subseteq eR$ that is $eR=fR$ and 2) is
proved.\hfill$\Box$

{\bf Proposition 2.5}. If the monoid $S$ is regularly linear
ordered, ${}_SA\in {\mathfrak R}$ and $a\in A$ then the set
$\{Sb\;|\;Sb\subseteq Sa\}$ is linear ordered by inclusion.

{\bf Proof}. Suppose ${}_SA\in {\mathfrak R}$, $a\in A$ and
$b_1,b_2\in Sa$. Since ${}_SSa\in {\mathfrak R}$ then there exists
the isomorphism $\varphi:{}_SSa\longrightarrow{}_SSe$, where
$e^2=e\in R$. Hence in view of regular linear ordered of the
monoid $S$ either $S\varphi(b_1)\subseteq S\varphi(b_2) $ or
$S\varphi(b_2)\subseteq S\varphi(b_1)$. Consequently, either
$Sb_1\subseteq Sb_2 $ or $Sb_1\subseteq Sb_2 $. \hfill$\Box$

By the left depth (or simply depth) of the semigroup $T$ we call
the greatest length of chain of principal left ideals of this
semigroup if it exists and finite, and the symbol $\infty$
otherwise. We will denote the left depth of the semigroup $T$ by
${\rm ld}(T)$.

{\bf Proposition 2.6}. If the depth ${\rm ld}(R)$ of core of
monoid $S$ is finite, then the kernel $K(R)$ is a rectangular band
of groups.

{\bf Proof}. The finiteness of ${\rm ld}(R)$ implies the existence
of minimal left ideal $Sa$. In view of the regularity of the
polygon ${}_SR$ on Corollary 2.1 there exist an idempotent $e\in
R$ and an isomorphism $\psi:{}_SSa\longrightarrow {}_SSe$ such
that $\psi(a)=e$. Then the ideal $Se$ is also minimal and the
element $e$ belongs to $K(R)$. On Proposition 1.4 we get that
$K(R)$ is a rectangular band of groups. \hfill $\Box$

\vskip 1cm

\medskip {\bf\centerline {\S~3. Information from Model Theory}}

\medskip

The initial information from the model theory, used in this
article, may be found in [EP], [ChK] and [Sac]. We will remind
some of it.

Let us fix some complete theory $T$ of language $L$ and a rather
large and saturate model $\mathfrak C$ of the theory $T$, which we
call it as a monster-model, because we suppose that all considered
models of the theory $T$ are its elementary submodels. All
elements and sets will also be taken from the monster-model. All
formulae, considered in this paragraph,  will have a language $L$.

The finite sequences are called the corteges. A set of the
corteges of a set $A$ are denoted by $A^{<\omega}$. A length of a
cortege $\bar a$  is denoted by $l(\bar a)$. The corteges of the
length $n$ are called $n$--corteges. For the simplicity instead of
a denotement $\bar a\in A^{<\omega}$ we will often use a
denotement $\bar a\in A$. If $\Phi(\bar x,\bar y)$ is a formula of
a language $L$, $\bar a$ is a cortege of the elements and $l(\bar
a)=l(\bar y)$, then by $\Phi(C,\bar a)$ we will denote a set
$\{\bar b\mid C\models\Phi(\bar b,\bar a)\}$.

Class $K$ of the structures is called axiomatizable if there exist
a language $L$ and a set of the sentences $Z$ of the language $L$
such that for any structure $\cal A$
\begin{equation}\tag{3.1}
  {\cal A}\in
K\Longleftrightarrow(\mbox{the language of}\;{\cal
A}\;\mbox{is}\;L\;\mbox {and}\;{\cal A}\models\Phi\;\mbox{for
all}\;\Phi\in Z).
\end{equation}
If the condition (3.1) holds for the class $K$, then $L$ is called
the language of $K$, and the set $Z$ is called the set of the
axioms for $K$ (we denote it as $K=K_L(Z)$). If all structures of
the class $K$ have a language $L$, then the set of the sentences
of the language $L$, which are true in all structures from $K$, is
called an elementary theory of the class $K$ or simply theory of
the class $K$ and denoted by $Th(K)$. If $K=\{\cal A\}$ then we
will write $Th(\cal A)$ instead $Th(K)$.

We will say that the class $K$ of structures is closed under the
elementary equivalence (the isomorphism, the subsystems, the
ultraproduct and others) if with structures ${\cal A}_i,\;i\in I$,
it contains all structures which are elementary equivalent to them
(isomorphic to them, subsystems of them, ultraproduct of the
structures ${\cal A}_i$ and others.)

When we will study the axiomatizability of the classes we will use
the following criterion.

\smallskip{\bf Theorem 3.1} [EP]. Class $K$ of structures
of the language $L$ is axiomatizable if and only if $K$ is closed
under elementary equivalence and ultraproducts.\hfill$\Box$

The set of the sentences of the language $L$, closed under
deducibility, is called the elementary theory or simply theory of
the language $L$. The structure, in which all sentences of the
theory $T$ are true, is called a model of the theory $T$. The
consistent theory $T$ of the language $L$ is called complete if
$\Phi\in T$ or $\neg\Phi\in T$ for any sentence $\Phi$ of the
language $L$. A substructure ${\cal A}$ of a structure ${\cal B}$
is called elementary (it denoted by ${\cal A\prec B}$), if for any
formula $\varphi(\bar x)$ of the language $L$ and any $\bar b\in
B$
$${\cal A}\models\varphi(\bar b)\Leftrightarrow{\cal
B}\models\varPhi(\bar b).$$ Note that in this definition the
condition "$\Leftrightarrow$" we can exchange to "$\Leftarrow$"
(it is necessary to go to the negation of the formula).

The consistent theory $T$ of the language $L$ is called model
complete if $$\cal A\subseteq B\Longrightarrow A\prec B $$ for any
models ${\cal A,B}$ of theory $T$ of the language $L$. Theory $T$
of the language $L$ is called the theory with elimination of
quantifiers if any formula $\Phi$ of the language $L$ is
equivalent in $T$ some quantifier-free formula $\Psi$. Obviously,
the consistent theory with elimination of quantifiers is model
complete.

The formula of a form $\exists \bar x \Psi(\bar x;\bar y)$
($\forall\bar x \Psi(\bar x;\bar y)$) for a quantifier--free
formula $\Psi(\bar x;\bar y)$ is called existential (universal).

Let ${\cal A}$ be a structure of the language $L$, $X\subseteq A$,
the language $L_X$ is obtained from $L$ by adding a new constant
symbol to $L$ for each element of the set $X$, ${\cal A}_X$ be
enrichment of the structure ${\cal A}$ up to the language $L_X$
with the natural interpretation of the new constant symbols. A set
$D({\cal A})$ of the atomic sentences of the language $L_A$, which
are true in a structure ${\cal A}_A$, is called a diagram of the
structure ${\cal A}$.

\smallskip{\bf Theorem 3.2} [ChK]. If $T$ is a theory of the language $L$, then the following conditions are equivalent:

1) the theory $T$ is model complete;

2) if ${\cal A}$ and ${\cal B}$ are the models of the theory $T$
and ${\cal A\subset \cal B}$, then for any existential formula
$\psi(\bar y)$ of the language $L$ and any $\bar b\in A$
$${\cal B}\models\psi(\bar b)\Longrightarrow{\cal A}
\models\psi(\bar b);$$

3) for any formula $\varphi(\bar x)$ of the language $L$ there
exists an existential formula $\psi(\bar x)$, which is equivalent
in the theory $T$ to a formula $\varphi(\bar x)$.

{\bf Proof.} The statement $1\Rightarrow 2$ is trivial. Denote by
3' the condition which obtains from the condition 3 by exchange
"existential"  on "universal". It is clear that the conditions 3
and 3' are equivalent (it is enough to go to the negation of the
formula). Since the existential formulae are preserved under the
extending of structures and the universal formulae are preserved
under substructures, then the statement $3 \Rightarrow 1$ holds.

Let us proof the statement $2\Rightarrow 3$. Suppose the condition
2 holds and $\varphi(\bar x)$ is some formula of the language $L$.
We proof  the property 3 by induction on the number of the
quantifiers in the formula $\varphi(\bar x)$. Since the
existential quantifier is expressing through the universal
quantifier and negation, we will consider  that the formula
$\varphi(\bar x)$ does not contain the existential quantifier.
Assume $\varphi(\bar x)$ has a form $\forall y\psi(y,\bar x)$. By
the induction supposition and the equivalence of the conditions 3
and 3', the formula $\psi(y,\bar x)$ is equivalent in the theory
$T$ to some universal formula. Hence, the formula $\varphi(\bar
x)$ is also equivalent in the theory $T$ to some universal
formula. Let a set $Q$ consists of all existential formulae
$\xi(\bar x)$ of the language $L$ with the condition
$(T\cup\{\xi(\bar x)\})\vdash\varphi(\bar x)$ and a set $Q^*$
consists of the negations of the formulae from $Q$. If the set
$(Q^*\cup\{\varphi(\bar x)\})$ is not consistent with the theory
$T$, then the formula $\varphi(\bar x)$ is equivalent in $T$ to
the disjunction of some formulae from $Q$, so to some existential
formula, and the property 3 is proved.

Suppose the set $(Q^*\cup\{\varphi(\bar x)\}\cup T)$ is
consistent. Assume for some model $\cal A$ of the theory $T$ and
$\bar b\in \cal A$ we have ${\cal A}\models\chi(\bar b)$ for all
formulae $\chi(\bar x)\in(Q^*\cup\{\varphi(\bar x)\})$. We claim
that
$$(D({\cal A})\cup T)\vdash\varphi(\bar c_{\bar b}).$$
If it is wrong then there exists a model $\cal B$ of a set $(T\cup
D({\cal A})\cup\{\neg\varphi(\bar c_{\bar b})\})$. Since $\cal B$
is a model of a set $D(\cal A)$, then we can suppose that ${\cal
A}\subseteq{\cal B}$. Since $\neg\varphi$ is equivalent in $T$ to
the existential formula and ${\cal A}\models \varphi(\bar b)$,
then it contradicts to the condition 2. Since $D({\cal A})$
consists of the quantifier--free formulae then we have
$(\{\Phi(\bar c_{\bar a};\bar c_{\bar b})\}\cup T)\vdash
\varphi(\bar c_{\bar b})$ for some quantifier--free formula $\Phi$
and cortege $\bar a\in \cal A$. Consequently, we have ${\cal
A}\models\exists\bar z\Phi(\bar z;\bar b)$ and $\exists\bar
z\Phi(\bar z;\bar x)\in Q$, but all formulae from the set $Q^*$
are true in $\cal A$ on the cortege $\bar b$, contradiction.
\hfill$\Box$

Consistent theory $T$ of the language $L$ is called submodel
complete if the theory $T\cup D({\cal B})$ of the language $L_B$
is complete for any substructure ${\cal B}$ of any model ${\cal
A}$ of the theory $T$. Note that any submodel complete theory is
model complete.

\smallskip{\bf Theorem 3.3} [Sac]. A theory $T$ is submodel
complete if and only if $T$ is the theory with elimination of
quantifiers.

{\bf Proof. Necessity.} Let $T$ be a submodel complete theory of
the language $L$ and $\varphi(\bar x)$ be an arbitrary formula of
the language $L$. Let $Q$ be a set of all quantifier--free
formulae $\xi(\bar x)$ of the language $L$ with the property
$(T\cup\{\xi(\bar x)\})\vdash\varphi(\bar x)$ and $Q^*$ be a set
of all negations of formulae from $Q$. If a set
$(Q^*\cup\{\varphi(\bar x)\})$ is not consists with the theory
$T$, then the formula $\varphi(\bar x)$ is equivalents in $T$ to
disjunction of some formulae from the set $Q$, consequence, to
some quantifier--free formula, that the necessity is proved.

Suppose the set $(Q^*\cup\{\varphi(\bar x)\}\cup T)$ is
consistent. Assume for some model $\cal A$ of the theory $T$ and
$\bar b\in \cal A$ we have ${\cal A}\models\chi(\bar b)$ for all
formulae $\chi(\bar x)\in(Q^*\cup\{\varphi(\bar x)\})$. Let ${\cal
B}$ be a substructure of a structure $\cal A$, generated in $\cal
A$ by elements of the cortege $\bar b$. Since the theory $T$ is a
submodel complete, then one of two cases are hold: 1) $(D({\cal
B})\cup T)\vdash\varphi(\bar c_{\bar b})$, 2) $(D({\cal B})\cup
T)\vdash\neg\varphi(\bar c_{\bar b})$.

Suppose the first case holds. Since $D({\cal B})$ is consists of
quantifier--free formulae and a structure $\cal B$ is generated by
a cortege $\bar b$, then we have $(\{\Phi(\bar c_{\bar b})\}\cup
T)\vdash \varphi(\bar c_{\bar b})$ for some quantifier--free
formula $\Phi$. Therefore, ${\cal A}\models\Phi(\bar b)$ and
$\Phi(\bar x)\in Q$, but in $\cal A$ all formulae from the set
$Q^*$ are true on the cortege $\bar b$, contradiction.

The second case is also impossible, because a structure $\cal A$
is a model of the set $(D({\cal B})\cup T)$ and we have ${\cal
A}\models\varphi(\bar b)$.

{\bf Sufficiency.} Since any quantifier--free formula is
equivalent to the disjunction of formulae which are the
conjunctions of the atomic formulae and their negations, and a set
$D(\cal A)$ contains $\varphi$ or its negation for any atomic
sentence $\varphi$ of the language $L_A$, then for any
quantifier--free sentence $\Phi$ of the language $L_A$ we have
$D(\cal A)\vdash\Phi$ or $D(\cal A)\vdash\neg\Phi$. Thus, if $T$
is the theory with elimination of quantifiers, then $T$ is
submodel complete.\hfill$\Box$

If $K$ is a class of the structures of the language $L$ then by
$K^\infty$ we denote the class of infinite structures from $K$.
Class $K$ is called complete (model complete) if the theory
$Th(K^\infty)$ of the class $K^\infty$ is complete (model
complete).

Let $T$ be a consistent theory of the language $L$, $X={\{x_i\mid
1\leq i\leq n \}},$ $L_n=L_X$. A set of all formulae of language
$L$ with free variables from $X$ and parameters from $A\subseteq
\mathfrak C$ is denote by $F_X(A)$.  Any set of the sentences $p$
of the language $L_n$ is called $n$--type of the language $L$. If
the theory $p\cup T$ is consistent, then $p$ is called $n$--type
over $T$. If $p$ is a complete theory, then $p$ is called complete
$n$--type of the language $L$. If additionally $T\subseteq p$,
then $p$ is called a complete $n$--type over $T$. A set of all
complete $n$--types over $T$ is denote by $S_n(T)$.

Let ${\cal A}$ be a structure of the language $L$, $X\subseteq A$,
$a\in A$. A set $tp(a,X)=\{\Phi(x)\mid{\cal A}_X\models\Phi(a)\}$
is called a type of an element $a$ over a set $X$. It is not
difficult to understand that $tp(a,X)$ is a complete $1$--type
over $Th({\cal A}_X)$. By $S_n({X})$ we denote $S_n(Th({\cal
A}_X))$. Often we will write $S({X})$ instead $S_1({X})$.

The theory $T$ is called stable in a cardinal $\kappa$ or
$\kappa$--stable if $|S(X)|\leq \kappa$ for any model ${\cal A}$
of the theory $T$ and any $X\subseteq A$ of cardinal $\kappa$. If
the theory $T$  is $\kappa$--stable for some infinite $\kappa$,
then $T$ is called stable. If the theory $T$ $\kappa$--stable for
all $\kappa\geq2^{|T|}$, then $T$ is called superstable. If the
theory $T$ is not stable, then $T$ is called unstable.

{\bf Theorem 3.4} [She]. A complete theory is unstable if and only
if there exists a formula $\Phi(\bar x,\bar y)$ on $2n$ variables,
a model ${\cal A}$ of the theory $T$ and $\bar a_i\in
A^n,\;i\in\omega$, such that for any $i,j,\;i\neq j$,
$$i<j\;\Longleftrightarrow\;{\cal A}\models\Phi(\bar a_i,\bar
a_j).$$\hfill$\Box$

{\bf Theorem 3.5}. If the theory $T$ is stable in a countable
cardinality ($\omega$-stable), then it stable in all infinite
cardinality.

{\bf Proof.} Let the theory $T$ is $\omega$-stable. Assume there
exists a subset $A$ of monster-model $\mathfrak C$ of the
cardinality $\lambda\geqslant\omega$, such that the cardinality of
the set $S(A)$ strictly more then a cardinal $\lambda$. From
$\omega$-stability it follows that there exists the countable
language $L'\subseteq L$, such that for any predicate or function
of the language $L$ there exists equivalent in theory $T$
predicate or function  accordingly of the language $L'$. So we can
consider that the language $L$ is countable. Then the set
$F_X(A)$, where $X=\{x\}$, has a cardinality $\lambda$. Let $S(A)$
be the Stone space of $A$. Remind that this space is defined by
the basis of open sets $\{U_\Phi\mid \Phi\in F_X(A)\}$, where
$U_\Phi=\{t\mid \Phi\in t\in S(A)\}$. For the topological space
$X$ we denote by $X'$ the derived space, that is the space, which
is obtained from the space $X$ removing the isolated points. By
induction on ordinal $\alpha$ we define the subspaces
$S^{(\alpha)}$ as follows: $S^{(0)}=S(A)$,
$S^{(\beta+1)}=(S^{(\beta)})'$ and
$S^{(\delta)}=\bigcap\{S^{(\beta)}\mid\beta<\delta\}$ for limited
$\delta$. Let $\gamma$ be the least ordinal, for which the
condition $(S^{(\gamma)})'=S^{(\gamma)}$ holds. Clearly, that for
any formula $\Phi\in F_X(A)$ there exists not more than one
ordinal $\beta<\gamma$ such that the set $(U_\Phi\cap
S^{(\beta)})$ consists of one point. So the cardinality of the
ordinal $\gamma$ is not more than the cardinality of the set
$F_X(A)$, which is equal to $\lambda$. Since the number of
isolated points of the space, whose basis of open sets has the
cardinality $\lambda$, is also not more than $\lambda$, and the
cardinality of $S(A)$ is strictly more than $\lambda$, it follows
that the cardinality of $S^{(\gamma)}$ is also more than
$\lambda$. Since the space $S^{(\gamma)}$ has not isolated points
and is Hausdorff, then for any cortege $\varepsilon\in
2^{<\omega}$ of elements of the set $\{0,1\}$ there exist nonempty
sets $X(\varepsilon)$ of the form $\Phi_\varepsilon(\mathfrak
C;\bar a_\varepsilon)$ for the formula $\Phi(x;\bar a)\in F_X(A)$
with the following properties:

1) $X(\varnothing)=\mathfrak C$;

2) $X(\varepsilon^\wedge L)\subseteq X(\varepsilon)$,
$L\in\{0,1\}$;

3) $(X(\varepsilon^\wedge 0)\cap X(\varepsilon^\wedge
1)=\varnothing$.

Write $A_0$ for the countable set $=\bigcup\{\bar
a_\varepsilon\mid \varepsilon\in 2^{<\omega}\}$. By the properties
1--3 the cardinality of the set $S(A_0)$ is equal to $2^\omega$
which contradicts with the $\omega$-stability of the theory
$T$.\hfill$\Box$

In [Mus] T.G.Mustafin gives the notion of stationary theory of
polygons, which we will use hereinafter. A complete theory $T$ of
polygons is called stationary if for any ${}_SM\models T$ and
$a,b\in{\cal C}\setminus M$
$$
a\in Sb\Longrightarrow M\cap Sa=M\cap Sb.
$$
Let ${}_SA$ be a polygon and $a\in {\mathfrak C}\setminus A$. An
element $c\in A$ is called input element from $a$ in ${}_SA$ if
$c\in Sa$ and $Sb\subseteq Sc$ for all $b\in A\cap Sa$.

{\bf Theorem 3.6} [Mus]. Let $T$ be a stationary theory,
${}_SM\models T$, $a,b\in\mathfrak C\setminus M$ and $c$ be an
input element from $a$ in ${}_SM$. The following conditions are
equivalent:

1) $tp(a,M)=tp(b,M)$;

2) $c$ is an input element from $b$ in ${}_SM$ and
$tp(a,\{c\})=tp(b,\{c\})$;

3) $M\cap Sa=M\cap Sb$ и $tp(a,(M\cap Sa))=tp(b,M\cap Sa)$.

{\bf Proof.} $1\Rightarrow 2$. Assume the condition 1 is
satisfied. The equality $tp(a,\{c\})=tp(b,\{c\})$ is obviously.
Let us prove $M\cap Sa=M\cap Sb$. Suppose $m\in M\cap Sa$. Then
$m=sa$ for some $s\in S$. Hence, $m=sx\in tp(a,M)$. Therefore,
$m=sx\in tp(b,M)$, that is $m=sb\in M\cap Sb$. Thus, $M\cap
Sa\subseteq M\cap Sb$. Similarly  $M\cap Sb\subseteq m\cap Sa$.
Consequently,  $M\cap Sa=M\cap Sb$. Since  $c$ is an input element
from $a$ in ${}_SM$ then $c$ is an input element from $b$ in
${}_SM$.

$2\Rightarrow 3$. Assume the condition 2 is satisfied. The
equalities $M\cap Sa=Sc$ and $M\cap Sb=Sc$ follow from the
definition of the input element. Hence, $M\cap Sa=M\cap Sb$. The
equality  $tp(a,\{c\})=tp(b,\{c\})$ implies the existence of an
identity on $\{c\}$ automorphism $\varphi$ of the polygon
$\mathfrak C$ such that $\varphi(a)=b$. If $d\in Sc$ and $d=sc$,
where $s\in S$, then $\varphi(d)=\varphi(sc)=s\varphi(c)=sc=d$.
Therefore, $\varphi$ is identity on $Sc$. Thus, $tp(a,(M\cap
Sa))=tp(b,M\cap Sa)$.

$3\Rightarrow 1$. Assume the condition 3 is satisfied. The
equality  $tp(a,M\cap Sa)=tp(b,M\cap Sa)$ implies the existence of
an identity on $M\cap Sa$ automorphism $\varphi$ of the polygon
$\mathfrak C$ such that $\varphi(a)=b$. Denote the set $\{d\in
\mathfrak C\mid d \;\mbox{connected with}\; u\; \mbox{out of}
\;M\}$ by $C_M(u)$, where $u\in \mathfrak C$.

We claim that $\varphi(C_M(a))\subseteq C_M(b)$. Let
$a_0,\ldots,a_n\in\mathfrak C\setminus M$, $a=a_0$ and $a_i\in
Sa_{i+1}$ or $a_{i+1}\in Sa_i$ for all $i,\;0\leq i\leq n-1$. By
the induction on $n$ let us prove that $\varphi(a_n)\in C_M(b)$.
For $n=0$ the statement follows from the condition $b\in\mathfrak
C\setminus M$. Suppose $\varphi(a_i)\not\in M$. If $a_i\in
Sa_{i+1}$, then $\varphi(a_i)\in S\varphi(a_{i+1})$ and
$\varphi(a_{i+1})\not\in M$. Assume $a_{i+1}\in Sa_i$ and
$\varphi(a_{i+1})\in M$. The definition of the stationary theory
implies $M\cap S\varphi(a_{i+1})=M\cap Sb$. Then
$\varphi(a_{i+1})\in M\cap Sb=M\cap Sa$. Since $\varphi$ acts on
$M\cap Sa$ identically, then $a_{i+1}\in M\cap Sa$, contradiction.
Thus, $\varphi(C_M(a))\subseteq C_M(b)$. Similarly,
$\varphi^{-1}(C_M(b))\subseteq C_M(a)$, that is $\varphi(C_M(a))=
C_M(b)$. The equality $\varphi(C_M(b))\subseteq C_M(a)$ is proved
the same way.

We construct the automorphism $\psi$ of polygon ${}_S\mathfrak C$
as follows: $\psi|(C_M(a)\cap C_M(b))=\varphi |(C_M(a)\cap
C_M(b))$ and $\psi$ is an identity mapping on the set $\mathfrak
C\setminus (C_M(a)\cap C_M(b))$. Clearly, that $\psi(a)=b$. Hence,
the condition 1 of lemma holds.\hfill$\Box$

{\bf Теорема 3.7} [Mus]. Each stationary theory is stable.

{\bf Proof.} Let ${}_SM$ be a model of the theory $T$,
$|M|=\kappa$, $\kappa=\kappa^{|T|}$. If $a\in\mathfrak C$, then
$|Sa|\leqslant|S|\leqslant|T|$ and $|M\cap Sa|\leqslant|T|$. So on
Lemma 3.1
$|S_1(M)|\leqslant|M|^{|T|}\cdot2^{|T|}=|M|^{|T|}=\kappa^{|T|}=\kappa$.
Hence, $T$ is $\kappa$--stable theory.\hfill$\Box$

Let $K$ be a class of polygons. Monoid $S$ is called
$K$--stabilizer ($K$--superstabilizer, $K$--$\omega$--stabilizer)
if $Th({}_SA)$ is stable (su\-per\-stable, $\omega$--stable
accordingly) for any polygon ${}_SA\in K$. If $K=S-Act$, then
$K$--stabilizer ($K$--su\-per\-stabilizer,
$K$--$\omega$--stabilizer) are called stabilizer
(su\-perstabilizer, $\omega$--stabilizer accordingly).

Let us give the characterization of stabilizer and
superstabilizer, which was got by T.G. Mustafin. A monoid $S$ is
called linearly ordered if the set $\{Sa\mid a\in S\}$ is linearly
ordered  by inclusion.

{\bf Theorem 3.8} [Mus]. Monoid $S$ is stabilizer if and only if
$S$ is linear ordered monoid.\hfill$\Box$

A linearly ordered monoid $S$ is called well--ordered if it
satisfies the ascending chain condition for principal left ideals.

{\bf Theorem 3.9} [Mus]. Let $S$ be a countable monoid. Monoid $S$
is superstabilazer if and only if $S$ well--ordered
monoid.\hfill$\Box$

\vskip 1cm
\medskip {\bf\centerline {\S~4. Axiomatizability of Class for Regular Polygons}}
\medskip

The main result in this paragraph is Theorem 4 which give the
characterization  for the monoids with the axiomatizable class of
the regular polygons. In particular we get the axiomatizability of
the class of regular polygons over the group (Corollary 4.3).

\smallskip {\bf Theorem 4.1} [Ste1].  Class ${{\mathfrak R}}$ for the regular
polygons axiomatizable is if and only if

1) the semigroup $R$ is satisfied the descending chain condition
for principal right ideals which are generated by the idempotens;

2) for any $n\geqslant 1,\;s_i,t_i\in S\;(1\leqslant i\leqslant
n)$ the set $\{x\in R\mid\bigwedge\limits_{i=1}^ns_ix=t_ix\}$ is
empty or finite generated as a right ideal of the semigroup $R$.

{\bf Proof. Necessity.} Let ${{\mathfrak R}}$ be an axiomatizable
class. Assume condition 1 is not hold. This means that there
exists a decreasing sequence of principal right ideals:
$$f_1S\supset f_2S\supset\ldots\supset f_nS\supset\ldots,$$
where $f_n\in R,\;f_n^2=f\;(n\geqslant 1)$. For any
$n,m,\;1\leqslant n\leqslant m$, the inclusion $f_nS\supseteq
f_mS$ implies the equality $f_nf_m=f_m$. Suppose $\bar
f=(f_n)_{n\in\omega}\in R^\omega$ and $D$ is an arbitrary
non--principal ultrafilter on $\omega$. Then  the equality
$f_n\cdot \bar f/D=\bar f/D$ is true in ${}_SR^\omega/D$ for any
$n\geqslant 1$. In view of the axiomatizability of the class
${{\mathfrak R}}$ by Theorem 3.1 we have
${}_SR^\omega/D\in{\mathfrak R}$. On Corollary 2.1 there exist an
idempotent $e\in R$ and an isomorphism $\varphi:{}_S(S\cdot\bar
f/D)\longrightarrow {}_SSe$ such that $\varphi(\bar f/D)=e$. Then
$e\cdot \bar f/D=\bar f/D$. The equality $f_ne=e$ implies the
equality $f_n\cdot\bar f/D=\bar f/D$ for any $n\geqslant 1$.
Consequently, there exists $m\geqslant 1$ such that $f_m=ef_m\in
eS\subseteq f_nS$ for any $n\geqslant 1$ that contradict to the
condition $f_{m+1}S\subset f_mS$. Thus condition 1 is proved.

Assume condition 2 is not hold. Then there exist $n\geqslant
1,\;s_i,t_i\in S\;(1\leqslant i\leqslant n)$ such that $X=\{x\in
R\mid\bigwedge\limits_{i=1}^ns_ix=t_ix\}$ is the non--empty set
and is not finite generated as a right ideal of $R$. So there are
the infinite ordinal $\gamma$ and $x_\tau\in X\;(\tau<\gamma)$
such that $X=\cup\{x_\tau R\mid\tau<\gamma\}$ and $x_\beta
R\not\subseteq\cup\{x_\tau R\mid \tau<\beta\}$ for all
$\beta<\gamma$. Suppose $\bar x=(x_\tau)_{\tau<\gamma}\in
R^\gamma$ and $D$ is the ultrafilter on $\gamma$ such that
$|Y|=\gamma$ for $Y\in D$. In vier of the axiomatizability of the
class ${{\mathfrak R}}$ by Theorem 3.1 we get
${}_SR^\gamma/D\in{\mathfrak R}$. On Corollary 2.1 there exist an
idempotent $e\in R$ and an isomorphism $\varphi:{}_SS\bar
x/D\longrightarrow {}_SSe$ such that $\varphi(\bar x/D)=e$. Since
$x_\tau\in X\;(\tau<\gamma)$ we have
$\bigwedge\limits_{i=1}^ns_i\bar x/D=t_i\bar x/D$ and $e\in X$.
Consequently, $eR\subseteq\cup\{x_\tau R\mid \tau<\gamma\}$ that
is $eR\subseteq x_{\tau_0}R$ for some $\tau_0<\gamma$. Since
$e=ee$ it follows that $\bar x/D=e\cdot\bar x/D$. In particular,
$x_\tau\in eR$ for some $\tau>\tau_0$ and $x_\tau R\subseteq
x_{\tau_0}R $. We get the contradiction. Thus condition 2 is
proved.

 {\bf Sufficiency.} Assume conditions 1,2 of this theorem
 hold.
 Suppose
 $n\geqslant 1$, $\bar s=\langle s_1,\ldots,s_n \rangle$, $\bar
t=\langle t_1,\ldots,t_n \rangle\in S^n$, $X_{\bar s\bar t}=\{x\in
R\mid\bigwedge\limits_{i=1}^ns_ix=t_ix\}$. Let us show that either
$X_{\bar s\bar t}=\emptyset$ or $X_{\bar s\bar t} =\cup\{e_iR\mid
1\leqslant i\leqslant k\}$ for some $k\geqslant 1 $ and
idempotents $e_i\in X_{\bar s\bar t}$ $(1\leqslant i\leqslant k)$.
Suppose $X_{\bar s\bar t}\neq\emptyset$. Under condition of the
theorem there exist $k\geqslant 1,\;r_i\in X_{\bar s\bar
t}\;(1\leqslant i\leqslant k)$ such that $X_{\bar s\bar
t}=\cup\{r_iR\mid 1\leqslant i\leqslant k\}$. We can consider that
$r_iR\not\subseteq r_jR\;(i\neq j)$. Fix $i,\;1\leqslant
i\leqslant k$. Since $r_i\in R$, on Corollary 2.1 there exist an
idempotent $e_i\in R$ and an isomorphism
$\varphi:{}_SSr_i\longrightarrow {}_SSe_i $ such that
$\varphi(r_i)=e_i$. Then $e_ir_i=r_i$. Since $r_i\in X_{\bar s\bar
t}$ we have $e_i\in X_{\bar s\bar t}=\cup\{r_iR\mid 1\leqslant
i\leqslant k\}$, that is, $e_i=r_js$ for some $j,1\leqslant
j\leqslant k,$ and $s\in R$. Consequently, $r_i=e_ir_i=r_jsr_i\in
r_jR$. In vier of the choice of the element $r_j\;(1\leqslant
j\leqslant k )$ this means that $r_i=r_j$. Since $r_i=e_ir_i$ it
follows that $r_i\in e_iR$. In view of $e_i=r_is$ we have $e_i\in
r_iS$. Consequently, $r_iS=e_iS$. On Proposition 2.4 (2),
$r_iR=e_iR$. Thereby, $X_{\bar s\bar t}=\cup\{e_iR\mid 1\leqslant
i\leqslant k\}$, where $e_i\in X_{\bar s\bar t}$.

Define a set of formulae $\Gamma$ as follows: for all $n\geqslant
1,\;\bar s=\langle s_1,\ldots,s_n \rangle$, $\bar t=\langle
t_1,\ldots,t_n \rangle\in S^n$
$$\neg\exists x\bigwedge\limits_{i=2}^ns_ix=t_ix\in\Gamma,\;
\mbox{if}\; {}_SR\models\neg\exists x(x\in X_{\bar s\bar t});$$
$$\forall x(\bigwedge\limits_{i=2}^ns_ix=t_ix\longrightarrow\bigvee
\limits_{j=1}^kx=e_jx)\in\Gamma,\; \mbox{if}\; {}_SR\models\exists
x(x\in X_{\bar s\bar t}),$$ where $X_{\bar s\bar t}=\{x\in
R\mid\bigwedge\limits_{i=1}^ns_ix=t_ix\}=\cup\{e_jR\mid1\leqslant
j\leqslant k\},\;e_j^2=e_j\in X_{\bar s\bar t}$. Let us show that
$${}_SA\in{{\mathfrak R}}\Longleftrightarrow{}_SA\models\Gamma.$$

Let ${}_SA\in{{\mathfrak R}}$. Suppose
$\bigwedge\limits_{i=1}^ns_ia=t_ia$ for some $a\in A$. On
Corollary 2.1 there exist an idempotent $f\in R$ and an
isomorphism $\varphi:{}_SSa\longrightarrow {}_SSf$ such that
$\varphi(a)=f$. Then $f\in X_{\bar s\bar t} $. Consequently,
$X_{\bar s\bar t}\neq\emptyset$. Assume $X_{\bar s\bar
t}=\cup\{e_jR\mid1\leqslant j\leqslant k\}$ and $e_j^2=e_j\in
X_{\bar s\bar t}$. Then ${}_SR\models\bigvee\limits_{j=1}^ke_jf=f$
and ${}_SA\models\bigvee\limits_{j=1}^ke_ja=a$.

Suppose ${}_SA\models\Gamma,\;a\in A$. Let us prove that
${}_SSa\cong {}_SSe$, where $e$ is some idempotent from $R$.
Suppose $sa=ta,\;s,t\in S$. Since ${}_SA\models\Gamma$ we have
${}_SR\models\exists x(sx=tx)$ and $a=fa$ for some idempotent
$f\in R$. Suppose  $\{f_\tau\mid \tau<\gamma\}=\{f\mid
f^2=f,\;fa=a,\;f\in R\}$. By induction on $\gamma$ we will show
that there exists $\gamma_0<\gamma$ such that
$$f_{\gamma_0}S=\cap\{f_\tau S\mid \tau<\gamma\}.$$

Let $\gamma$ be a limit ordinal, $\tau_0<\gamma$. Under suggestion
of the induction there exists $\beta_0<\tau_0$ such that
$f_{\beta_0}S=\cap\{f_\tau S\mid \tau<\tau_0\}$. If
$f_{\beta_0}S\neq\cap\{f_\tau S\mid \tau<\gamma\}$  then there are
$\beta_1,\tau_1,\;\beta_1<\tau_1<\gamma$ such that
$$f_{\beta_0}S\supset f_{\beta_1}S=\cap\{f_\tau S\mid
\tau<\tau_1\}$$ and etc. Since $f_{\beta_n}S\supset
f_{\beta_{n+1}}S$ we have
$f_{\beta_{n+1}}=f_{\beta_{n}}f_{\beta_{n+1}}$. Consequently,
$f_{\beta_n}R\supseteq f_{\beta_{n+1}}R$ and on Proposition 2.4
(2) $f_{\beta_n}R\supset f_{\beta_{n+1}}R,\;n\geqslant 0$. Under
condition of the theorem the decreasing chain of ideals
$f_{\beta_0}R\supset f_{\beta_1}R\supset\ldots\supset f_{\beta_n}
R\supset\ldots$ is become stabilize. On Proposition 2.4 (2)
 the decreasing chain of  ideals $f_{\beta_0}S\supset
f_{\beta_1}S\supset\ldots\supset f_{\beta_n} S\supset\ldots$ is
become stabilize also that is $f_{\beta_k}S=\cap\{f_\tau S\mid
\tau<\gamma\}$ for some $k\geqslant 0$.

Let $\gamma$ be a non--limit ordinal. Assume there exists
$\beta_0<\gamma-1$ such that $f_{\beta_0}S=\cap\{f_\tau
S\mid\tau<\gamma -1\}$. Then ${}_SA\models a=f_{\beta_0}a\wedge
a=f_{\gamma-1} a$. Since ${}_SA\models\Gamma$ we have
${}_SR\models\exists x(x=f_{\beta_0}x\wedge x=f_{\gamma-1} x)$ and
there exists $f\in R$ such that
$a=fa,\;f=f_{\beta_0}f=f_{\gamma-1} f$. Consequently,
$$f=f_{\gamma_0},\;\gamma_0<\gamma,\;f_{\gamma_0}S\subseteq
f_{\beta_0}S\cap f_{\gamma-1} S,\;f_{\gamma_0}S=\cap
\{f_\tau\mid\tau<\gamma\}.$$ We put $e=f_{\gamma_0}$. Then $ea=a$
and the equality $ga=a$ implies $eS\subseteq gS$ for any
idempotent $g\in R$ that is $e=ge$. Let us show that the mapping
$\varphi:Sa\longrightarrow Se$ such that $\varphi(sa)=se$ for any
$s\in S$, is a polygon isomorphism. Suppose $ra=ka,\;r,k\in S$.
Since ${}_SA\models \Gamma$ then there exists an idempotent $g\in
R$ such that $rg=kg$ and $ga=a$. So $g e=e$ and $re=ke$. Suppose
$re=ke,\;r,k\in S$. Since $ea=a$ we get $ra=ka$. Thus, in view of
arbitrary of the choice of the element $a$, ${}_SSa\cong {}_SSe$
and ${}_SA\in{\mathfrak R}$.\hfill$\Box$

From proof of sufficiency it follows

\smallskip {\bf Corollary 4.1.} Let ${\mathfrak R}$ be an axiomatizable
class and  $X_{\bar s\bar t}=\{x\in
R\mid\bigwedge\limits_{i=1}^ns_ix=t_ix\}$ be a non--empty set,
where $n\geqslant 1$, $\bar s=\langle s_1,\ldots,s_n \rangle$,
$\bar t=\langle t_1,\ldots,t_n \rangle\in S^n$. Then the set
$X_{\bar s\bar t}$ is finite generated as a right ideal of the
semigroup $R$ if and only if $X_{\bar s\bar t} =\cup\{e_iR\mid
1\leqslant i\leqslant k\}$ for some $k\geqslant 1 $ and some
idempotants $e_i\in X_{\bar s\bar t}$ $(1\leqslant i\leqslant
k)$.\hfill$\Box$

\smallskip {\bf Corollary 4.2}. If the class ${\mathfrak R}$ of the regular polygons is
axiomatizable then $R=\bigcup\{e_iR\mid 1\leqslant i\leqslant n
\}$ for some $n\geqslant 1$, $e_i\in R$, $e^2_i=e_i$ $(1\leqslant
i\leqslant n)$.

{\bf Proof} follows from Theorem 4.1, Corollary 4.1 and the
equality $R=\{x\in R\mid x=x\}$. \hfill$\Box$

The following statement is obviously.

\smallskip {\bf Corollary 4.3}. The class ${\mathfrak R}$ for regular polygons over
group is axiomatizable.\hfill$\Box$

\vskip 1cm

\medskip{\bf\centerline {\S~5. Model Completeness of Class for
Regular Polygons}}

\medskip

In this paragraph it is formulated and proved the criterion of the
model completeness for an axiomatizable class of regular polygons
(Theorem 5.2). As a consequence we derive a model completeness of
the class of regular polygons over infinite group (Corollary 5.1).

\smallskip {\bf Lemma 5.1}. Let a monoid $S$ be regularly linearly ordered,
${}_SB\subseteq{}_SA\in{{\mathfrak R}},\;a_i\in A$ $(1\leqslant
i\leqslant k)$. Then for any $i,\;1\leqslant i\leqslant k$, the
following conditions hold:

1) $\cap\{Sa_j\mid Sa_j\cap Sa_i\neq\emptyset\}\neq\emptyset$;

2) if $Sa_i\cap B\neq\emptyset$ then $B\cap\bigcap\{Sa_j\mid
Sa_j\cap Sa_i\neq\emptyset\}\neq\emptyset$.

{\bf Proof}. Suppose that the conditions of the lemma hold. Let us
prove the statement 2 (the statement 1 is proved similarly).
Assume that $\{j_0,\ldots,j_s\}=\{j\mid Sa_j\cap
Sa_i\neq\emptyset\}$, where $j_0=i$;
$$B_n=B\cap\bigcap\{Sa_j\mid j\in\{j_0,\ldots,j_n\}\}\;\;(n\leqslant s).$$
It is sufficiently to show that for any $n,\;n<s$, the inequality
$B_n\neq\emptyset$ implies the inequality $Sa_{j_{n+1}}\cap
B_n\neq\emptyset$. Let $c\in B_n\subseteq Sa_i,\;b\in
Sa_{j_{n+1}}\cap Sa_i$. Since $c,b\in Sa_i$, on Proposition 2.5
either $Sc\subseteq Sb$ or $Sb\subseteq Sc$, that is either $c\in
Sa_{j_{n+1}}\cap B_n$ or $b\in Sa_{j_{n+1}}\cap B_n$.\hfill$\Box$

\smallskip {\bf Lemma 5.2}. Let ${\mathfrak R}$ be an axiomatizable class
and for any $a\in R$ and an idempotent $e\in R$ from the inclusion
$Sa\subseteq Se$ follows the existence an idempotent $f\in R$ such
that $Sa=Sf$ and $fS\subseteq eS$. Then for any idempotent $g\in
R$ there exists an idempotent $h\in R$ such that $Sh\subseteq Sg$,
the polygon ${}_SSh$ is minimum by inclusion and the right ideal
$hS$ is minimum for the principal right ideals of the monoid $S$,
generated by idempotents.

{\bf Proof}. Suppose the conditions of the lemma hold and $g\in
R,\;g=g^2$. Assume there exist the infinitely decreasing chain of
polygons $${}_SSg\supset {}_SSk_1\supset\ldots\supset
{}_SSk_n\supset\ldots,$$ where $k_i\in R\;(i\geqslant 1)$. On the
condition there exist the idempotents $h_i\in R\;(i\geqslant0)$
such that $h_0=g,\;Sk_i=Sh_i\;(i\geqslant 1)$ and
$$gS=h_0S\supseteq h_1S\supseteq\ldots\supseteq
h_nS\supseteq\ldots.$$ In view of the axiomatizability of the
class ${{\mathfrak R}}$ and Theorem 4.1 there exists $n\in\omega$
such that $h_iS=h_jS$ for all $i\geqslant n,\;j\geqslant n$. On
Proposition 1.2 the inclusion $Sh_{n+1}\subset Sh_n$ and the
equality $h_{n+1}S=h_nS$ imply $h_{n+1}=h_n$ which contradicts
with the suggestion $Sh_{n+1}\neq Sh_n$. Consequently, there
exists a minimum by inclusion poligon ${}_SSh$, such that
$Sh\subseteq Sg$. On Proposition 1.3 the right ideal $hS$ is
minimum for the principal right ideals of the monoid $S$,
generated by idempotents.\hfill$\Box$

\smallskip {\bf Theorem 5.1} [Ste1]. Let the class ${{\mathfrak R}}$ for regular polygons
be axiomatizable. The class ${{\mathfrak R}}$ is model complete if
and only if the following conditions hold:

1) $S$ is a regularly linearly ordered monoid;

2) for any idempotent $e\in R$ and for $a,a_i\in S,\;1\leqslant
i\leqslant m$ if $Sa\subset Se$ and
$e\not\in\cup\{a_iS\mid1\leqslant i\leqslant m\}$ then there exist
the idempotents $e_j\in R\;(j\in\omega) $ such that
$$e_j\neq e_k,\;Sa=Se_j,\;e_j\in
eS\setminus\cup\{a_iS\mid1\leqslant i\leqslant m\}$$ for any
$j,k\in \omega,\;j\neq k$;

3) $|eSf|\geqslant\omega$ for any idempotents $e,f\in R$.

{\bf Proof.  Necessity.} Suppose that the class ${{\mathfrak R} }$
is  model complete.

Let us prove condition 1. Assume $e^2=e\in R$. Since any cyclic
regular polygon is isomorphic to a subpolugon of the polygon
${}_SS$, generated by idempotent, it follows that it is enough to
show that the set $\{Sa\mid Sa\subseteq Se \}$ is linearly ordered
by inclusion. Let $Sa_1\subseteq Se$ and $Sa_2\subseteq Se$. Then
$a_1e=a_1,\;a_2e=a_2$ and ${}_SSe\models\exists x(a_1x=a_1\wedge
a_2x=a_2)$. In view of model completeness of the class
${{\mathfrak R}}$ we have ${}_S(Sa_1\cup Sa_2)\prec {}_SSe$, that
is
$${}_S(Sa_1\cup Sa_2)\models\exists x(a_1x=a_1\wedge a_2x=a_2).$$
Suppose $a_1c=a_1,\;a_2c=a_2$, where $c\in Sa_1\cup Sa_2$. If for
instance $c\in Sa_1$ then $a_2=a_2c\in Sa_1 $. Thus,
$Sa_2\subseteq Sa_1$.

Let us prove condition 2. Suppose $Sa\subset
Se,\;e\not\in\cup\{a_iS\mid1\leqslant i\leqslant m\}$, $e^2=e\in
R$. Assume
$$\Phi(y,a)\rightleftharpoons a=ay\wedge\bigwedge\limits_{i=1}^m
\neg\exists x(y=a_ix)\wedge y=ey.$$ Then ${}_SSe\models\Phi(e,a)$,
that is ${}_SSe\models\exists y\Phi(y,a)$. Since the class
${{\mathfrak R}}$ is model complete, we have ${}_SSa\prec {}_SSe$
and ${}_SSa\models\exists y\Phi(y,a)$. Consequently,
${}_SSa\models\Phi(e_1,a)$ for some $e_1\in Sa\subseteq R$. Hence
$a=ae_1$ and $Sa=Se_1$. Let $e_1=ka$. Since $e_1=ee_1$ it follows
that $e_1\in eS$. We claim that  $e_1$ is an idempotent. Let
$e_1e_1=kae_1=ka=e_1$. Furthermore,
${}_SSa\models\bigwedge\limits_{i=1}^m\neg\exists x(e_1=a_ix)$,
that is $e_1\not\in\cup\{a_iS\mid 1\leqslant i\leqslant m\}$.

Suppose there exist the idempotents $e_1,\ldots,e_k\in R$
satisfying the conditions: $$Sa=Se_i,\;e_i\in eS\setminus\cup\{a_r
S\mid 1\leqslant r\leqslant m\},\;e_i\neq e_j,$$ for any
$i,j,\;i\neq j,\;1\leqslant i,j\leqslant k$. We claim that there
exists an idempotent $e_{k+1}\in R$ satisfying the same conditions
with substitution $i$ on $k+1$, that is the idempotent $e_{k+1}\in
R$ such that ${}_SSe\models\Phi(e_{k+1})$, and $e_{k+1}$ is not
equal to $e_1,\ldots ,e_k$. Under condition $Se_j= Sa\subset Se$,
so $e\neq e_j\;(1\leqslant j\leqslant k)$. So
${}_SSe\models\exists y\Psi(y)$, where
$$\Psi(y)\rightleftharpoons\bigwedge\limits_{j=1}^k\neg
y=e_j\wedge\Phi(y).$$ Consequently, ${}_SSa\models\Psi(e_{k+1})$
for some $e_{k+1}\in Sa$. Thus, ${}_SSe\models\Phi(e_{k+1})$ and
$e_{k+1}$ is not equal to $e_1,\ldots ,e_k$. As for $e_1$ it is
proved that $e_{k+1}$ is an idempotent.

Let us prove condition 3. Suppose $e^2=e\in R$, $f^2=f\in R$,
${}_SSf_i\;(i\in\omega)$ are mutually disjoint copies of the
polygon ${}_SSf$. Since ${}_SSf\models e(ef)=ef$ it follows that
${}_SSf\models\exists x(ex=x)$ and
${}_SSf\sqcup\coprod\limits_{i\in\omega}{{}_S}Sf_i\models\exists
x_1,\ldots,x_n(\bigwedge\limits_{i\neq j}x_i\neq x_j\wedge
\bigwedge\limits_{i\leqslant n}ex_i=x_i)$ for any $n\geqslant1$.
In view of the model completeness of the class ${\mathfrak R}$ and
inclusion $Sf\subseteq Sf\sqcup\coprod\limits_{i\in\omega}Sf_i$ we
have ${}_SSf\models\exists x_1,\ldots,x_n(\bigwedge\limits_{i\neq
j,}x_i\neq x_j\wedge \bigwedge\limits_{i\leqslant n}ex_i=x_i)$ for
any $n\geqslant 1$. Consequently, $|eSf|\geqslant\omega$.

{\bf Sufficiency.} Let conditions 1--3 of the theorem hold.
Suppose
$${}_SA,{}_SB\in{{\mathfrak R}},\;{}_SB\subseteq{}_SA,\;\bar d=\langle
d_1,\ldots,d_{r}\rangle\in B^r,$$
$${}_SA\models\exists \bar x\bigwedge\limits_{j=1}^4\Phi_j(\bar x,\bar d),$$
where $\bar x=\langle x_1,\ldots x_{k}\rangle$,
$$\Phi_1(\bar x,\bar d)\rightleftharpoons\wedge\{nx_i=mx_j\mid\langle i,j,n,m\rangle\in
L_1\},$$
$$ \Phi_2(\bar x,\bar
d)\rightleftharpoons\wedge\{nx_i=md_j\mid\langle i,j,n,m\rangle\in
L_2\},$$
$$\Phi_3(\bar x,\bar d)\rightleftharpoons\wedge\{\neg nx_i=mx_j\mid\langle i,j,n,m\rangle\in
L_3\},$$ $$\Phi_4(\bar x,\bar d)\rightleftharpoons\wedge\{\neg
nx_i=md_j\mid\langle i,j,n,m\rangle\in L_4\},$$
$$L_t\subseteq \tilde{k}\times \tilde{k}\times S\times S\;(t\in\{1,3\}),\;
L_t\subseteq \tilde{k}\times \tilde{r}\times S\times
S\;(t\in\{2,4\}),\;\tilde{k}=\{1,\ldots,k\},\;\tilde{r}=\{1,\ldots,r\}.$$
 moreover, if $\langle i,j,n,m\rangle\in L_r$ then
 $\langle j,i,m,n\rangle\in L_r\;(r\in\{1,3\})$.
On the definition of the model completeness of the theory and on
Theorem 3.2 to proving ${}_SB\prec{}_SA$ it is enough to show
$${}_SB\models\exists \bar x\bigwedge\limits_{j=1}^4\Phi_j(\bar
x,\bar d).$$

Suppose ${}_SA\models \bigwedge\limits_{j=1}^4\Phi_j(\bar a,\bar
d),\;\bar a=\langle a_1,\ldots a_{k}\rangle\in A^k$. We can
consider that for any $i,j$, which are satisfying condition
$Sa_i\cap Sa_j\neq\emptyset$, there exist $n,m\in S$ such that
$\langle i,j,n,m\rangle\in L_1$. Since $a_i\in
A,\;{}_SA\in{\mathfrak R}$ and on Proposition 2.1 there exist an
idempotent $f_i\in R$ and an isomorphism
$\varphi_i:{}_SSa_i\longrightarrow {}_SSf_i\;(1\leqslant
i\leqslant k)$. If $Sa_i\subset Sa_j$ then on condition 2 we can
consider $Sf_i=S\varphi_i(a_i)\subseteq Sf_j$ and
$\varphi_i=\varphi_j|_{Sa_i}$. If $Sa_i=Sa_j$ then we suppose
$f_i=f_j$ and $\varphi_i=\varphi_j$.

For $i\;(1\leqslant i\leqslant k)$ such that $Sa_i\cap
B\neq\emptyset$ suppose
$$Sc_i=\max(\{Smd_j\mid md_j\in Sa_i,\;\langle j',j,m',m
\rangle\in L_2\;\mbox{for some}\;j',\;1\leqslant j'\leqslant
k,\;m'\in S\}\cup$$$$\cup\{Sma_j\mid ma_j\in B\cap Sa_i,\langle
j,j',m,m' \rangle\in L_1 \;\mbox{for some}\;j',\;1\leqslant
j'\leqslant k,\;m'\in S \}\cup\{Sb\}),$$ where $b\in
B\cap\bigcap\{Sa_j\mid Sa_j\cap Sa_i\neq\emptyset\}$. The
correctness of the definition of the set $Sc_i$ follows from Lemma
5.1 and Proposition 2.5.

Let us fix $i,j,\;1\leqslant i,j\leqslant k$. Assume $B\cap
Sa_i\cap Sa_j\neq\emptyset$, $$Sd_{ij}=\max (\{Smd_{j'}\mid
md_{j'}\in Sa_i\cap Sa_j,\;\langle i',j',n,m \rangle\in
L_2\;\mbox{for some}\;i',\;1\leqslant i'\leqslant k,\;n\in
S\}\cup$$$$\cup\{Sma_{i'}\mid ma_{i'}\in Sa_i\cap Sa_j,\;\langle
i',j',m,n \rangle\in L_1\; \mbox{for some}\;j',\;1\leqslant
j'\leqslant k,\;n\in S \}.$$ In view of Proposition 2.5 this
denotement is correct. Clearly, $Sd_{ij}=Sd_{ji}$. Suppose
$d_{ij}=d_{ji}$.

Renumber (if it is necessary) the elements of the set
$\{a_1,\ldots,a_k\}$ so as for any $i,\;1\leqslant i< k$, we have:
\begin{equation}\tag{5.1}
Sd_{{i+1}i}=\max\{Sd_{ji}\mid i<j\leqslant k\}.
\end{equation}

The proof of the sufficiency will contain of some lemmas.

\smallskip {\bf Lemma 5.3}. For any $i,j$, $1\leqslant i<k,\;0\leqslant j<k-i$,
we have

1) $Sd_{{i+j+1}\;i}\subseteq Sd_{i+j\;i}$;

2) $Sd_{{i+j+1}\;i}\subseteq\cap\{Sd_{{i+s+1}\;i+s}\mid0\leqslant
s\leqslant j\}$.

{\bf Proof}. Let $i,j$ be any numbers such that $1\leqslant
i<k,\;0\leqslant j<k-i$.

1) We will prove by the induction on $j$. If $j=0$ then the
inclusion $Sd_{i+1\;i}\subseteq Sd_{ii}$ follows from the
definition of $Sd_{ii}$. Assume the statement 1 of this lemma is
proved for all $j'<j$, that is
$$Sd_{{i+j}\;i}\subseteq Sd_{{i+j-1}\;i}\subseteq\ldots\subseteq
Sd_{{i+1}\;i}\subseteq Sd_{ii}.$$ We claim that
$Sd_{{i+j+1}\;i}\subseteq Sd_{{i+j}\;i}$.

Suppose the contraries. In view of the regular linear ordering of
monoid $S$ and on Proposition 2.5 we have $Sd_{{i+j}\;i}\subset
Sd_{{i+j+1}\;i}$. Let $a$ be an element of the form $md_{i'}$ or
$ma_{i'}$ from the definition $Sd_{{i+j+1}\;i}$ such that
$Sd_{{i+j+1}\;i}=Sa$ and $a\in Sa_{i+j+1}\cap Sa_i $. Then
$a\not\in Sa_{i+j}$. By the induction on $r,\;0\leqslant
r\leqslant j-1$, we will prove that $a\not\in Sa_{i+j-r}$. Suppose
$a\not\in Sa_{i+j-r+1}$, $a\in Sa_{i+j-r}$. Since $a\in
Sa_{i+j+1}$, on (5.1) and the definition of
$Sd_{{i+j-r}\;{i+j-r+1}}$ we derive $a\in
Sd_{{i+j-r}\;{i+j-r+1}}$. Consequently, $a\in Sa_{i+j-r+1}$, which
contradicts with the induction suggestion. Thus, $a\not\in
Sa_{i+j-r}$ for any $r,\;0\leqslant r\leqslant j-1$. In
particular, $a\not\in Sa_{i+1}$. Since $a\in Sa_i\cap Sa_{i+j+1}$
and on (5.1) we have $i=i+j+1$ that impossible.

2) On the statement 1 of this lemma $Sd_{{i+j+1}\;i}\subseteq
Sd_{{i+j}\;i}\subseteq\ldots\subseteq Sd_{{i+1}\;i}\subseteq
Sd_{ii}$. Suppose the inclusion $Sd_{{i+j+1}\;i}\subseteq
Sd_{{i+1}\;i}\cap\ldots\cap Sd_{{i+r}\;{i+r-1}}$ is proved for
$0\leqslant r\leqslant j-1,\;1\leqslant j<k $. We claim that
$Sd_{{i+j+1}\;i}\subseteq Sd_{{i+r+1}\;{i+r}}$. Let
$Sd_{{i+j+1}\;i}=Sb$, where $b$ is an element of the form
$md_{i'}$ or $ma_{i'}$ from the definition of $Sd_{{i+j+1}\;i}$,
and $b\in Sa_{i+j+1}\cap Sa_i$. Since $Sd_{{i+j+1}\;i}\subseteq
Sd_{{i+r}\;{i+r-1}}$ it follows than $b\in Sa_{i+r}$. The
inequality $i+j+1>i+r+1$ and definition of $Sd_{i+r\;{i+r+1}}$
imply $b\in Sd_{{i+r}\;{i+r+1}}$, that is
$Sd_{{i+j+1}\;i}\subseteq Sd_{i+r\;i+r+1}$.\hfill$\Box$

\smallskip {\bf Lemma 5.4}. For any $i,j,j',\;1\leqslant i,j,j'\leqslant k $,
there exists an idempotent $e_{ij}\in R$ such that
$Se_{ij}=S\varphi_i(d_{ij})$, moreover, if $Se_{ij}\subseteq
Se_{ij'}$ then $e_{ij}S\subseteq e_{ij'}S\subseteq f_iS$.

{\bf Proof}. On condition 2 of this theorem there exist the
idempotents $e'_{ij}\in R$ such that
$Se'_{ij}=S\varphi_i(d_{ij})$. Let $1\leqslant i\leqslant k$. In
view of the regular linear ordering of the monoid $S$ the set
$\{Se'_{ij}\;|\;1\leqslant j\leqslant k\}$ is linear ordering by
inclusion. Assume $Se'_{ij_1}\subseteq Se'_{ij_2}\subseteq
\ldots\subseteq Se'_{ij_k}\subseteq Sf_i$, where
$\{j_1,\ldots,j_k\}=\{1,\ldots,k\}$. Suppose there exist the
idempotents $e_{ij_t},e_{ij_{t+1}},\ldots,e_{ij_k}\in R$ such that
$Se_{ij_t}=Se'_{ij_t},\;Se_{ij_{t+1}}=Se'_{ij_{t+1}},\ldots,Se_{ij_k}=Se'_{ij_k}$,
$Se_{ij_t}\subseteq Se_{ij_{t+1}}\subseteq \ldots\subseteq
Se_{ij_k}\subseteq Sf_i$, where $1<t\leqslant k$. We claim that
there exists an idempotent $e_{ij_{t-1}}\in R$ such that
$Se_{ij_{t-1}}=Se'_{ij_{t-1}}$ and $e_{ij_{t-1}}S\subseteq
e_{ij_{t}}S$. If $Se'_{ij_{t-1}} = S e_{ij_{t}}$, we suppose
$e_{ij_{t-1}}=e_{ij_{t}}$. If $Se'_{ij_{t-1}}\subset Se_{ij_t}$
then $e_{ij_{t}}S\not\subseteq e'_{ij_{t-1}}S$ (otherwise, on
Proposition 1.2 $e_{ij_{t-1}}= e'_{ij_{t}}$) and on condition 2 of
this theorem there exists an idempotent $e_{ij_{t-1}}\in R$ such
that $Se_{ij_{t-1}}=Se'_{ij_{t-1}}$ and $e_{ij_{t-1}}S\subseteq
e_{ij_{t}}S$.\hfill$\Box$

Since $\varphi_i$ is an isomorphism and
$\varphi_i(d_{ij})=\varphi_i(d_{ij})e_{ij}\in
\varphi_i(d_{ij})e_{ij}Se_{ij}$ we derive
$d_{ij}\in\varphi_i(d_{ij})e_{ij}Sd_{ij}$. Since $\varphi_j$ is
the isomorphism we have
$\varphi_j(d_{ij})\in\varphi_i(d_{ij})e_{ij}S\varphi_j(d_{ij})=\varphi_i(d_{ij})e_{ij}Se_{ji}.$
Choose $t_{ij}\in e_{ij}Se_{ji}$ such that
$\varphi_j(d_{ji})=\varphi_i(d_{ij})t_{ij}$. If
$\varphi_i|_{Sd_{ij}}=\varphi_j|_{Sd_{ji}}$ then we put
$t_{ij}=e_{ij},\;t_{ji}=e_{ji}$, in particular, $t_{ii}=e_{ii}$.

\smallskip {\bf Lemma 5.5}. For any $i,j,\;1\leqslant i,j\leqslant k,$ the follows are hold:

1) $\varphi_j(x)=\varphi_i(x)t_{ij}$ for any $x\in Sd_{ij}$;

2) $t_{ij}\cdot t_{ji}=e_{ij}$;

3) $\varphi_i^{-1}(y)=\varphi_j^{-1}(yt_{ij})$ for any $y\in
Se_{ij}$.

{\bf Proof}. 1) Let $x\in Sd_{ji}$. Then $x=sd_{ji}=sd_{ij}$ for
some $s\in S$. Consequently,
$$\varphi_j(x)=\varphi_j(sd_{ji})=s\varphi_j(d_{ji})=s\varphi_i(d_{ij})t_{ij}=
\varphi_i(sd_{ij})t_{ij}=\varphi_i(x)t_{ij}.$$

2) Since $t_{ij}\in e_{ij}S$ we have $e_{ij}t_{ij}=t_{ij}$. Since
$Se_{ij}=S\varphi_i(d_{ij})$ we derive $e_{ij}=\varphi_i(x)$ for
some $x\in Sd_{ij}$. On the statement 1 of this lemma
$\varphi_j(x)=\varphi_i(x)t_{ij}$ and
$\varphi_i(x)=\varphi_j(x)t_{ji}$. Consequently, $$t_{ij}\cdot
t_{ji}=e_{ij}t_{ij}t_{ji}=\varphi_i(x)t_{ij}\cdot
t_{ji}=\varphi_j(x)t_{ji}=\varphi_i(x)=e_{ij}.$$

3) Let $y\in Se_{ij}=S\varphi_i(d_{ij})$. Then
$\varphi_i^{-1}(y)\in Sd_{ij}$. On the statement 1 of this lemma
$\varphi_j(\varphi_i^{-1}(y))=\varphi_i(\varphi_i^{-1}(y))t_{ij}=yt_{ij}$,
that is $\varphi_i^{-1}(y)=\varphi_j^{-1}(yt_{ij})$.\hfill$\Box$

\smallskip {\bf Lemma 5.6}. Let $Sa_i\cap B\neq\emptyset$ and
$c_i\in Sa_j$ for all $i,j,\;1\leqslant i,j\leqslant k$. Then
there exist the idempotents $g_1,\ldots,g_k\in R$ such that for
any $i,j,\;1\leqslant i,j\leqslant k$, the following conditions
hold:

1) $Sg_i=S\varphi_i(c)$;

2) $g_i\in Se_{ij}$;

3) $g_i\in e_{i\;i-1}S$;

4) $g_i=t_{i\;i-1}g_{i-1}t_{i-1\;i}$;

5) $e_{i-1\;i}g_{i-1}=t_{i-1\;i}g_it_{i\;i-1}$;

6) $e_{i1}g_i=t_{i1}g_1t_{1i}$;

7) $g_1=t_{1i}g_it_{i1}$.

{\bf Proof}. Suppose conditions of this lemma hold. Then
$Sc_i=Sc_j\;(1\leqslant i,j\leqslant k)$. Let $Sc=Sc_1$. So
$Sc\subseteq Sd_{ij}$ for all $i,j,\;1\leqslant i,j\leqslant k$.
Choose an idempotent $g_1$ so that $S\varphi_1(c)=Sg_1$ and
$g_1\in e_{1i}S\;(1\leqslant i\leqslant k)$. It is possible to do
it because condition 2 of this theorem hold and by Lemma 5.4 the
set $\{e_{1j}S\;|\;1\leqslant j\leqslant k\}$ is linear ordering
by inclusion. Then
$$g_1\in S\varphi_1(c)\subseteq S\varphi_1(d_{1i})=Se_{1i}\;(1\leqslant i\leqslant k)$$
and for $i=1$ conditions 1,2 of this lemma hold. For $i=1$
conditions 6,7 of this lemma hold trivially, since
$t_{11}=e_{11}$.

Suppose $i\geqslant 2$ and for the idempotents
$g_1,\ldots,g_{i-1}$ conditions 1-7 hold. We put
$g_i=t_{i\;i-1}g_{i-1}t_{i-1\;i}$ and prove that $g_i$ is an
idempotent. On Lemma 5.5 (2) $t_{i-1\;i}\cdot
t_{i\;i-1}=e_{i-1\;i}$. Since $g_{i-1}\in Se_{i-1\;i}$ we have
$g_{i-1}e_{i-1\;i}=g_{i-1}$. Consequently,
$$g_ig_i=t_{i\;i-1}g_{i-1}t_{i-1\;i}t_{i\;i-1}g_{i-1}t_{i-1\;i}
=t_{i\;i-1}g_{i-1}e_{i-1\;i}g_{i-1}t_{i-1\;i}=t_{i\;i-1}g_{i-1}t_{i-1\;i}=g_i.$$

We claim condition 1 of this lemma, that is $Sg_i=S\varphi_i(c)$.
Since $Sg_{i-1}=S\varphi_{i-1}(c)$ it follows that
$g_{i-1}=\varphi_{i-1}(x),\;\varphi_{i-1}(c)=rg_{i-1}$, where
$x\in Sc\subseteq Sd_{i\;i-1},\;r\in S $. On Lemma 5.5 (1) we have
$$g_i=t_{i\;i-1}g_{i-1}t_{i-1\;i}=t_{i\;i-1}\varphi_{i-1}(x)t_{i-1\;i}=t_{i\;i-1}\varphi_i(x)\in
S\varphi_i(c);$$
$$\varphi_i(c)=\varphi_{i-1}(c)t_{i-1\;i}=rg_{i-1}t_{i-1\;i}=rg_{i-1}g_{i-1}t_{i-1\;i}=$$$$
=r\varphi_{i-1}(x)g_{i-1}t_{i-1\;i}=r\varphi_i(x)t_{i\;i-1}g_{i-1}t_{i-1\;i}=r\varphi_i(x)g_i\in
Sg_i.$$

Since $g_i\in S\varphi_i(c)\subseteq S\varphi_i(d_{ij})=Se_{ij}$
for any $j,\;1\leqslant j\leqslant k$ condition 2 of this lemma
hold.

On the definition of $g_i$ and the building of $t_{i\;i-1}$ (on
the building we have $t_{i\;i-1}\in e_{i\;i-1}S$) condition 3
hold.

We claim condition 5 of this lemma. On Lemma 5.5 (2) and the
equality $g_{i-1}e_{i-1\;i}=g_{i-1}$, which is true on the
suggestion of the induction, we have
$$t_{i-1\;i}g_{i}t_{i\;i-1}=t_{i-1\;i}t_{i\;i-1}g_{i-1}t_{i-1\;i}t_{i\;i-1}
=e_{i-1\;i}g_{i-1}e_{i-1\;i}=e_{i-1\;i}g_{i-1}.$$

We claim condition 6 of this lemma. Since $e_{i1}\in
S\varphi_i(d_{i1})$ and $g_1\in S\varphi_1(c)$ it follows that
$e_{i1}=\varphi_i(y),\;g_1=\varphi_1(z)$, where $y\in
Sd_{i1},\;z\in Sc$. On Lemma 5.3 $Sd_{i1}\subseteq
\cap\{Sd_{2+r\;1+r}\mid 0\leqslant r\leqslant i-2\}$; on the
definition of $Sc$ we have $Sc\subseteq\cap\{Sd_{jr}\mid
1\leqslant j,r\leqslant k\}$. Consequently, on Lemmas 5.5 (1)
$$e_{i1}t_{i1}=\varphi_i(y)t_{i1}=\varphi_1(y)=\varphi_2(y)t_{21}=\ldots
=\varphi_i(y)t_{i\;i-1}\cdot\ldots \cdot
t_{21}=e_{i1}t_{i\;i-1}\cdot\ldots\cdot t_{21};$$
$$g_1t_{1i}=\varphi_1(z)t_{1i}=\varphi_i(z)=\varphi_{i-1}(z)t_{i-1\;i}=
\ldots=\varphi_1(z)t_{12}\cdot\ldots\cdot
t_{i-1\;i}=g_1t_{12}\cdot\ldots\cdot t_{i-1\;i}.$$ Since condition
4 of this lemma is true for all index, not greater than $i$ we
have
$$t_{i1}g_1t_{1i}=e_{i1}t_{i1}g_1t_{1i}=e_{i1}t_{i\;i-1}
\ldots t_{21}g_1t_{12}\ldots t_{i-1\;i}=e_{i1}g_i.$$

The follow equalities, which are true on Lemma 5.5 (2) and
condition 6 of this lemma, imply condition
7:$$g_1=e_{1i}g_1e_{1i}=e_{1i}t_{1i}t_{i1}g_1t_{1i}t_{i1}=e_{1i}t_{1i}e_{i1}g_it_{i1}.$$
Thus, we build the idempotents $g_1,\ldots,g_k\in R$ for which
conditions 1--7 of this lemma hold.\hfill$\Box$

\smallskip {\bf Remark 5.1}. The idempotent $g_1$ is chose arbitrarily with the regard for the conditions
\begin{equation}\tag{5.2}
S\varphi_1(c)=Sg_1\; \mbox{ш}\; g_1\in e_{1i}S\; \mbox {for all}
\;i,\;1\leqslant i\leqslant k.
\end{equation}

\emph{In lemmas 5.7--5.9 we suggest that the elements
$a_1,\ldots,a_k$ satisfy the conditions: $Sa_i\cap B\neq\emptyset$
and $Sa_i\cap Sa_j\neq \emptyset$ for all $i,j\;1\leqslant
i,j\leqslant k$.} On the set $\{a_1,\ldots,a_k\}$ we define the
following binary relation:
$$a_i\sim
a_j\Longleftrightarrow c_i,c_j\in Sa_i\cap Sa_j.$$ On the
definition of the sets $Sc_i\;(1\leqslant i\leqslant k)$ this
relation is an equivalence. For each class of this equivalence we
build the idempotent $g_i$ satisfying Lemma 5.6.

\smallskip {\bf Lemma 5.7}. ${}_SB\models\Phi_1(\bar b,\bar d)$, where $\bar
b=\langle b_1,\ldots,b_k\rangle\in B^k$,
$b_i=\varphi_i^{-1}(\varphi_i(a_i)g_i)\;(1\leqslant i\leqslant
k)$.

{\bf Proof}. Let $\langle i,j,n,m\rangle\in L_1$. Suppose $a_i\sim
a_j$.  We claim that ${}_SB\models nb_i=mb_j$ for
$b_r=\varphi_r^{-1}(\varphi_r(a_r)g_r),\;r\in\{i,j\}$. Let
$i=j+l$, where $l\geqslant0$. Since $na_i=ma_j\in Sa_i\cap Sa_j$
we have $na_i\in Sd_{ij}=Sd_{j+l\;j}$. On Lemma 5.3
$$Sd_{j+l\;j}\subseteq \bigcap\{Sd_{j+r+1\;j+r}\mid0\leqslant r\leqslant l-1\}$$
for $l\geqslant 1$. Consequently, on Lemma 5.5 (1,3) and Lemma 5.6
(4) we have
$$mb_j=m\varphi_j^{-1}(\varphi_j(a_j)g_j)=\varphi_j^{-1}(\varphi_j(ma_j)g_j)=\varphi_{j+1}^{-1}(\varphi_{j+1}(ma_j)t_{i+1\;j}g_jt_{j\;j+1})=
$$ $$=\varphi_{j+1}^{-1}(\varphi_{j+1}(ma_j)g_{j+1})=\ldots=\varphi_{j+l}^{-1}(\varphi_{j+l}(ma_j)g_{j+l})=
\varphi_i^{-1}(\varphi_i(na_i)g_i)=n\varphi_i^{-1}(\varphi_i(a_i)g_i)=nb_i,$$
that is $mb_j=nb_i$. If $l=0$ then the equality $mb_j=nb_i$ is
obviously.

Suppose $a_i\sim a_j$ is false. Let, for instance, $c_i\not\in
Sa_j$. We claim that $Sa_i\cap Sa_j\subseteq B$. Let $a\in
Sa_i\cap Sa_j$. On the definition of $Sc_i$ we have $Sc_i\subseteq
Sa_i\cap B$. On Proposition 2.5 either $Sa\subseteq Sc_i$ or
$Sc_i\subseteq Sa$. If $Sc_i\subseteq Sa$, then $c_i\in Sa_j$
which is contradicts with the suggestion. Consequently, $Sa_i\cap
Sa_j\subseteq Sc_i\subseteq B$.

Since $na_i=ma_j$ we have $na_i,ma_j\in Sa_i\cap Sa_j\subseteq
Sc_i$. Consequently, $ma_j\in B$, $ma_j\in Sc_j$ and $na_i\in
Sc_i$. Thus,
$$\varphi_i(na_i)\in S\varphi_i(c_i)=Sg_i,\;\varphi_j(ma_j)\in
S\varphi_j(c_j)=Sg_j,$$$$nb_i=n\varphi_i^{-1}(\varphi_i(a_i)g_i)=
\varphi_i^{-1}(\varphi_i(na_i)g_i)=\varphi_i^{-1}(\varphi_i(na_i))=na_i.$$
Similarly, $mb_j=ma_j$. Hence $nb_i=na_i=ma_j=mb_j$, and lemma is
proved.\hfill$\Box$

\smallskip {\bf Lemma 5.8} ${}_SB\models\Phi_2(\bar b,\bar d)$, where $\bar
b=\langle b_1,\ldots,b_k\rangle\in B^k$,
$b_i=\varphi_i^{-1}(\varphi_i(a_i)g_i)\;(1\leqslant i\leqslant
k)$.

{\bf Proof}. Let $\langle i,j,n,m\rangle\in L_2$. We claim that
${}_SB\models nb_i=md_j$ for
$b_i=\varphi_i^{-1}(\varphi_i(a_i)g_i)$. On the definition of
$Sc_i$ we have $md_j\in Sc_i$. Since $na_i=md_j$ it follows that
$na_i\in Sc_i$. On Lemma 5.6 (1) $S\varphi_i(c_i)=Sg_i$.
Consequently, $\varphi_i(na_i)g_i=\varphi_i(na_i)$ and
$$nb_i=n\varphi_i^{-1}(\varphi_i(a_i)g_i)=\varphi_i^{-1}(\varphi_i(na_i)g_i)
=\varphi_i^{-1}(\varphi_i(na_i))=na_i=md_j,$$, that is
$nb_i=md_j$. \hfill$\Box$

On Remark 5.1 for each class of the equivalence by the relation
$\sim$ there exists some freedom of the choice of one of the
idempotents $g_i$, which were build for this class. Let
$g_{i_1},\ldots,g_{i_s}$ be all such idempotents.

\smallskip {\bf Lemma 5.9}. The idempotents $g_{i_1},\ldots,g_{i_s}$
can be choose such that for them it is held  condition (5.2) with
the change $g_1$ on $g_{i_t}\;(1\leqslant t\leqslant s)$ and
${}_SB\models\bigwedge\limits_{i=1}^4\Phi_i(\bar b,\bar d)$, where
$\bar b=\langle b_1,\ldots,b_k\rangle\in B^k$,
$b_i=\varphi_i^{-1}(\varphi_i(a_i)g_i)\;(1\leqslant i\leqslant
k)$.

{\bf Proof}. Let $\langle i,j,n,m\rangle\in L_3$. Note that
$\varphi_i(b_i)=\varphi_i(\varphi_i^{-1}(\varphi_i(a_i)g_i))=\varphi_i(a_i)g_i\in
Sg_i$, that is $b_i\in Sc_i\subseteq B$. Furthermore, $b_i\in
Sa_i$. Add an element $b_i$ to the set $\{d_1,\ldots,d_r\}$. We
can consider that $b_i=d_{i}'\in\{d_1,\ldots,d_r\}$. Add an
element $b_i$ to the set $\{a_1,\ldots,a_r\}$. We can consider
that $b_i=a_{i'}\in \{a_1,\ldots,a_k\}$. Similarly, we can
consider that $b_j=d_{j}'\in\{d_1,\ldots,d_r\}$ and $b_j=a_{j'}\in
\{a_1,\ldots,a_k\}$. Denote the formula $\Phi_1(\bar x,\bar
d)\wedge nx_{i'}=nd_{i}'\wedge mx_{j'}=md_{j}'$ by $\Phi_1(\bar
x)$  once again. Since $b_i\in Sc_i$, the sets $Sc_i$ and
$Sd_{il}$ are not changed, where $a_{l}\sim a_i$. On the same
reason the sets $Sc_j$ and $Sd_{jl}$ are not changed too, where
$a_{l}\sim a_j$. Since $a_{i'}\in Sa_i$, it follows that
$\varphi_{i'}=\varphi_i|Sa_{i'}$, $\varphi_{i'}(a_{i'})\in
Sg_{i'}'$ and
$b'_{i'}=\varphi^{-1}_{i'}(\varphi_{i'}(a_{i'})g'_{i'})=a_{i'}=b_i$
under any choice of idempotents $g_i'$ satisfying Lemma 5.6.
Similarly, $b'_{j'}=b_j$. If $na_i=md'_j$ and $ma_j=nd'_i$ then
$nb'_{i}=n\varphi^{-1}_{i}(\varphi_{i}(a_{i})g'_{i})=a_{i'}=\varphi^{-1}_{i}(\varphi_{i}
(a_{i})g'_{i})=na_i$ under any choice of idempotents $g_i'$
satisfying Lemma 5.6; similarly, $mb_j'=ma_j$, that is $nb'_i\neq
mb'_j$. So we will consider only the case: $na_i\neq md_{j}'$ or
$ma_{j}\neq nd_{i}'$. Without the restriction of generality we can
consider $ma_{j}\neq nd_{i}'$. Denote the formula $\Phi_4(\bar
x,\bar d)\wedge \neg mx_{j}= nd_{i}'$ by $\Phi_4(\bar x,\bar d)$
once again. If $na_i\neq md_{j}'$ or $ma_{j}\neq nd_{i}'$ denote
the formula, getting from the formula $\Phi_3(\bar x,\bar d)$ by
cancellation of subformulae $\neg mx_j=nx_i$ and $\neg nx_i=mx_j$,
by $\Phi_3(\bar x,\bar d)$. The same way we do with all $\langle
i,j,n,m\rangle\in L_3$. As a result, in the set $L_3$ it is
remained only such collections of $\langle i,j,n,m\rangle$, for
which $nb'_i\neq mb'_j$ under any choice of the idempotents $g_i'$
satisfying Lemma 5.6. If we will choose the idempotents
$g_{i_1},\ldots,g_{i_s}$ (under which it is built all another
$g_i\;,1\leqslant i\leqslant k $) satisfying (5.2) with
substitution $g_1$ on $g_{i_t}$ $(1\leqslant t\leqslant s)$ such
that ${}_SB\models\Phi_4(\bar x,\bar d)$ it follows that as it was
note above and on Lemmas 5.7 and 5.8,
${}_SB\models\bigwedge\limits_{i=1}^4\Phi_i(\bar b,\bar d)$, where
$\bar b=\langle b_1,\ldots,b_k\rangle\in B^k$.

Without loss of generality we can consider that $g_{i_1}=g_1$.
Show, how can be choose the element $g_1$ (another idempotents
$g_{i_2},\ldots,g_{i_s}$ are chosen similarly). Let
$K_1=\{1,2,\ldots,k_1\}=\{i\mid a_i\sim a_1\}$ and Lemma 5.3 holds
with substitution $k$ on $k_1$. Note that $\langle
i,j,n,m\rangle\in L_4$, $i\in K_1$ and $md_j\not\in Sa_i$ imply
$nb_i\neq md_j$ under any choice $g_1$, satisfying conditions
(5.2). For all $i\in K_1$ and $j\in K_1$ we will build $e'_{ij}$
and $t'_{ij}$ satisfying Lemma 5.4 with substitution $e_{ij}$ on
$e_{ij}'$, Lemma 5.5 with substitution $t_{ij}$ on $t_{ij}'$ and
the condition
$$\varphi_l(na_l)t'_{ll-1}\cdot t'_{l-1l-2}\cdot\ldots\cdot t'_{i+1i}e'_{ij}\neq
\varphi_l(md_t)t'_{ll-1}\cdot t'_{l-1l-2}\cdot\ldots\cdot
t'_{i+1i}e'_{ij},$$ where $md_t\in Sa_l,\;l\geqslant i,\; \langle
l,t,n,m\rangle\in L_4$. Suppose the elements
$e_{k_1j}',e_{k_{1-1}j}',\ldots,e_{i+1j}'$ and
$t_{k_1j}',t_{k_{1-1}j}',\ldots,t_{i+1j}'$ are built for all $j\in
K_1$. Let
$$Y_{i}^l=\{x\in
R\mid\varphi_l(na_l)t'_{ll-1}\cdot t'_{l-1l-2}\cdot\ldots\cdot
t'_{i+1i}x= \varphi_l(md_t)t'_{ll-1}\cdot
t'_{l-1l-2}\cdot\ldots\cdot t'_{i+1i}x,\;md_t\in Sa_l,\; \langle
l,t,n,m\rangle\in L_4 \},$$
$$Y_i=\cup\{Y_{i}^l\mid l\geqslant i\}.$$ Let us prove that $f_i\not\in Y_i $.
Let $f_i\in Y_i^l $. If $l=i$ then
$\varphi_i(na_i)=\varphi_i(na_i)f_i=\varphi_i(ma_j)f_i=\varphi_i(ma_j)=$,
that is $na_i=ma_j$, a contradiction. On the definition,
$t'_{i+1i}\in Se_{ii+1}'\subseteq Sf_i$, that is
$t'_{i+1i}f_i=t'_{i+1i}$. Consequently,
$\varphi_l(na_l)t'_{ll-1}\cdot t'_{l-1l-2}\cdot\ldots\cdot
t'_{i+1i}= \varphi_l(md_t)t'_{ll-1}\cdot
t'_{l-1l-2}\cdot\ldots\cdot t'_{i+1i}$. Then
$\varphi_l(na_l)t'_{ll-1}\cdot t'_{l-1l-2}\cdot\ldots\cdot
t'_{i+1i}\cdot t'_{ii+1}= \varphi_l(md_t)t'_{ll-1}\cdot
t'_{l-1l-2}\cdot\ldots\cdot t'_{i+1i}\cdot t'_{ii+1}$. On Lemma
5.5 (2) $\varphi_l(na_l)t'_{ll-1}\cdot t'_{l-1l-2}\cdot\ldots\cdot
t'_{i+2i+1}\cdot e'_{i+1i}= \varphi_l(md_t)t'_{ll-1}\cdot
t'_{l-1l-2}\cdot\ldots\cdot t'_{i+2i+1}\cdot e'_{i+1i}$ which
contradicts with the choice of the element $e'_{i+1i}$.  Let
$Y_i\neq\emptyset$. In view of the axiomatizability of the class
${\mathfrak R}$ we have $Y_i=\cup\{k_sR\;|\;1\leqslant s\leqslant
t\}$ for some $t> 0,\;k_s\in S \;(1\leqslant s\leqslant t)$. We
put $X_i=\cup\{k_sS\;|\;1\leqslant s\leqslant t\}$. Clearly,
$$h\in Y_i\Longleftrightarrow h\in X_i
$$ for any idempotent $h\in R$.
Consequently, $f_i\not\in X_{i}$. On Lemmas 5.3 and 5.4 the set
$\{Se_{ij}\;|\;j\in K_1\}$ is linearly ordered by inclusion. Let
$$
Se_{ij_1}\subseteq Se_{ij_2}\subseteq\ldots\subseteq
Se_{ij_r}\subseteq Sf_{i},
$$ where $\{j_1,\ldots,j_r\}=\{j\mid j\in K_1\}$.
Let $e'_{ij_{r+1}}=f_i$. Suppose the idempotents
$e_{ij_l}',e_{ij_{l+1}}',\ldots,e_{ij_r}'$ are already built,
moreover $Se_{ij_l}'=Se_{ij_l}$. Since $e'_{ij_l}\not\in X_i$, by
condition 2 of this theorem there exists an idempotent
$e_{ij_{l-1}}'$ such that
$Se_{ij_{l-1}}'=Se_{ij_{l-1}},\;e_{ij_{l-1}}'S\subseteq
e_{ij_{l}}S$ and $e_{ij_{l-1}}'\not\in X_i$. By the idempotents
$e_{ij}'$ we build $t'_{ij'}\in S$, satisfying conditions of Lemma
5.5. On the building, $g_1\in Se_{1j}=Se_{1j}'$, where $j\in K_1$.
Then on condition 2 of this theorem there exists an idempotent
$g_1'$ such that $Sg_1'=Sg_1,\;g_1'S\subseteq e_{1j}S$ for all
$j\in K_1$ and $g_1'\not\in X_1$. Similarly we build the
idempotents $g_{i_2},\ldots,g_{i_s}$. Another idempotents $g_i$ we
build such that they satisfy Lemma 5.6. The elements $b_i'$ are
defined as previously: $b_i'=\varphi_i^{-1}(\varphi_i(a_i)g_i')$,
where $i\in K_1$. Then, as it is noticed above,
${}_SB\models\bigwedge\limits_{i=1}^2\Phi_i(b_1',\ldots,b_k')$.

Let us prove that $nb'_i\neq md_j$, where $\langle
i,j,n,m\rangle\in L_4,\;i\in K_1$. Let $nb'_i= md_j$. Then
$\varphi_i(na_i)g'_i=\varphi_i (md_j)$. On Lemma 5.6 (4),
$\varphi_i(na_i)t'_{ii-1}\cdot g'_{i-1}\cdot t'_{i-1i}=\varphi_i
(md_j)$. Then $\varphi_i(na_i)t'_{ii-1}\cdot g'_{i-1}\cdot
t'_{i-1i}\cdot t'_{ii-1}=\varphi_i (md_j)\cdot t'_{ii-1}$. On
Lemma 5.5 (2), $\varphi_i(na_i)t'_{ii-1}\cdot g'_{i-1}\cdot
e'_{i-1i}=\varphi_i (md_j)\cdot t'_{ii-1}$. On Lemma 5.6 (2),
$\varphi_i(na_i)t'_{ii-1}\cdot g'_{i-1}=\varphi_i (md_j)\cdot
t'_{ii-1}$. Continuing this process, we have
$\varphi_i(na_i)t'_{ii-1}\cdot t'_{i-1i-2}\cdot\ldots\cdot
t'_{21}\cdot g'_1= \varphi_i(md_j)t'_{ii-1}\cdot
t'_{i-1i-2}\cdot\ldots\cdot t'_{21}$, that is $g'_1\in Y_1$.
Consequently, $g_1'\in X_1$ which contradicts with the building of
the element $g_1'$. Thus, ${}_SB\models\Phi_4(\bar b,\bar d)$.
\hfill$\Box$

Let $\{a_{i_1},\ldots,a_{i_s}\}$ be an arbitrary maximum subset of
the set $\{a_1,\ldots,a_k\}$ such that
\begin{equation}\tag{5.3}
Sa_i\cap B=\emptyset,\;Sa_i\cap Sa_j\neq\emptyset\;\mbox{for
any}\;i,j\in\{i_1,\ldots,i_s\},
\end{equation}
$b^1$ be an arbitrary element of $B$. On the definition of the
regular polygon there exist an idempotent $f^1\in R$ and an
isomorphism $\psi_1:{}_SSb^1\longrightarrow{}_SSf^1$. On Lemma 5.2
there is an idempotent $f\in R$ such that $Sf\subseteq Sf^1$ and
${}_SSf$ is a minimum by inclusion polygon. We put:
$Sb=S\psi_1^{-1}(f_1),\;\psi=\psi_1|_{Sb}$. Since $Sa_i\cap
Sa_j\neq\emptyset$ for any $i,j\in\{i_1,\ldots,i_s\}$, on Lemma
5.1 (1) there exists $a^1\in S$ such that $a^1\in\cap\{Sa_i\mid
i\in\{i_1,\ldots,i_s\}\}$. In view of the regularity of the
polygon ${}_SA$ there exist an idempotent $e^1\in R$ and an
isomorphism $\varphi:{}_SSa^1\longrightarrow{}_SSe^1$. On Lemma
5.2 there exists an idempotent $e\in R$ such that $Se\subseteq
Se^1,$ the polygon ${}_SSe$ is minimum and the right ideal $eS$ is
minimum for the principal right ideals of the monoid $S$,
generated by idempotents. Suppose $B'=Sa=S\varphi^{-1}(e)$, where
$a\in A$. Then $Sa_i\cap B'\neq\emptyset$ for any
$i\in\{i_1,\ldots,i_s\}$. Replacing $B$ by $B'$ in all previous
argumentations, which concern to the proof of sufficiency, we have
in particular $b^i=\varphi_i^{-1}(\varphi_i(a_i)g_i)\in B'$, where
$i\in\{i_1,\ldots,i_s\}$. Let $\alpha$ be an arbitrary isomorphism
from ${}_SSa$ to ${}_SSe$ (for instance, $\alpha=\varphi|_{Sa}$).
Denote the element $\psi^{-1}(\alpha(b^i)f)$ by $b_i$. Clearly,
$b_i\in Sb\subseteq B\;(i\in\{i_1,\ldots,i_s\})$. For the element
$a_i\in\{a_1,\ldots,a_k\}$ such that $Sa_i\cap B\neq\emptyset$, we
build the element $b_i\in B$ as in Lemma 5.9.

Let us prove that ${}_SB\models\bigwedge\limits_{i=1}^4\Phi_i(\bar
b,\bar d)$, where $\bar b=\langle b_1,\ldots,b_k\rangle\in B^k$.

If $\langle i,j,n,m\rangle\in L_2$ then $Sa_i\cap B\neq\emptyset$
and on Lemma 5.9 ${}_SB\models nb_i=md_j$.

If $\langle i,j,n,m\rangle\in L_1$ and $Sa_i\cap B\neq\emptyset$
then $Sa_i\cap Sa_j\neq\emptyset$. On condition 1 of this theorem
$Sa_j\cap B\neq\emptyset$; on Lemma 5.9 ${}_SB\models nb_i=md_j$.

Suppose $\langle i,j,n,m\rangle\in L_1$, $Sa_i\cap B=\emptyset$,
$Sa_j\cap B=\emptyset$. Then on Lemma 5.9 and choice $b^i,b^j\in
B'$ we have ${}_SSa\models nb^i=mb^j$. Consequently,
$$nb_i=n\psi^{-1}(\alpha(b^i)f)=\psi^{-1}(\alpha(nb^i)f)=
\psi^{-1}(\alpha(mb^j)f)=m\psi^{-1}(\alpha(b^j)f)=mb_j.$$

If $\langle i,j,n,m\rangle\in L_3$ or $\langle i,j',n,m\rangle\in
L_4$ and $Sa_i\cap Sa_j\neq\emptyset$, $Sa_i\cap B\neq\emptyset$,
$Sa_j\cap B\neq\emptyset$ then on Lemma 5.9 $ nb_i\neq mb_j$ or $
nb_i\neq md_{j'}$ accordingly.

Let $i,\;1\leqslant i\leqslant k$, be such that $Sa_i\cap
B=\emptyset$; the set $\{a_{i_1},\ldots,a_{i_s}\}$ be a maximum
subset of the set $\{a_1,\ldots,a_k\}$ with properties (5.3),
where $i_1=i$. Above for the set $\{a_{i_1},\ldots,a_{i_s}\}$ it
is defined the elements $a\in\cap\{Sa_j\mid
j\in\{i_1,\ldots,i_s\}\},\;e\in S$ and it is chose an arbitrary
isomorphism $\alpha:{}_SSa\longrightarrow {}_SSe$, moreover the
polygon ${}_SSe$ is minimum and the right ideal $eS$ is minimum
for the principal right ideals of the monoid $S$ generated by
idempotents.

Suppose $\langle i,j,n,m\rangle\in L_3$, $j\in\{i_1,\ldots,i_s\}$.
On Lemma 5.9 $ nb^i\neq mb^j$. Since $\alpha(nb^i),\alpha(mb^j)\in
Se$ we have $e\not\in X=\{x\mid\alpha(nb^i)x=\alpha(mb^j)x\}$.
Clearly, $efS\subseteq eS $. Since $Sef\subseteq Sf$, on condition
2 of this theorem, $efS=gS$ for some idempotent $g\in S$. In view
of the minimality of the right ideal $eS$ and the inclusion
$gS\subseteq eS$ we have the equality $eS=gS=efS$, that is $e=efr$
for a some $r\in S$. If $f\in X$ then $ef\in X$ and $e\in X$.
Consequently, $f\not\in X$ and
$$nb_i=n\psi^{-1}(\alpha(b^i)f)= \psi^{-1}(\alpha(nb^i)f)\neq
\psi^{-1}(\alpha(mb^j)f)=m\psi^{-1}(\alpha(b^j)f)=mb_j.$$

Suppose $\langle i,j,n,m\rangle\in L_3$ and
$j\not\in\{i_1,\ldots,i_s\}$ or $\langle i,j,n,m\rangle\in L_4$.
Define a homomorphism $\alpha_t:{}_SSe\longrightarrow{}_SSe\;(t\in
S)$ as follows: $\alpha_t(se)=sete$ for any $s\in S$. In view of
the minimality of the polygon ${}_SSe$ the equality $Se=Sete$
hold, that is $e=lete$ for some $l\in S$ and $\alpha_t$ is an
epimorphism. The minimality of the right ideal $eS$ implies the
equality $eS=eteS$, that is $e=etel'$ for some $l'\in S$ and
$\alpha_t$ is an isomorphism. Note that $\alpha_t\neq \alpha_s$
implies $\alpha_t(c)\neq \alpha_s(c)$ for any $c\in Se$. Really,
let $cete= \alpha_t(c)= \alpha_s(c)=cese$. Since $Se=Sce$, that is
$e=c_1ce$ for some $c_1\in S$, we have $ete=c_1cete=c_1cese=ese$,
which contradicts with suggestion. On condition 3 of this theorem,
$|eSe|\geqslant\omega$. Consequently, there exists
$t_i\;(i\in\omega)$ such that
$\alpha_{t_i}\neq\alpha_{t_j}\;(i\neq j)$. In view of the
minimality of the polygon $Sf$ the equality $Sef=Sf$ hold, that is
\begin{equation}\tag{5.4}
f=uef
\end{equation}
for some $u\in S$. For $i\in\{i_1,\ldots,i_s\}$ and $n\in S$ we
put: $$X_i^n=\{x\in Se\mid x=\psi(mb_j)ue,\;\mbox{where}\;mb_j\in
Sb,\;\langle i,j,n,m\rangle\in L_3\},$$
$$Y_i^n=\{x\in Se\mid x=\psi(md_j)ue,\;\mbox{where}\;md_j\in Sb,\;
\langle i,j,n,m\rangle\in L_4 \}.$$ If $i\in\{i_1,\ldots,i_s\}$
and $n\in S$ such that either the set $X_i^n$ or the set $Y_i^n$
is not determined, then we assume either $X_i^n=\emptyset$ or
$Y_i^n=\emptyset$ accordingly. Since the sets $X_i^n$ and
$Y_i^{n}$ are finite for all $i\in\{i_1,\ldots,i_s\}$ and $n\in
S$, but the set of the different isomorphisms $\alpha_t$ is
infinite, then there exists $t\in S$ such that
$\alpha_t\alpha(nb^i)\not\in X_i^n\cup Y_i^{n}$ for any
$i\in\{i_1,\ldots,i_s\},\;n\in S$. For the isomorphism $\alpha$ we
will take the isomorphism
$\alpha_t\alpha:{}_SSa\longrightarrow{}_SSe$, which we denote by
$\alpha$. Then $\alpha(nb^i)\not\in X_i^n\cup Y_i^{n}$ for any
$i\in\{i_1,\ldots,i_s\},\;n\in S$.

Suppose $\langle i,j,n,m\rangle\in L_3$ and
$j\not\in\{i_1,\ldots,i_s\}$. If $nb_i=mb_j$ then
$\psi^{-1}(\alpha(nb^i)f)=mb_j$ by the definitions of $b_i$ and
$mb_j\in Sb$, that is
$$\alpha(nb^i)ef=\alpha(nb^i)f=\psi(mb_j)=\psi(mb_j)f.$$
The equality (5.4) implies $\alpha(nb^i)ef=\psi(mb_j)uef$. In view
of the minimality of the right ideal $eS$ and the equality
$efS=eS$ we have $\alpha(nb^i)e=\alpha(nb^i)=\psi(mb_j)ue$, that
is $\alpha(nb^i)\in X_i^n$. We have a contradict with the choice
of the isomorphism $\alpha$. Consequently, $nb_i\neq mb_j$.

If $\langle i,j,n,m\rangle\in L_4$ then the inequality $nb_i\neq
md_j$ is proved similarly. Theorem is proved.\hfill$\Box$

\smallskip {\bf Corollary 5.1}. The class ${\mathfrak R}$ for the regular
polygons over the infinite group is axiomatizable and model
complete.

{\bf Proof} follows from Theorem 5.1 and Corollary
4.3.\hfill$\Box$

\vskip 1cm
\medskip {\bf\centerline {\S~6. Completeness of Class for Regular
Polygons}}
\medskip

In this paragraph we investigate the monoids with the complete
class of regular polygons. It is known that the model completeness
of the class $\Re$ for regular polygons implies the completeness
of this class (Lemma 6.1). Theorems 6.1 and 6.2 assert that the
completeness of the class $\Re$ implies the model completeness of
this class if we put some conditions on the monoid: the class
$\Re$ satisfies the condition of formula definability of
isomorphic orbits or the monoid is linearly ordered and has depth
2. Theorems 6.3 and 6.4 assert that this conditions in Theorems
6.1 and 6.2 are essential.

If for element $a\in {}_SA$ and $b\in {}_SB$ there exists an
isomorphism $f: {}_SSa\to {}_SSb$ such that $f(a)=b$, then this
fact will be denoted by ${}_SSa \mathrel{\widetilde{\to}} {}_SSb$.

We will say that a class ${}\Re$ satisfies the condition of
formula definability of isomorphic orbits if for each idempotent
$e\in R$ there exists a formula $\Phi_e(x)$ such that for any
regular polygon ${}_SA$ and any $a\in {}_SA$ the following is
true:
$${}_SA\models\Phi_e(a)\iff {}_SSa \mathrel{\widetilde{\to}} {}_SSe.$$

{\bf Theorem 6.1} [Ov2]. The following conditions are equivalent:

1) the axiomatizable class ${}\Re$ is model complete and satisfies
the condition of formula definability of isomorphic orbits;

2) the axiomatizable class ${}\Re$ is complete and satisfies the
condition of formula definability of isomorphic orbits;

3) the semigroup $R$ is a rectangular band of infinite groups, and
the set $I(R)$ is finite certainly.

{\bf Proof} contains some lemmas.

{\bf Lemma 6.1.} If the class of regular polygons $\Re$ is model
complete, then $\Re$ is complete.

{\bf Proof}. If ${}_SA,{}_SB\in \Re$, then the polygon
${}_SC\rightleftharpoons {}_SA\sqcup {}_SB$ is also regular. In
view of model completeness of the class ${}\Re$ the polygons
${}_SA$ and ${}_SB$ are elementary submodels of the polygon
${}_SC$. Then ${\rm Th}({}_SA)={\rm Th}({}_SC)={\rm Th}({}_SB)$,
that is ${}_SA\equiv {}_SB$.\hfill $\Box$

{\bf Lemma 6.2.} If the class $\Re$ is complete and satisfies the
condition of formula definability of isomorphic orbits then for
any idempotent $e\in R$ and any element $a\in R$ the condition
$Sa\subseteq Se$ implies the equality $Sa=Se$.

{\bf Proof}. Let $e$ be an idempotent belonging to $R$, and
formula $\Phi_e(x)$ be such that for any ${}_SA\in \Re$ and any
$a\in A$ the following is true:
$${}_SA\models\Phi_e(a)\iff {}_SSa \mathrel{\widetilde{\to}} {}_SSe.$$
Assume we have $Sa\subset Se$ for some element $a\in R$. Since
$a\in T$, the polygon ${}_SSa$ is regular. In view of completeness
of the class $\Re$ we get ${}_SSa\equiv{}_SSe$. Since
${}_SSe\models\Phi_e(e)$ then ${}_SSe\models\exists x\:\Phi_e(x)$
and, consequently, ${}_SSe\models\exists x\:\Phi_e(x)$. Let $b$ be
a realization of the formula $\Phi_e(x)$ in the polygon ${}_SSa$.
Since $Sa\subset Se$, then $b\ne e$ and $Sb\subset Se$. Then
$be=b$ and
$${}_SSe\models\Phi_e(b)\wedge be=b\wedge\Phi_e(e).$$
Consequently,
$${}_SSe\models\Phi_e(b)\wedge\exists y\:(by=b\wedge\Phi_e(y)).$$
Denote by $\Psi(x)$ the formula $\Phi_e(x)\to\exists
y\:(by=x\wedge\Phi_e(y))$. Obviously, ${}_SSe\models\Psi(b)$.

We will show that ${}_SSe\not\models\Psi(e)$. Since
${}_SSe\models\Phi_e(e)$, it is enough to check that
$${}_SSe\not\models\exists y\:(by=e\wedge\Phi_e(y))$$
or
$${}_SSe\models\forall y\:(\Phi_e(y)\to by\ne e).$$
Assume there exists an element $s\in Se$ such that
${}_SSe\models\Phi_e(c)$ and $b\cdot c=e$. Since
${}_SSe\models\Phi_e(c)$, then ${}_SSc \mathrel{\widetilde{\to}}
{}_SSe$. Let $f:{}_SSc\to {}_SSe$ be an isomorphism for which
$f(c)=e$. Since $c\in Se$ then $ce=c$. From the equality $b\cdot
c=e$ we get $cbc=c$. Then $cbf(c)=f(c)$ or $cbe=e$, but $be=b$.
Consequently, $cb=e$. So $e\in Sb$ which contradicts with the
condition $Sb\subset Se$. Thus ${}_SSe\not\models\Psi(e)$ and,
consequently, ${}_SSe\not\models\forall x\:\Psi(x)$.

On induction we will build now a regular polygon ${}_SA$, on which
the formula $\forall x\:\Psi(x)$ is true.

1. We set ${}_SA_0\rightleftharpoons {}_SSe$. Obviously, the
polygon ${}_SA_0$ is regular and connected.

2. If a regular connected polygon ${}_SA_k$ is built, then we
denote by $X_k$ the set
$$\{x\mid x\in A_k,\: {}_SSx \mathrel{\widetilde{\to}} {}_SSe,\:
{}_SA_k\models\forall y\:(\Phi_e(y)\to by\ne x)\}.$$ For each
$x\in X_k$ we consider the regular polygon ${}_SSa_x$ such that
${}_SSa_x \mathrel{\widetilde{\to}} {}_SSe$. Herewith without loss
of generality we can assume that $Sa_x\cap Sa_y=\emptyset$ for
$x\ne y$ and $Sa_x\cap A_k=\emptyset$ for any $x\in X_k$. Let
$f_x:{}_SSa_x\to {}_SSe$ be an isomorphism for which $f(a_x)=e$.
We set
$${}_S\tilde{A}_{k+1}\rightleftharpoons
\left({}_SA_k\bigsqcup\left(\bigsqcup\limits_{x\in X_k}
{}_SSa_x\right)\right)\Biggm/\theta,$$ where $\theta$ is a
congruence generated by the set $\{(x,ba_x)\mid x\in X_k\}$. Since
${}_SSb={}_SS(be) \mathrel{\widetilde{\to}} {}_SSe$ and ${}_SSa_x
\mathrel{\widetilde{\to}} {}_SSe$, then
$${}_SS(ba_x) \mathrel{\widetilde{\to}} {}_SS(be)
\mathrel{\widetilde{\to}} {}_SSe.$$ By the definition of the set
$X_k$ for any $x\in X_k$ it is executed ${}_SSx
\mathrel{\widetilde{\to}} {}_SSe$. Since there exists an
isomorphism $f_x:{}_SSx \to {}_SSba_x$ such that $f_x(x)=b\cdot
a_x$, then
$$\langle x_0,y_0\rangle\in\theta\iff$$
$$\iff x_0=y_0,\mbox{ or }x_0=sx\mbox{ and }y_0=sba_x \mbox{ for some }s\in S,x\in X_k.$$
Since for any $x\in X_k$ we have $x\in A_k$ and $ba_x\in Sa_x$,
then $x\not\sim ba_x$. Consequently, $sx\not\sim sba_x$ and
$\theta$ is an amalgam congruence. Since ${}_SA_k$, ${}_SSa_x$ is
a regular polygon, then on Proposition 2.21 the polygon
${}_S\tilde{A}_{k+1}$ is regular.

We will show that $c\sim d$ for any $c,d\in{}_SA_{k+1}$. Let
$c_1$, $d_1$ be some elements, for which $c=c_1/\theta$,
$d=d_1/\theta$. If $c_1,d_1\in A_k$, then $c\sim d$ in view of
connectivity of the polygon ${}_SA_k$. If $c_1\in A_k$ and $d_1\in
Sa_x$, then $d_1\sim a_x\sim ba_x$ and in view of connectivity of
the polygon ${}_SA_k$, $c_1\sim x$. Then $c_1/\theta\sim x/\theta$
and $d_1/\theta\sim ba_x/\theta$. Since $x/\theta=ba_x/\theta$,
then $c_1/\theta\sim d_1/\theta$, that is $c\sim d$. But if
$c_1\in Sa_x$ and $d_1\in Sa_y$ for some $x,y\in X_k$ then for any
$b_1\in A_k$ it is executed $c_1/\theta\sim b_1/\theta$ and
$d_1/\theta\sim b_1/\theta$. By the transitivity of relation
$\sim$ we get that $c_1/\theta\sim d_1/\theta$, that is i.e.
$c\sim d$. Thereby, ${}_S\tilde{A}_{k+1}$ is a connected polygon.

We will show that for any $x_0\in X_k$ it is executed
${}_S\tilde{A}_{k+1}\models\Psi\left(x_0/\theta\right)$. Since the
polygons ${}_SA_k$ and ${}_SSa_{x_0}$ are connected and $\theta$
is an amalgam congruence, then the restrictions on ${}_SA_k$ and
${}_SSa_{x_0}$ of a natural homomorphism corresponding to $\theta$
are isomorphic embeddings. Since ${}_SSx_0
\mathrel{\widetilde{\to}} {}_SSe$ and ${}_SSa_{x_0}
\mathrel{\widetilde{\to}} {}_SSe$, then
$${}_SS\left(x_0/\theta\right) \mathrel{\widetilde{\to}} {}_SSe,\mbox{ }
{}_SS\left(a_{x_0}/\theta\right) \mathrel{\widetilde{\to}}
{}_SSe,\mbox{
}{}_S\tilde{A}_{k+1}\models\Phi_e\left(x_0/\theta\right),\mbox{ }
{}_S\tilde{A}_{k+1}\models\Phi_e\left(a_{x_0}/\theta\right).$$
Since $b\left(a_{x_0}/\theta\right)=(ba_{x_0})/\theta=x_0/\theta,$
then ${}_S\tilde{A}_{k+1}\models
b\left(a_{x_0}/\theta\right)=x_0/\theta$ and
${}_S\tilde{A}_{k+1}\models\Phi_e\left(x_0/\theta\right)\to
\exists
y\:\left(by=\left(x_0/\theta\right)\wedge\Phi_e(y)\right),$ that
is ${}_S\tilde{A}_{k+1}\models\Psi\left(x_0/\theta\right).$

Since ${}_SA_k$ is isomorphically embeddable in
${}_S\tilde{A}_{k+1}$, then there exists a polygon ${}_SA_{k+1}$
such that ${}_SA_k\subseteq{}_SA_{k+1}$ and
${}_SA_{k+1}\cong{}_S\tilde{A}_{k+1}$. Then for any $x_0\in X_k$
it is executed ${}_SA_{k+1}\models\Psi\left(x_0/\theta\right)$.

As a result of the induction process we get a chain of regular
polygons
$${}_SA_0\subseteq{}_SA_1\subseteq\ldots\subseteq{}_SA_k\subseteq
{}_SA_{k+1}\subseteq\ldots$$ We set
${}_SA\rightleftharpoons\bigcup\limits_{i\in\omega}{}_SA_i$. As
the class $_S\Re$ is closed relative to unions, the polygon
${}_SA$ is regular.

We will show that ${}_SA\models\forall x\:\Psi(x)$. Notice that
since for any regular polygon ${}_S\hat{A}$ and any $a\in\hat{A}$
it is executed
$${}_S\hat{A}\models\Phi_e(a)\iff {}_SSa \mathrel{\widetilde{\to}} {}_SSe,$$
then for any element $a\in A$ and for any subpolygon ${}_SB$,
containing the element $a$ and any extended polygon ${}_SC$ we
have
$${}_SB\models\Phi_e(a),\mbox{ }{}_SC\models\Phi_e(a).$$
Since
$$\Psi(x)=\Phi_e(x)\to\exists y\:(by=x\wedge\Phi_e(y)),$$
then from ${}_SA\models\Psi(a)$ we get ${}_SC\models\Psi(a)$ for
any extended regular polygon ${}_SC$.

Assume there exists an element $d\in A$ such that
${}_SA\not\models\Psi(d)$. Choose the number $k\in\omega$, for
which $d\in A_k$. Since ${}_SA\not\models\Psi(d)$ and ${}_SA_k$ is
a subpolygon of the polygon ${}_SA$, then
${}_SA_k\not\models\Psi(d)$, that is ${}_SA_k\models\Phi_e(d)$ and
${}_SA_k\models\forall y\:(\Phi_e(y)\to by\ne d)$. Then $d\in
X_k$, but by the construction of the polygon ${}_SA_{k+1}$ we get
${}_SA_{k+1}\models\Psi(d)$ and, consequently,
${}_SA\models\Psi(d)$. Thereby, ${}_SA\models\forall x\:\Psi(x)$.

Since ${}_SSe\not\models\forall x\:\Psi(x)$, then one get a
contradiction to completeness of the class $\Re$.\hfill $\Box$

\bigskip
{\bf Lemma 6.3.} If ${\rm ld}(R)=1$, then the class $\Re$
axiomatizable if and only if the set $I(R)$ is finite.

{\bf Proof}. Assume the class $\Re$ is axiomatizable. Then on
Theorem 4.1 the set $\{x\in R\mid 1\cdot x=x\}$ is equal to $R$
and finitely generated as a right ideal of the semigroup $R$.
Since by Proposition 2.6 the semigroup $R$ is a rectangular band
of groups, then principal right ideals in $R$ are pairwise
disjoint and $R$ is their union. Since the set of all principal
right ideals in $R$ has the same cardinality as the set $I(R)$, we
conclude that the set $I(R)$ is finite.

For a proof of sufficiency we notice that the number of principal
right ideal coincides with the number of element of the set $I(R)$
and, consequently, is finite. Then, obviously, the conditions of
Theorem 4.1 are satisfied and the class $\Re$ axiomatizable.
\hfill $\Box$

\bigskip
{\bf Lemma 6.4.} If the class $\Re$ is complete, ${}_SA\in \Re$
and ${}_SA\models\exists x\:\Phi(x)$, then
${}_SA\models\exists^{\geq\omega} x\:\Phi(x)$.

{\bf Proof}. Let ${}_SA_i$, $i\in\omega$, be disjoint isomorphic
copies of the polygon ${}_SA$. Then
$\bigsqcup\limits_{i\in\omega}{}_SA_i\models\exists^{\geq\omega}x\:\Phi(x)$
and in view of $\bigsqcup\limits_{i\in\omega}{}_SA_i\equiv{}_SA$
we get ${}_SA\models\exists^{\geq\omega} x\:\Phi(x)$.\hfill $\Box$

\bigskip
{\bf Lemma 6.5.} If the class $\Re$ is complete then
$|eSg|\geq\omega$ for any idempotents $e\in S$ and $g\in R$.

{\bf Proof}. Denote by $a$ the element $eg$. Obviously, $a\in eSg$
and $Sg\models ea=a$. Consequently, $Sg\models\exists x\:(ex=x)$
and, on Lemma 6.4, $Sg\models\exists^{\geq\omega}x\:(ex=x)$. Since
the condition $x\in eSg$ is equivalent to $ex=x$ and $x\in Sg$,
then $|eSg|\geq\omega$.\hfill $\Box$

{\bf Lemma 6.6.} If the kernel $K(R)$ is not empty and the class
$\Re$ is complete, then $K(R)$ is a rectangular band of infinite
groups.

{\bf Proof}. Proposition 1.4 and the regularity of polygon ${}_SR$
imply that $K(R)$ is a rectangular band of groups. On Remark 1.1
it is enough to prove that the group $G_e$ is infinite for any
idempotent $e\in K(R)$.

Let $e$ be an arbitrary idempotent from $K(R)$. Then
${}_SSe\models (ee=e)$ and, consequently, ${}_SSe\models\exists
x\: (ex=x)$. On Lemma 6.4 we get
${}_SSe\models\exists^{\geq\omega} x\: (ex=x)$. Since $K(R)$ is a
rectangular band of groups, then the ideal $Se$ is minimal, and,
consequently, for any $a\in Se$ we have $Se\models (ea=a)$ if and
only if $a\in G_e$. Thereby, for any $e\in K(R)$ the group $G_e$
is infinite and the kernel $K(R)$ is a rectangular band of
infinite groups. \hfill $\Box$

{\bf Lemma 6.7.} If $H$ is a subgroup of the monoid $S$, $e$ is a
unit in $H$, then ${}_SSa \mathrel{\widetilde{\to}} {}_SSe$ for
any $a\in H$.

{\bf Proof}. Define the map $\varphi:{}_SSa\to {}_SSe$ by the rule
$\varphi(x)=x\cdot a^{-1}$, where $a^{-1}$ is inverse for $a$ in
the group $H$. Then $\varphi(a)=a\cdot a^{-1}=e$ and
$$
\varphi(sx)=(sx)\cdot a^{-1}=s(x\cdot a^{-1})=s\varphi(x)$$ for
any $s\in S$ and $x\in Sa$. Consequently, $\varphi$ is a
homomorphism.

If $x,y\in Sa$ and $\varphi(x)=\varphi(y)$, then $x\cdot
a^{-1}=y\cdot a^{-1}$. So $xa^{-1}a=ya^{-1}a$, $xe=ye$, $x=y$.
Thus $\varphi$ is an one-to-one correspondence. If $x\in Se$, then
$x=xe=xa\cdot a^{-1}=\varphi(xa)$. Thereby, $\varphi$ is an
isomorphism. \hfill $\Box$

{\bf Lemma 6.8.} If $R$ is a rectangular band of groups then for
any element $a\in R$ and for any idempotent $e\in R$ the condition
${}_SSa \mathrel{\widetilde{\to}} {}_SSe$ is true if and only if
$e\cdot a=a$.

{\bf Proof}. If ${}_SSa \mathrel{\widetilde{\to}} {}_SSe$, then
the equality $e\cdot a=a$ is obvious.

Let the equality $e\cdot a=a$ be true and $R$ be a rectangular
band of groups $S_{ij}$ with units $e_{ij}$, $i\in I$, $j\in J$.
Then $e=e_{ij}$ for some $i\in I$ and $j\in J$. Since $e_{ij}\cdot
a=a$, then $e_{ij}T=aT$ and on Remark 1.1 we get $a\in S_{ik}$ for
some $k\in J$.

We will prove now that ${}_SSe_{ik} \mathrel{\widetilde{\to}}
{}_SSe_{ij}$. Define the map $\varphi$ from $Se_{ik}$ to $Se_{ij}$
by the rule $\varphi(x)=x\cdot e_{ij}$. Assume $s\in S$, $x\in
Se_{ik}$. Then $sx\in Se_{ik}$ and
$$
\varphi(sx)=(sx)\cdot e_{ij}=s(x\cdot e_{ij})=s\varphi(x),
$$
that is $\varphi$ is a homomorphism. If $a,b\in Se_{ik}$ and
$\varphi(a)=\varphi(b)$, then $ae_{ij}=be_{ij}$. So
$ae_{ij}e_{ik}=be_{ij}e_{ik}$. Consequently, $ae_{ik}=be_{ik}$,
$a=b$. Thus the map $\varphi$ is injective. If $c\in Se_{ij}$,
then
$$c=c\cdot e_{ij}=c\cdot e_{ik}e_{ij}=\varphi(ce_{ik}),$$
that is $\varphi$ is surjective. Thereby, $\varphi$ is an
isomorphism and $\varphi(e_{ik})=e_{ik}\cdot e_{ij}=e_{ij}$.
Consequently, ${}_SSe_{ik} \mathrel{\widetilde{\to}} {}_SSe_{ij}$.
Since ${}_SSa \mathrel{\widetilde{\to}} {}_SSe_{ik}$ on Lemma 6.7
and $e_{ij}=e$, then ${}_SSa \mathrel{\widetilde{\to}} {}_SSe$.
\hfill $\Box$

The implication $1\Rightarrow 2$ results from Lemma 6.1.

$2\Rightarrow 3$. Assume the class $\Re$ is complete. Then on
Lemma 6.2 we have ${\rm ld}(R)=1$ and on Lemma 6.3 $R$ is a
rectangular band of groups.

On Lemma 6.3 the set $I(R)$ is finite, and on lemma 6.8 groups
forming $R$ are infinite.

$3\Rightarrow 1$. Since in rectangular band of  groups any
principal right ideal in $R$ is minimal, then $S$ satisfies
condition 1 of Theorem 4.1. Since the set $I(R)$ is finite, the
set of principal right ideal is also finite. Consequently, the
monoid $S$ satisfies condition 2 of Theorem 4.1, and the class
$\Re$ is axiomatizable.

Since in rectangular band of  groups any left ideal in $R$ is
minimal, then $S$ satisfies conditions 1 and 2 of Theorem 5.1. If
idempotents $e$ and $g$ belong to $R$, then $e=e_{ij}$ and
$g=e_{kl}$ for some $i,k\in I(R)$, $j,l\in J(R)$. Then using
properties of the rectangular bands of the groups we have
$$
eSg=e_{ij}Te_{kl}=\left(\bigcup\limits_{p\in
J(T)}S_{ip}\right)e_{kl}= \bigcup\limits_{p\in
J(T)}\left(S_{ip}e_{kl}\right)=S_{il}.
$$
In view of infinity of the group $S_{il}$ we get condition 3 of
Theorem 5.1. So the class $\Re$ is model complete.

The condition of formula definability of isomorphic orbits follows
from Lemma 6.8. \hfill $\Box$

{\bf Corollary 6.1} [Ov2]. If a semigroup $R$ contains a finite
number of idempotents, then the following conditions are
equivalent:

1) the class $\Re$ is model complete;

2) the class $\Re$ is complete;

3) the semigroup $T$ is a rectangular band of a finite number of
infinite groups.

{\bf Proof}. Since in any regular polygon ${}_SA$ for any element
$a\in A$ there exists an idempotent $e\in R$ such that ${}_SSa
\mathrel{\widetilde{\to}} {}_SSe$, and the semigroup $R$ contains
only finite number of idempotents then for each idempotent $e\in
R$ there exists a formula $\Phi_e(x)$ such that
$${}_SA\models\Phi_e(x)\iff {}_SSa \mathrel{\widetilde{\to}} {}_SSe,$$
that is the class $\Re$ satisfies the condition of formula
definability of isomorphic orbits. Thus, on Theorem 6.1 the
corollary holds. \hfill $\Box$

{\bf Corollary 6.2.} If a semigroup $R$ contains single
idempotent, then the following conditions are equivalent:

1) the class $\Re$ is model complete;

2) the class $\Re$ is complete;

3) $R$ is an infinite group. \hfill $\Box$

{\bf Corollary 6.3.} If $S$ is a commutative monoid, then the
following conditions are equivalent:

1) the class $\Re$ is model complete;

2) the class $\Re$ is complete;

3) $R$ is an infinite Abelian group.

{\bf Proof}. Considering Corollary 6.2 it is enough to show that
if $\Re$ is a complete class of regular polygons over commutative
monoid, then $R$ contains single idempotent. Let $e_1$, $e_2$ be
some idempotents from $R$. Then for arbitrary $x\in Se_1$ we have
$xe_1=x$. In view of commutativity of monoid $S$ we get $e_1x=x$.
Consequently, ${}_SSe_1\models\forall x\: (e_1x=x)$. Since class
$\Re$ is complete, then ${}_SSe_2\models\forall x\: (e_1x=x)$ and
so $e_1e_2=e_2$. Similarly we get $e_2e_1=e_1$. Thereby, the
commutativity of monoid $S$ implies $e_1=e_2$. \hfill $\Box$

{\bf Theorem 6.2} [Ov3]. If $S$ is a linearly ordered monoid of
the depth 2 and $\Re$ is an axiomatizable class, then the
following conditions equivalent:

1) the class $\Re$ is model complete;

2) the class $\Re$ is complete;

3) the semigroup $K(S)$ is a rectangular band of infinite number
of infinite groups.

{\bf Proof} contains some lemmas.

{\bf Lemma 6.9.} If a monoid $S$ has a finite depth, then
$K(S)=K(R)$.

{\bf Proof}. The condition $K(R)\subseteq K(S)$ is obvious. Since
the depth ${\rm ld}(S)$ is finite, then the depth ${\rm ld}(R)$ is
finite too. So on Proposition 2.6 the kernel $K(R)$ is a
rectangular band of groups. Let $e$ be an element from $E(S)\cap
K(R)$. Then $e$ belongs to $K(S)$ and, consequently, on
Proposition 1.4 the kernel $K(S)$ is a rectangular band of groups.
Thence we get $K(S)\subseteq R$ and, obviously, $K(S)\subseteq
K(R)$. \hfill $\Box$

{\bf Lemma 6.10.} If a monoid $S$ has a depth 2, then $R=S$.

{\bf Proof}. As the monoid $S$ has the depth 2 then for any
element $a\in S$ the condition $Sa\ne S1$ implies $Sa\subset S1$
and $Sa\subseteq K(S)$. On Lemma 6.9 and Proposition 2.6 the
kernel $K(S)$ is a rectangular band of groups. Consequently,
$K(S)\subseteq R$. If $a\in S\setminus K(S)$, then $Sa=S=S1$. So
the polygon ${}_SSa$ belongs to the class $\Re$, that is $a\in R$.
Thereby, $R=S$. \hfill $\Box$

{\bf Lemma 6.11.} If $S$ is a rectangular band of groups and
$|J(S)|=1$ then for some group $G$ the semigroup $S$ is isomorphic
to a semigroup with the universe $S'=\{\langle b,i\rangle\mid b\in
G,i\in I(S)\}$ and the operation $\cdot$, defined by equalities
$$
\langle c,j\rangle\cdot\langle d,k\rangle=\langle c\cdot
d,j\rangle
$$
for any elements $\langle c,j\rangle$, $\langle c,j\rangle\in S'$.

{\bf Proof}. Let $S$ be a rectangular band of groups $G_i$ with
idempotents $e_i$ for $i\in I(S)$. We fix an element $i_0\in I(S)$
and we set $G\rightleftharpoons G_{i_0}$, $S'\rightleftharpoons
G\times I(S)$. Define a map $\varphi:S'\to S$, acting by the rule
$\varphi(\langle b,i\rangle)\rightleftharpoons e_i\cdot b$, $b\in
G$, $i\in I(S)$. On Remark 1.1 (3) for any $i\in I(S)$ the group
$G_i$ is isomorphic to the group $G$, as well as to the group
$G'_i$ with the universe $\{\langle b,i\rangle\mid b\in G\}$ and
the operation, defined by the rule $\langle
b_1,i\rangle\cdot\langle b_2,i\rangle=\langle b_1\cdot
b_2,i\rangle$. Herewith, obviously, the restriction of the map
$\varphi$ on the set $G'_i$ realizes an isomorphism between groups
$G'_i$ and $G_i$. Besides it is directly checked that the map
$\varphi$ itself is a bijection.

Define the operation $\cdot$ on elements of the set $S'$ by the
following rule: $ \langle c,j\rangle\cdot\langle
d,k\rangle=\langle c\cdot d,j\rangle, $ $\langle c,j\rangle$,
$\langle c,j\rangle\in S'$. We will show that the map $\varphi$
conserves the operation $\cdot$. Really, using properties of the
rectangular band of groups, for any $j,k\in I(S)$ and $c,d\in G$
we get
$$
\varphi(\langle c,j\rangle)\cdot\varphi(\langle d,k\rangle)=
e_jc\cdot e_kd=e_j(ce_{i_0})\cdot e_kd=
$$
$$
=e_jc(e_{i_0}e_k)d=e_jce_{i_0}d=e_j(ce_{i_0})d=e_j(cd)=
\varphi(\langle c\cdot d,j\rangle).
$$
Thereby, the map $\varphi$ realizes an isomorphism between $S'$
and $S$, and the map $\varphi^{-1}:S\to S'$ is a required
isomorphism. \hfill $\Box$

Hereinafter under $J(S)=1$ the rectangular band of the groups $S$
will be identified with the semigroup $S'$.

{\bf Lemma 6.12.} If a monoid $S$ has a finite depth and $R=S$
then $1\in aS$ for any element $a\in S$ such that $Sa=S$.

{\bf Proof}. Denote by $G_1$ the set $\{x\in S\mid Sx=S\}$. We
will show that $G_1\cap E(S)=\{1\}$. Really if $e\in G_1\cap E(S)$
then $1\cdot e=1$. Since $1\cdot e=e$ it follows that $e=1$. On
the other hand, obviously, the unit belongs to $G_1\cap E(S)$.
Thereby, $G_1\cap E(S)=\{1\}$.

Denote by $n$ the depth ${\rm ld}(S)$. Since for any $a\in S$ it
is executed $Sa\subseteq S\cdot 1$, then any chain of the length
$n$ of principal left ideals contains the set $S\cdot 1$ as
maximal element, that is $S\cdot 1$ is single principal left ideal
of the depth $n$. Since $R=S$ then for any element $a\in G_1$
there exists an idempotent $g\in S$ such that
${}_SSa\:\widetilde{\to}\:{}_SSg$. Since  $Sa=S$ and the depth of
$Sa$ is equal to $n$, we get ${}_SSa\:\widetilde{\to}\:{}_SS1$.
Let $\varphi$ be an isomorphism from ${}_SSa$ to ${}_SS1$ such
that $\varphi(a)=1$. Since $1\in Sa$ then the element $\varphi(1)$
is defined. So $a\cdot\varphi(1)=\varphi(a\cdot 1)=\varphi(a)=1$,
that is $1\in aS$. \hfill $\Box$

{\bf Lemma 6.13.} If $2\leq{\rm ld}(R)<\infty$ and the class $\Re$
is complete, then $|I(T)|\geq\omega$.

{\bf Proof}. Since ${\rm ld}(R)\geq 2$, then there exist an
idempotent $g\in R$ and an element $a\in R$ such that $Sa\subset
Sg$. In view of finiteness of ${\rm ld}(R)$ one can suppose that
$Sa$ is a minimal left ideal. Then the element $a$ belongs to
$K(R)$. On Proposition 2.6 the kernel $K(R)$ is a rectangular band
of groups, and in view of properties of the rectangular band of
groups we have $Sa=\bigcup\limits_{i\in I(R)}S_{ij_0}$ for some
$j_0\in J(R)$.

Assume the set $I(R)$ is finite. Then the set $\{e_{ij_0}\mid i\in
I(R)\}$ is finite and
$$
{}_SSa\models\forall x\: \left(\bigvee\limits_{i\in
I(R)}e_{ij_0}x=x\right).
$$
But since $Se_{ij_0}\subset Sg$, then $e_{ij_0}g=e_{ij_0}\ne g$
for all $i\in I(R)$ and ${}_SSg\not\models e_{ij}g=g$.
Consequently, ${}_SSg\not\equiv {}_SSa$ that contradicts to the
completeness of the class $\Re$. \hfill $\Box$

The implication $1\Rightarrow 2$ results from Lemma 6.1.

$2\Rightarrow 3$. On Lemmas 6.9 and 2.6  $K(S)$ is a rectangular
band of groups. On Lemma 6.13 the set $I(K(S))$ is infinite, that
is there are infinitely many groups forming $K(S)$. The infinity
of the groups follows from Lemma 6.6.

$3\Rightarrow 1$. Since the class $\Re$ is axiomatizable, we can
use Theorem 5.1.

On Lemma 6.10 we have $R=S$. Since $S$ is a linearly ordered
monoid of depth 2, there exists single chain of principal left
ideals $Sa\subset S\cdot 1$. Consequently, condition 1 of Theorem
5.1 holds.

Since ${\rm ld}(S)=2$, the condition $Sa\subset Se$ implies $e=1$,
where $a\in S$, $e\in E(S)$. Then $Sa\subseteq K(S)$. In view of
linear ordering of $S$ the equality $Sa=K(S)$ holds, that is on
Lemmas 6.9 and 2.6 semigroup $Sa$ is a rectangular band of the
groups. Thereby, kernel $K(S)$ consists of one principal left
ideal and $J(S)=1$. Let us fix the group $\langle
G;\:\cdot\rangle$ with unit $e'$ such that all groups from $Sa$
are isomorphous to it. Then on Lemma 6.11 elements from $K(S)$ can
be present as $\langle b,i\rangle$, where $b\in G$, $i\in
I(K(S))$. Consequently,
$$
S=\{x\mid Sx=S\cdot 1\}\cup\{\langle b,i\rangle\mid b\in G,i\in
I(K(S))\}.
$$

Assume $1\not\in\bigcup\limits_{i=1}^m a_iS$ for some $a_i\in S$,
$1\leq i\leq m$. Then we have $1\not\in a_iS$ for each
$i\in\{1,\ldots,m\}$, and in view of Lemma 6.12 for any
$i\in\{1,\ldots,m\}$ there exist $b_i\in G$, $k_i\in I(K(S))$ such
that $a_i=\langle b_i,k_i\rangle$. In view of properties of the
rectangular band of groups the equality $|J(K(S))|=1$ implies
$\langle b_i,k_i\rangle\: S=S_{k_i}$ for all $i\in\{1,\ldots,m\}$,
where $ S_{k_i}=\{\langle b,k_i\rangle\mid b\in G\}$. Since for
any $j\in I(K(S))\setminus\{k_i\}$ it is satisfied $\langle
e',j\rangle\not\in S_{k_i}$, $S\langle e',j\rangle=S\langle
b_i,k_i\rangle$ and $\langle e',j\rangle=1\cdot\langle
e',j\rangle\in 1\cdot S$, then the set $E_{k_i}\rightleftharpoons
E(S)\cap\{x\mid x\in (1\cdot S)\setminus a_iS\}$ coincides with
the set $\{\langle e',j\rangle\mid j\in
I(K(S))\setminus\{k_i\}\}$. On Lemma 6.13 the set $I(K(S))$ is
infinite, and, consequently, the set $E_{k_i}$ is also infinite
for any $i\in I(K(S))$. Since
$$
(1\cdot S)\setminus\bigcup\limits_{i=1}^m a_i S=
\bigcap\limits_{i=1}^m\left((1\cdot S)\setminus a_i S\right),
$$
then
$$
E(S)\cap\{x\mid x\in (1\cdot S)\setminus\bigcup\limits_{i=1}^m
a_iS\}= \{\langle e',j\rangle\mid j\in I(S)\setminus\{k_i\mid
1\leq i\leq m\}\}
$$
and in view of infinity of the set $I(K(S))$ condition 2 of
Theorem 5.1 holds.

For the checking of condition 3 of Theorem 5.1 we will consider
all possible variants for values of idempotents $e$, $g$.

If $e=1$ and $g=1$, then $|eSg|=|1S\cdot 1|=|S|\geq\omega$.

If $e=1$ and $g=\langle e',i\rangle$ for some $i\in I(K(S))$, then
$$
|eSg|=|1S\langle e',i\rangle|=|S\langle
e',i\rangle|=|K(S)|\geq\omega.
$$

If $e=\langle e',i\rangle$ for some $i\in I(K(S))$ and $g=1$ then
on Lemma 6.6 we have
$$
|eSg|=|\langle e',i\rangle S\cdot 1|=|\langle e',i\rangle
S|=|S_i|\geq\omega.
$$

If $e=\langle e',i\rangle$ and $g=\langle e',j\rangle$ for some
$i,j\in I(K(S))$ then by properties of the rectangular band of
groups it is satisfied
$$
|eSg|=|\langle e',i\rangle S\langle e',j\rangle|= |S_i\langle
e',j\rangle|=|S_i|\geq\omega.
$$

Thus, all conditions of Theorem 5.1 hold and the class $\Re$ is
model complete.

The theorem is proved. \hfill $\Box$

The following two theorems show that there exist monoids over
which classes for all regular polygons are complete, but  not
model complete. Here we will only build such monoids but detailed
proofs of their properties can be found in [Ov4, Ov5]

{\bf Theorem 6.3.} There exists a not linearly ordered monoid of
depth 2, over which the class of all regular polygons is complete,
but not model complete.

{\bf Proof}. For an Abelian group ${\cal G}=\langle G,+\rangle$,
nonempty sets $I, J$ and a function $\varphi: I\times J\to G$ we
call by $\langle{\cal G},I,J,\varphi\rangle$-band the semigroup
$\langle{\cal G}\times I\times J,\ast\rangle$, in which the
operation $\ast$ is defined by the following way:
$$
\langle a,i,j\rangle\ast\langle b,k,l\rangle\rightleftharpoons
\langle a+b+\varphi(k,j),i,l\rangle.
$$

For arbitrary elements $i\in I$ and $j\in J$ we denote by $S_{ij}$
the set $\{\langle a,i,j\rangle\mid a\in G\}$. Then the algebra
${\cal S}_{ij}=\langle S_{ij},\ast\rangle$ is a subgroup of
$\langle{\cal G},I,J,\varphi\rangle$-band with the idempotent
$\langle -\varphi(i,j),i,j\rangle$ as the unit element. Since
$S_{ij}\ast S_{kl}\subseteq S_{il}$, then $\langle{\cal
G},I,J,\varphi\rangle$ is a rectangular band of groups, which will
be denoted by ${\rm RB}\langle{\cal G},I,J,\varphi\rangle$.

Remind [KM] that for any group $G$ with its unit $e$ and any
ordinal $\alpha$ by $G^\alpha$ one denote the direct degree of
group $G$, that is the set of all function $f:\alpha\to G$ such
that the set $\{x\mid f(x)\ne e\}$ is finite.

Let us consider the group ${\bf Z}_2^\omega$, where ${\bf
Z}_2=\langle\{0,1\},+\rangle$. Elements of the group ${\bf
Z}_2^\omega$ will be denoted by $\bar{a},\bar{b},\ldots$, where
$\bar{a}=(a_0,a_1,\ldots)$, $\bar{b}=(b_0,b_1,\ldots)$. We will
denote zero of the group ${\bf Z}_2^\omega$ by $\bar{0}$. For
arbitrary element $\bar{a}\in{\bf Z}_2^\omega$ we introduce the
following notations:
$$
h(\bar{a})\rightleftharpoons\left\{
\begin{array}{rl}
(0), & \mbox{если } \bar{a}=\bar{0}, \\ (a_0,\ldots,a_{n-1}), &
\mbox{если } a_{n-1}=1 \mbox{ and } a_k=0 \mbox{ for all } k\geq
n,
\end{array}\right.
$$
$$
l(\bar{a})\rightleftharpoons\left\{
\begin{array}{rl}
0, & \mbox{если } \bar{a}=\bar{0}, \\ n, & \mbox{если }
h(\bar{a})=(a_0,\ldots,a_{n-1}),
\end{array}\right.
$$
$$
r(\bar{a})\rightleftharpoons\sum_{i=0}^{l(\bar{a})}a_i\cdot 2^i.
$$

If $a=\langle\bar{a},i,j\rangle\in S$, then we set
$l(a)\rightleftharpoons l(\bar{a})$, $r(a)\rightleftharpoons
r(\bar{a})$.

Define the function $\psi$ from $\omega$ to ${\bf Z}_2^\omega$ by
the following way:
$$
\psi(n)=r^{-1}\left(\left[\sqrt{n}\right]\right),
$$
where $[x]$ means the integer part of real number $x$.

Define the multiplying operation of elements of group ${\bf
Z}_2^{\omega}$ on 0 and 1 as follows:
$$
\bar{a}\cdot 0\rightleftharpoons \bar{0},
$$
$$
\bar{a}\cdot 1\rightleftharpoons \bar{a}.
$$
We set ${\cal G}\rightleftharpoons{\bf Z}_2^{\omega}$,
$I\rightleftharpoons\omega$, $J\rightleftharpoons\{0,1\}$,
$\varphi(i,j)\rightleftharpoons\psi(i)\cdot j$.

Define a monoid $S$ by equality
$$
S=\langle{\rm RB}({\cal G},I,J,\varphi)\cup{\bf
Z}_2^\omega,\ast\rangle,
$$
where action of operation $\ast$ between elements of the same
nature is defined by natural way, but between element from ${\rm
RB}({\cal G},I,J,\varphi)$ and ${\bf Z}_2^\omega$ by the following
rule:
$$
\langle\bar{a},i,j\rangle\ast\bar{b}=\bar{b}\ast\langle\bar{a},i,j\rangle=
\langle\bar{a}+\bar{b},i,j\rangle.$$

The unit of the monoid $S$ is the element $\bar{0}$. All
idempotents of monoid $S$ form the set
$$ E(S)=\{\langle\varphi(i,j),i,j\rangle\mid i\in I,j\in
J\}\cup\{\bar{0}\}.$$ The monoid $S$ is a union of the groups,
isomorphic ${\bf Z}_2^\omega$. Consequently, $R=S$. \hfill $\Box$

\bigskip
{\bf Theorem 6.4.} There exists a linearly ordered monoid of depth
3, over which the class of all regular polygons is complete, but
not model complete.

{\bf Proof}. Let $G=\langle G;\ast,0\rangle $ be a countable
group.

Define a monoid $ S=\langle S;\cdot \rangle $ as follows: $
S=(\cup_{i\in\omega}G_{1i})\cup (\cup_{i\in\omega}G_{2i})\cup
(\cup_{i\in\omega}G_{3i})\cup G_{4}\cup G_5,\ $ where

$$G_{1i}=\{[1,i,\bar j]\ |\ \bar j \in G^{\omega}\},\ \
G_{2i}=\{[2,i,\bar j]\ |\ \bar j\in G^{\omega}\},$$
$$G_{3i}=\{[3,i,\bar j]\ |\ \bar j\in G^{\omega}
\},\ \ G_4=\{[4,\bar i]\ |\ \bar i\in G^{\omega}\},$$
$$G_5=\{[5,\bar i]\ |\ \bar i\in G^{\omega}\}.$$

An action of operation $\ \cdot \ $ will be defined by the
following table:

\begin{center}
\begin{tabular}{r|lllll}$\cdot $ & $[1,m,\bar k]$ & $[2,m,\bar k]$ &
$[3,m,\bar k]$ & $[4,\bar k]$ & $[5,\bar k]$ \\
\hline $[1,i,\bar j]$ & $[1,i,\bar j\ast \bar k]$ & $[1,i,\bar
j\ast \bar k]$ & $[1,i,\bar j\ast \bar k]$
& $[1,i,\bar j\ast \bar k]$ & $[1,i,\bar j\ast \bar k]$\\
$[2,i,\bar j]$ & $[2,i,\bar j\ast \bar k]$ & $[2,i,\bar j\ast \bar
k]$ & $[2, i,\bar j\ast \bar k]$
& $[2,i,\bar j\ast \bar k]$ & $[2, i,\bar j\ast \bar k]$ \\
$[3,i,\bar j]$ & $[3,i,\bar j\ast \bar k]$ & $[3,i,\bar j\ast \bar
k]$ & $[3,i,\bar j\ast \bar k]$
& $[3,i,\bar j\ast \bar k]$ & $[3,i,\bar j\ast \bar k]$\\
$[4,\bar j]$ & $[1,m,\bar j\ast \bar k]$ & $[2,m,\bar j\ast \bar
k]$ & $[2,m,\bar j\ast \bar k]$ &
$[4,\bar j\ast \bar k]$ & $[4,\bar j\ast \bar k]$\\
$[5,\bar j]$ & $[1,m,\bar j\ast \bar k]$ & $[2,m,\bar j\ast \bar
k]$ & $[3,m,\bar j\ast \bar k]$ &
$[4,\bar j\ast \bar k]$  & $[5,\bar j\ast \bar k]$  \\

\end{tabular}
\end{center}

Note that the monoid $S$ is a union of the groups. The groups,
forming $S,$ are isomorphic to the group $G^{\omega}$. The unit of
the monoid is an element $[5,\bar 0]$. The polygon ${}_SS$ has
three different orbits $S[1,0,\bar 0]\subset S[4,\bar 0]\subset
S[5,\bar 0]$, each of which is generated by an idempotent.
Consequently, $R=S$ and ${\rm ld}(S)=3$.\hfill $\Box$

\vskip 1cm {\bf\centerline {\S~7. Stability of Class for Regular
Polygons}}
\medskip

In this paragraph we give the characterization of $\cal
R$--stabilizer whose regular core is presented as a union of the
finite number of principal right ideals (Theorem 7.1). As the
consequence corollary it is proved that for the axiomatizable
class of all regular polygons the stability of this class is the
equivalent to the regularly linearly ordering of monoid S. The
regular linearly ordering of monoid is the necessary (Proposition
7.1) but not sufficient (Example 7.1) condition for the stability
of any regular polygon. The Example 7.2 shows in particular that
there exists a regular linearly ordering monoid with
non--axiomatizable class of all regular polygons over which all
regular polygons are stable.

{\bf Proposition 7.1} If $S$ is a $\cal R$--stabilizer then $S$ is
a regularly linearly ordered monoid.

{\bf Proof.} Assume  $S$ is a $\cal R$--stabilizer which is not a
regularly linearly ordered monoid, that is, there exist $a,b,c\in
R$ such that $Sb\not \subseteq Sc\subseteq Sa$ and
$Sc\not\subseteq Sb\subseteq Sa$. Then $b=ta$ and $c=sa$ for some
$t,s\in S$. Let $K=\{\langle i,j\rangle\mid j\leqslant
i<\omega\}$; ${}_SA_{ij}$ is a copy of the polygon ${}_SSa$
($\langle i,j\rangle\in K$), and $A_{ij}\cap A_{kl}=\varnothing$
if $\langle i,j\rangle \neq \langle k,l\rangle$; $d_{ij}$ is a
copy of the element $d\in Sa$ in $A_{ij}$. Write ${}_SA$ for a
polygon $\bigcup\limits_{\langle i,j\rangle\in
K}{}_SA_{ij}/\theta$, where $\theta$ is a congruence on
$\bigcup\limits_{\langle i,j\rangle\in K}{}_SA_{ij}$ generated by
the set $\{\langle b_{ij},b_{il}\rangle\mid\langle i,j\rangle\in
K,\langle i,l\rangle\in K\}\cup\{\langle
c_{ij},c_{lj}\rangle\mid\langle i,j\rangle\in K,\langle
l,j\rangle\in K\}$. Denote by $b_i$ an equivalence class of
$\theta$ with a representative $b_{ij}$ $(\langle i,j\rangle\in
K)$, and by $c_j$ an equivalence class of $\theta$ with a
representative $c_{ij}$ $(\langle i,j\rangle\in K)$. Write
$\varphi(x,y)$ to abbreviate the formula $\exists z(x=tz\wedge
y=sz)$. Since a restriction of $\theta$ to the polygon
${}_SA_{ij}$ ($\langle i,j\rangle\in K$) is zero congruence,
${}_SA\in {\mathfrak R}$. Note that $b_i=ta_{ij}/\theta$ and
$c_j=sa_{ij}/\theta$. Moreover,
$$ {}_SB\models\varphi(b_i,c_j)\Leftrightarrow i\geqslant j, $$
which contradicts the stability of $ {\rm Th}({}_SA)$ on Theorem
3.5\hfill$\Box$

\smallskip
{\bf Theorem 7.1} [Ste2]. Let
$$
R=\bigcup\limits_{i=0}^n a_iR
$$
for some $n\geqslant 0$ and $a_i\in R$ $(0\leqslant i\leqslant
n)$. The monoid $S$ is an $\mathfrak R$--stabilizer if and only if
$S$ is a regularly linearly ordered monoid.

{\bf Proof.} Assume the theorem condition is satisfied.

{\bf Necessity} follows form Proposition 7.1.

{\bf Sufficiency.} Let $S$ be a regularly linearly ordered monoid
and ${}_SA\in {\mathfrak R}$. We claim that ${\rm Th}({}_SA)$ is
stationary. Suppose ${}_SM\equiv {}_SA$, $a,b\in {\mathfrak
C}\setminus M$, $a=sb$, and $c=rb$, $c\in M$. Suppose  $d\in A$.
On Corollary 2.1, there exist an idempotent $e\in R$ and an
isomorphism $\varphi: Sd\rightarrow Se$ such that $\varphi(d)=e$.
On (7.1), there exist $i$, $0\leqslant i\leqslant n$, and
 $u \in R$ such that $e=a_iu$. Therefore, $e=a_iu e$.
Consequently, ${}_SSd\models\exists y(d=a_iy)$, that is,
${}_SA\models\forall x\exists y\bigvee\limits_{i=0}^n (x=a_iy)$.
Since ${}_SM\equiv{}_SA$ and ${}_SM\prec {}_S{\mathfrak C}$, we
have $b=a_iv$ for some $i$, $0\leqslant i\leqslant n$, and $v\in
{\mathfrak C}$. Since the set $\{Sa\mid Sa\subseteq Sa_i\}$ is
linearly ordered, one of the inclusions $Ssa_i\subseteq Sra_i$ or
$Sra_i\subseteq Ssa_i$ holds. If $Ssa_i\subseteq Sra_i$ then
$a=sb=sa_iv\in Sra_iv=Srb=Sc\subseteq M$, which contradicts the
choice of $a$. Consequently, $Sra_i\subseteq Ssa_i$. In this case
$c=rb=ra_iv\in Ssa_iv=Ssb=Sa$, that is, the theory ${\rm
Th}({}_SA)$ is stationary. On Theorem 3.7, ${\rm Th}({}_SA)$ is
stable, that is, $S$ is an $\mathfrak R$--stabilizer.\hfill$\Box$

\smallskip
{\bf Corollary 7.1.} Let class $\mathfrak R$ for regular polygons
be axiomatizable. Class $\mathfrak R$ is stable if and only if $S$
is a linearly ordered monoid.

The {\bf proof} follows from Theorem 7.1 and Corollary
4.2.\hfill$\Box$

From the proof of Theorem 7.1 it follows

\smallskip
{\bf Corollary 7.2.} If the condition of Theorem 7.1 holds and $S$
is a linearly ordered monoid then the thoty of any regular
polygons is stationary.\hfill$\Box$

\smallskip
{\bf Example 7.1.} We construct a non $\mathfrak R$--stabilizer
which is a regularly linearly ordered monoid.

Let $S$ be a union of five disjoint sets:
$$ S=A\cup B\cup C\langle\alpha,\beta\rangle\cup\{1\}, $$
where $A=\{a_i\mid i\in\omega\},\;B=\{b_i\mid
i\in\omega\},\;C=\{c_{ij}\mid i \geqslant j,\; i,j\in\omega\}$ and
$\langle\alpha,\beta\rangle$ is a free two-generated semigroup
with $\alpha$ and $\beta$ generators. Equip $S$ with a binary
operation defined as follows: 1 is an identity element of $S$;
$a_ix=a_i$, $b_ix=b_i$, $c_{ij}x=c_{ij}$,
$(\gamma\alpha)c_{ij}=a_i$, $(\gamma\beta)c_{ij}=b_j$,
$(\gamma\alpha)y=y$, and $(\gamma\beta)y=y$ for any $x\in S$,
$y\in S\setminus (C\cup\langle\alpha,\beta\rangle)$,
$\gamma\in\langle\alpha,\beta\rangle\cup\{1\}$, $a_i\in A,\;b_i\in
B,\;c_{ij}\in C$. It is easy to verify that $S$ is a monoid under
the operation given. Furthermore, $R=A\cup B\cup C$ and
$Sa_i=Sb_i=Sc_{ij}=R$ for any $a_i\in A,\;b_i\in B,\;c_{ij}\in C$.
Hence $S$ is  a regularly linearly ordered monoid. Let
$\Phi(x,y)\rightleftharpoons\exists z(x=\alpha z\wedge y=\beta
z)$. Then
$$ {}_SS\models\exists z(a_i=\alpha z\wedge b_j=\beta z)\Leftrightarrow i\geqslant j, $$
and so on Theorem 3.5 ${\rm Th}({}_SR)$ is not stable, that is,
$S$ is not $\mathfrak R$--stabilizer.\hfill$\Box$

\smallskip
{\bf Example 7.2.} We construct a regularly linearly ordered
$\mathfrak R$--stabilizer which does not satisfy the condition of
Theorem 7.1.

Let $S=\bigcup\limits_{i\in\omega}Z_i\cup\langle\alpha,
\beta\rangle\cup\{1\}$, where $Z_i=\{n_i\mid n\in Z\}$ is a copy
of the set $Z$ of integers on which addition is defined naturally,
$Z_i\cap Z_j=\varnothing$ $(i\neq j)$, and
$\langle\alpha,\beta\rangle$ is a free two-generated commutative
semigroup with $\alpha$ and $\beta$ generators. Equip $S$ with a
binary operation defined as follows: 1 is an identity element of
$S$, $n_i\cdot m_j=(n+m)_{\min(i,j)}$, $\alpha n_i=n_i
\alpha=(n+3)_i$, $\beta n_i=n_i \beta = (n+2)_i$, and
$(\gamma_1\gamma_2)n_i=\gamma_1(\gamma_2n_i)$, where
$\gamma_1,\gamma_2\in \langle\alpha,\beta\rangle$, $i,j\in\omega$.
It is easy to verify that $S$ is a monoid under the operation
given, $\{1,0_0,0_1,\ldots,0_i,\ldots\}$ is the set of all
idempotents of $S$ and $S\cdot n_i=n_i\cdot
S=\bigcup\limits_{j=0}^iZ_j$ for all $i\in\omega$, $n\in Z$.
Consequently, $R=\bigcup\limits_{i\in\omega}Z_i$ and $S$ is a
regularly linearly ordered monoid which does not satisfy the
condition of Theorem 7.1.

We claim that $S$ is an $\mathfrak R$--stabilizer. It suffices to
prove that ${\rm Th}({}_SA)$ is stationary for any ${}_SA\in
{\mathfrak R}$. Let ${}_SA\in {\mathfrak R}$, ${}_SM\equiv{}_SA$,
$c=sb$, $a=tb$, $a,b\in{{\mathfrak C}} \setminus M$, and $c\in M$.
If $s,t\in R$, that is, $s=n_i$ and $t=m_j$, then the equalities
$n_i=(n-m)_im_j$ and $c=(n-m)_im_jb=(n-m)_ia$ hold for all
$i\leqslant j$, and $c\in Sa$; for $i>j$, likewise we derive $a\in
Sc$, which is impossible. Let $s\in\langle\alpha,\beta\rangle$.
Note that
\begin{equation}
\tag{7.2} S\cdot0_i\models\forall x\exists ! y(x=sy)
\end{equation}
for each $i\in \omega$. Let us prove that ${}_SA\models\forall
x\exists ! y(x=sy)$. In fact, on Corollary 2.1 for any $d\in A$
there exist $i\in\omega$ such that ${}_SSd\cong {}_SS0_i$. Since
the formula $\forall x\exists ! y(x=sy)$ is true in ${}_SS0_i$
then this formula is true in ${}_SSd$. Hence it is true in
${}_SA$. Let $d=sd_1$ and $d=sd_2$ for some $d,d_1,d_2\in A$. On
Corollary 2.1 there exist $i\in\omega$ and isomorphism
$\varphi:Sd_1\rightarrow S\cdot0_i$. Then $\varphi d_1,\varphi
d\in Z_i$ and $S\varphi d_1=S\varphi d$. Consequently, $Sd_1=Sd$.
Similarly $Sd_2=Sd$, that is, $Sd_1=Sd_2\cong S\cdot0_i$. In view
of (7.2), $d_1=d_2$. Since ${}_SM\equiv{}_SA$, it follows that
${}_SM\models\exists !y(c=sy)$ and ${}_S{\mathfrak
C}\models\exists !y(c=sy)$; on the other hand, $c=sb$, $b\not\in
M$, which is impossible. Thus $s\in R$ and
$t\in\langle\alpha,\beta\rangle$. Let $s=n_i$. There exists a
$k\in\omega$, $k\geqslant 2$, such that ${}_SR\models\forall
x(n_ix=(n-k)_itx)$. In view of regularity of polygon ${}_SA$ on
Corollary 2.1 ${}_SA\models\forall x(n_ix=(n-k)_itx)$. Since
${}_SA\equiv{}_S{\mathfrak C}$, it follows that ${}_S{{\mathfrak
C}}\models\forall x(n_ix=(n-k)_itx)$; in particular,
$c=sb=n_ib=(n-k)_itb=(n-k)_ia$, that is, $c\in Sa$. Consequently,
${Th}({}_SA)$ is stationary, and $S$ is an ${\mathfrak
R}$--stabilizer in view of Theorem 3.7 \hfill$\Box$

\vskip 1cm {\bf\centerline {\S~8. Superstability of Class for
Regular Polygons}}
\medskip

In this paragraph we give the characterization of an $\cal
R$--superstabilizer whose regular core is presented as the union
of finite number of principal right ideals (Theorem 8.1). As in
stable case it may be weakened to the condition when the regular
core is presented as a union of finite number of principal right
ideals replacing it by the axiomatizing of the class $\mathfrak R$
(Corollary 8.1). Since the class regular polygons over group is
axiomatizable, it is superstable (Corollary 8.2). Proposition 8.1
gives necessary condition for the superstability  of all regular
polygons. But as Example 8.1 shows, this condition is not
sufficient. In this paragraph we also construct the example of an
$\mathfrak R$--superstabilizer, but not stabilizer (Example 8.3)
and the example of $\mathfrak R$--stabilizer, but not an
$\mathfrak R$--superstabilizer (Example 8.4).

{\bf Proposition 8.1.} If $S$ is an $\mathfrak R$--superstabilizer
then $S$ is a regularly linearly ordered monoid and the semigroup
$Sa$ satisfies the ascending chain condition for left ideals,
where $a\in R$.

{\bf Proof.} Let $S$ be an $\mathfrak R$--superstabilizer. Then
$S$ is an $\mathfrak R$--stabilizer, which is regularly linearly
ordered by Proposition 7.1.
 Assume  $a,a_i\in R$ are such that
$Sa_i\subset Sa_{i+1}\subseteq Sa$, $a_i=s_ia$ $(s_i\in S$,
$i\in\omega)$. Let $T$ be the theory of $\mathfrak R$; $\kappa$ be
an arbitrary cardinal, $\kappa>2^{|T|}$;
$Q=\{\eta\in\kappa^\omega\mid\exists n<\omega$ $\forall m>n
(\eta(m)=0)\}$; $\hat{0}\in\kappa^\omega$ is such that
$\hat{0}(m)=0$ for each $m\in\omega$, and
$r(\eta)=\min\{n\in\omega\mid\forall m\geqslant n (\eta(m)=0)\}$.
On $Q$, we define the following relation:
$$ \eta<\varepsilon\Leftrightarrow r(\eta)\leqslant r(\varepsilon);$$
moreover, if $r(\eta)=r(\varepsilon)$, then there exists a
$k\in\omega$ such that $\eta|_k=\varepsilon|_k$, but $\eta
(k)<\varepsilon (k)$. The set $Q$ equipped with a relation
$\leqslant$ conforming to the relation $<$ above is well ordered.
Denote the element from $\kappa^{k-1}$ by $\eta_k$ if
$\eta_k|_k=\eta|_k$, where $\eta\in Q$, $k\in\omega$. We put
$\eta_0=\varnothing$.
 For
any $\eta\in Q$ and $k\in\omega$, we construct an $S$-polygons
${}_SN_\eta$, ${}_SN'_\eta$ and elements $b_\eta$, $b_{\eta_k}\in
N_\eta$.
  Let $N_{\hat{0}}=Sa$, $b_{\hat{0}}=a$, and
$b_{\hat{0}_k}=a_k$, $k\in\omega$. Assume the polygons
${}_SN_\xi$, ${}_SN'_\xi$ and the elements $b_\xi$, $b_{\xi_k}$
are constructed for all $\xi\in Q$, $\xi<\eta$, $k\in\omega$. Note
that element $b_{\eta_{r-1}}$ has constructed on the previous
steps. If it exists the largest element $\varepsilon\in Q$ such
that $\varepsilon<\eta$, then ${}_SN'_\eta={}_SN_\varepsilon$; in
the other ways,
${}_SN'_\eta=\bigcup\limits_{\varsigma<\eta}{}_SN_\varsigma$. We
put ${}_SN_\eta ={}_S(N'_\eta\sqcup Sa)/\Theta_\eta$, where
$\Theta_\eta$ is a congruence on the polygon ${}_S(N'_\eta\sqcup
Sa)$ generated by the pair $\langle
b_{\eta_{r-1}},a_{r-1}\rangle$; in this case elements of $N'_\eta$
are identified with congruence classes $\eta$ whose
representatives they are of. Redenote the element $a/\theta_\eta$
by $b_\eta$, and elements $a_i/\theta_\eta$ by $b_{\eta_i}$, where
$i>r-1$. Write ${}_SN$ for $\bigcup\limits_{\eta\in Q}{}_SN_\eta$.
Clearly, ${}_SN\in{\mathfrak R}$. Let
$A=\{b_{\eta_k}\mid\eta\in\kappa^\omega,k\in\omega\}$. Since
$$ b_{\eta_k}=s_kx\in {\rm
tp}(b_\varepsilon,A)\Leftrightarrow\eta(\xi)=\varepsilon(\xi)
\;\mbox{ for all }\;\xi\leqslant k, $$
 where $k\geqslant0$, we have ${\rm tp}(b_\eta,A)\neq {\rm
tp}(b_\varepsilon,A)$, $\eta\neq\varepsilon$, $\eta,\varepsilon\in
Q$. Since $|A|= \sum\limits_{k\in\omega}\kappa^k$, we obtain $|
S(A)|\geqslant|\{b_\eta\mid\eta\in\kappa^\omega\}|
=\kappa^\omega>\kappa$. Thus ${\rm Th}({}_SN)$ is not superstable,
and consequently $S$ fails as an $\mathfrak
R$--superstabilizer.\hfill$\Box$

\smallskip
{\bf Theorem 8.1} [Ste2] Let
$$ R=\bigcup\limits_{i=0}^n a_iR $$
for some $n\geqslant 0$ and $a_i\in R$ $(0\leqslant i\leqslant
n)$. The monoid $S$ is an $\mathfrak R$--superstabilizer if and
only if $S$ is a regularly linearly ordered monoid and  the
semigroup $Sa$ satisfies the ascending chain condition for left
ideals, where $a\in R$.

{\bf Proof.} Assume  the theorem condition is satisfied.

{\bf Necessity} follows from Proposition 7.1.

{\bf Sufficiency.}
 Let $S$ be a regularly linearly ordered monoid,  the semigroup
$Sa$ satisfies the ascending chain condition for left ideals,
where $a\in R$; ${}_SA\in{\mathfrak R}$; ${}_SM\equiv{}_SA$;
$a\in{{\mathfrak C}} \setminus M$; $Sa\cap M=\varnothing$.

We claim that there exists an entering element from $a$ into $M$,
that is, there exists a $c\in M$ for which $c\in Sa$ and
$Sb\subseteq Sc$ with all $b\in M\cap Sa$. Let $d$ be an arbitrary
element of $A$ and $\Phi\rightleftharpoons\forall
x\bigvee\limits_{i\leqslant n}\exists y(x=a_iy)$. Since ${}_SA$ is
regular on Corollary 2.1 it exists an idempotent $e\in R$ such
that $Sd\cong Se$. Under conditions theorem ${}_SSe\models\Phi$.
Consequently, ${}_SSd\models\Phi$. In view of the arbitrary of
element $d$ we have ${}_SA\models\Phi$. Since the polygons
${}_SA,\;{}_SM$ and ${}_S{\mathfrak C}$ are elementary equivalents
it follows that ${}_SM\models\Phi$ and ${{}_S{\mathfrak
C}}\models\Phi$. Consequently, $a=a_ia'$ for some $i\leqslant n$
and $a'\in{\mathfrak C}$. Assume $m,m'\in Sa\cap M$. Then
$m=sa=sa_ia'$ and $m'=ta=ta_ia'$ for some $s,t\in S$. Since $S$ is
a regularly linearly ordered monoid, either $Ssa_i\subseteq Sta_i$
or $Sta_i\subseteq Ssa_i$. Hence either $sa_i=rta_i$ or
$ta_i=rsa_i$, $r\in S$, that is, either $m=rm'$ or $m'=rm$.
Consequently, either $Sm\subseteq Sm'$ or $Sm'\subseteq Sm$.
Suppose that there exists no entering element from $a$ into $M$,
which means that there exist $m_i\in Sa\cap M$, $i\in\omega$, such
that $Sm_i\subset Sm_{i+1}$ and any polygon ${}_SSm,\;m\in Sa\cap
M$, coincides with some polygon ${}_SSm_i$. For any $j\in\omega$,
there then exists an $s_j\in S$ for which $m_j=s_ja=s_ja_ia'$. Let
$Ss_{j+1}a_i\subseteq Ss_{j}a_i$, that is $s_{j+1}a_i=ks_ja_i$ for
some $k\in S$. Hence $s_{j+1}a_ia'=ks_ja_ia'$. Then
$s_{j+1}a=ks_ja$, that is $m_{j+1}=km_j$ and $Sm_{j+1}\subseteq
Sm_j$, which contradicts the assumption. In view of the regularity
linear ordering $S$ we claim that $Ss_ja_i\subset Ss_{j+1}a_i$,
which contradicts  the ascending chain condition for left ideals
of semigroup $Sa_i$. We have thus proved that the desired entering
element from $a$ into $M$ exists.

Let $m$ be an entering element from $a$ into $M$,
$a,b\in{\mathfrak C}$. On Corollary 7.4, the theory ${\rm
Th}({}_SM)$ is stationary. On Theorem 3.6,
$$ {\rm tp}(a,M)={\rm tp}(b,M)\Leftrightarrow {\rm tp}(a,\{m\})={\rm tp}(b,\{m\}),$$
 and $m$ is an entering element from $b$ into $M$. Let $T={\rm
Th}({}_SM)$, $|M|=\kappa\geqslant2^{| T|}$, $Q_0=\{{\rm
tp}(a,M)\mid a\in M\}$, $Q_1=\{{\rm tp}(a,M)\mid a\in{\mathfrak
C}\setminus M,M\cap Sa=\varnothing\}$, and $Q_2=\{{\rm
tp}(a,M)\mid a\in{\mathfrak C} \setminus M,M\cap
Sa\neq\varnothing\}$. Then $|Q_0|=\kappa$;
$|Q_1|\leqslant|S(\varnothing)|\leqslant 2^{|T|}\leqslant\kappa$;
$|Q_2|\leqslant|M|\cdot2^{|T|}=\kappa\cdot2^{|T|}=\kappa$, that
is, $|S(M)|=| Q_0\cup Q_1\cup Q_2|=\kappa$. Consequently, $T$ is
superstable and $S$ is an  ${\mathfrak
R}$--superstabilizer.\hfill$\Box$

Note that  the ascending chain condition for left ideals of $Sa$
($a\in R$) in Theorem 8.1 may be replaced by the ascending chain
condition for left ideals of $Sa_i$ ($1\leqslant i\leqslant n$).
It follows from the proof of Theorem 8.1.

\smallskip
{\bf Corollary 8.1.} Let the class  ${{\mathfrak R}}$ for regular
polygons is axiomatizable. The class ${{\mathfrak R}}$ is
superstable if and only if $S$ is a regularly linearly ordered
monoid and  the semigroup $Sa$ satisfies the ascending chain
condition for left ideals, where $a\in R$.

The {\bf proof} follows from Theorem 8.1 and Corollary
4.2\hfill$\Box$

\smallskip
{\bf Corollary 8.2.} The class  ${{\mathfrak R}}$ for regular
polygons over group is superstable.

The {\bf proof} follows from  Corollaries 4.3 and 8.1.\hfill$\Box$

\smallskip
{\bf Example 8.1.} Example 7.1 exemplifies a regularly linearly
ordered, non--${{\mathfrak R}}$--superstabilizer $S$ such that the
semigroup $Sa$ satisfies the ascending chain condition for left
ideals, where $a\in R$. \hfill$\Box$

\smallskip
{\bf Example 8.2.} We construct a regularly linearly ordered
superstabilizer $S$, such that the semigroup $Sa$ satisfies the
ascending chain condition for left ideals, where ($a\in R$), for
which the assumption of Theorem 8.1 fails.

Let $S=\omega\cup\langle\alpha\rangle\cup\{1\}$, where
$\langle\alpha\rangle$ is a one-generated free semigroup with
generator $\alpha$. Equip $S$ with a binary operation defined as
follows: 1 is an identity element of $S$; $n\cdot m=n$ and
$\alpha^i\cdot n=n\cdot\alpha^i=n$ for any $n,m\in\omega$,
$i\geqslant 1$. It is easy to verify that $S$ is a monoid under
the operation given, $Sn=\omega=R$ and $nS=\{n\}$, where
$n\in\omega$. Then $S$ is a regularly linearly ordered monoid for
which the assumption of Theorem 8.1 fails and the semigroup $Sa$
satisfies the ascending chain condition for left ideals, where
$a\in R$. Since $Sm\subset S\alpha^{i+1} S\subset\alpha^i\subset
S$ for any $m\in\omega,\;i\geqslant 1$ then $S$ is a
superstabilizer by Theorem 3.9\hfill$\Box$

\smallskip
{\bf Example 8.3.} We construct an ${{\mathfrak R}}$--stabilizer
$S$ which is not a stabilizer.

The monoid $S$ is obtained from the monoid specified in Example
7.2 by replacing the ordinal $\omega$ by $\{0\}$ and preserving
the operation. It follows that $R=Z_0=0_0\cdot R$, and by Theorem
8.1, $S$ is a regular superstabilizer. Since $S$ is not linearly
ordered, $S$ is ${{\mathfrak R}}$--superstabilizer. Since $S$ is
not linearly ordered then by Theorem 3.8 $S$ is not a
stabilizer.\hfill$\Box$

\smallskip
{\bf Example 8.4.} We construct an ${{\mathfrak R}}$--stabilizer
$S$ which is not an ${{\mathfrak R}}$--superstabilizer.

The monoid $S$ is obtained via the monoid in Example 7.2 by
replacing the ordinal $\omega$ by $\omega\cup\{\omega\}$ and
preserving the operation. It follows that
$R=\bigcup\limits_{i\in\omega\cup\{\omega\}}Z_i=Z_\omega$ and
$S0_i\subset S0_{i+1}\subset S0_{\omega}$ for any $i\in\omega$.
 That
is, $S$ fails as an ${{\mathfrak R}}$--superstabilizer by
Proposition 8.1. That $S$ is an ${{\mathfrak R}}$--stabilizer is
explicated as in Example 7.2.\hfill$\Box$

\vskip 1cm

\medskip {\bf\centerline {\S~9. $\omega$--stability of Class for
Regular Polygons}}

\medskip

In [FG] there were considered the questions of the stabilities of
the theories, which are the model companion of the theory of all
polygons. Such theory exists if a monoid $S$ is left coherent. In
this work it is proved that each complete type of this theory is
characterized by triple consisting of a left ideal of the monoid
$S$, a left congruence on the monoid $S$ and a polygon
homomorphism. In this paragraph in the case of an axiomatizability
and a model completeness of the class ${\mathfrak R}$ each
complete type of the theory of ${\mathfrak R}$ is characterized by
the triple consisting of a left ideal of the semigroup $Se_i$,
left congruency of the semigroup $Se_i$ and a polygon
homomorphism, where the idempotents $e_i$ are taken from Corollary
4.2 (Lemma 9.1). It is this result which allows us to translate
model theoretic properties of the class ${\mathfrak R}$ into
algebraic properties of the monoid $S$. The consequence is that we
can easily find upper bounds for the number of the full types of
the class ${\mathfrak R}$ theory. Theorem 9.1 gives us a criterion
of superstability and $\omega$--stability of axiomatizable and a
model complete class ${\mathfrak R}$ of regular polygons. In this
theorem the stability of the class ${\mathfrak R}$ it proved. At
the end of the paragraph there is an example of  ${\mathfrak
R}$--superstabilizer, but not $\omega$--${\mathfrak
R}$--stabilizer and not superstabilizer.

Assume ${}_SA\in{\mathfrak R}$, $R=\bigcup \{ e_iR\mid 1\le i\le
n\}$, where $n\in\omega$, $e^2_i=e_i\in R$ $(1\le i\le n)$. Write
$Tr^i(A)$ $(1\le i\le n)$ for the set of the triples
$\langle\theta,I,\alpha\rangle$, satisfying the following
conditions:

(1) $\theta$ is a congruency of a polygon ${}_SSe_i$;

(2) ${}_SI$ is a subpolygon of a polygon ${}_SSe_i$; $Se_i$;

(3) $\alpha : {}_SI \to {}_SA$, $\alpha$ is a polygon
homomorphism;

(4) ${}_SI$ is a $\theta$--saturated polygon that is $\langle
a,b\rangle \in \theta$ and $a\in I$ imply $b\in I$;

(5) $Ker \alpha = \theta \cap (I\times I)$;

(6) ${}_S(Se_i/\theta)\in {\mathfrak R}$.

Suppose $Tr(A)=\cup \{Tr^i(A)\mid 1\le i\le n\}$.

\smallskip {\bf Lemma 9.1}. Let ${\mathfrak R}$ be an axiomatizable and
model complete class, $R=\bigcup \{ e_iR\mid 1\le i\le n\}$, where
$n\in\omega$, $e^2_i=e_i\in R$ $(1\le i\le n)$, ${}_SA\in
{\mathfrak R}$. There exists a surjective mapping
$\varphi:Tr(A)\to S_1(A)$ such that $\varphi |_{ Tr^i(A)}$ is an
injective mapping for any $i$, $1\le i\le n$, where $S_1(A)$ is
the set of 1--types over $A$.

{\bf Proof}. Assume ${}_SA\in {\mathfrak R}$ and
$\langle\theta,I,\alpha\rangle\in Tr^i(A)$ for some
$i,\;1\leqslant i\leqslant n$. Write $p'(x)$ for the following
type over $A$:
$$\{sx=a\mid \alpha(se_i) =a,\;se_i\in I\}\cup\{sx\neq a\mid se_i\not\in I,\;a\in
A\}\cup$$
$$\cup\{ sx=tx\mid\langle se_i,
te_i\rangle\in\theta\}\cup\{ sx\neq tx\mid\langle se_i,
te_i\rangle\not\in\theta\}.$$ We claim that the type $p'(x) $ is
consist with the theory $T$ of the class ${{\mathfrak R}}$. We can
assume that $A\cap Se_i=\emptyset$ (otherwise we replace $A$ by
the isomorphic copy $A',\;A'\cap Se_i=\emptyset$, a homomorphism
$\alpha$ by $\alpha'=\beta\alpha$, where $\beta:A\longrightarrow
A'$ is an isomorphism). On the set $A\cup Se_i$ we define a
relation $\sigma$:
$$\sigma=\{\langle \alpha(t),t\rangle\mid t\in I\}.$$
Suppose $\tau=\varepsilon\cup\sigma\cup\sigma^{-1}\cup\theta$,
where $\varepsilon$ is a zero congruence of the polygon
${}_S(A\cup Se_i)$. We claim that $\tau$ is a congruence of the
polygon ${}_S(A\cup Se_i)$. For this it is enough to prove that
$\tau$ is a transitive relation. Assume $x,y,z\in A\cup
Se_i,\;x\tau y\tau z,\;x\neq y,\;y\neq z$.

If $x\sigma y$ then $x=\alpha(y),\;y\in I$. So $z=\alpha(y)$ or
$y\theta z$. In the first case we have $x=y$, in the second case
in view of $\theta$--saturation of the polygon ${}_SI$ we have
$z\in I$, and the equality $Ker\alpha=\theta\cap(I\times I)$
implies $\alpha(y)=\alpha(z)$, that is, $x\sigma z$.

If $x\sigma^{-1} y$ then $y=\alpha(x),\;y\in A$, consequently,
$y=\alpha(z)$, $x,y\in I$ and $x\theta z$.

If $x\theta y$ then either $y\theta z$, or $z=\alpha(y),\;y\in I$,
that is, $z=\alpha(x)$ and $z\sigma x$.

Thus, $\tau$ is a congruence of the polygon ${}_S(A\cup Se_i)$.
Let ${}_SC={}_S(A\cup Se_i)/\tau,\;c=e_i/\tau$. Then an element
$c$ realizes the type $p'(x)$ in ${}_SC$. Furthermore, since
${}_SA\in{{\mathfrak R}}$ and ${}_S(Se_i/\theta)\in{{\mathfrak
R}}$ it follows that ${}_SC\in{\mathfrak R}$. Consequently, the
type $\Sigma(x)=p'(x)\cup T\cup D({}_SA)$ is realized by element
$c\in C$,  where $D({}_SA)$ is a diagram of the polygon ${}_SA$.
Hence this type is consistent.

We claim that the type $p(x)=tp(c,A)$ of an element $c$ over a set
$A$ is a unique complete type, containing the type $\Sigma(x)$. It
is clear that $\Sigma(x)\subseteq p(x)$. Let $\Phi(x)$ be an
atomic formula of the language $L_S(A)$, which is an enriching the
language $L_S$ by adding in constant symbols for all elements $A$.
If $\Phi(x)$ is a sentence then $\Phi(x)\in D({}_SA)$ or
$\neg\Phi(x)\in D({}_SA)$. If $x$ is a free variable of the
formula $\Phi(x),\;\Phi(x)\not\in\Sigma$ and
$\neg\Phi(x)\not\in\Sigma$ then the definition of $p'(x)$ implies
that $\Phi(x)$ has the form of $sx=a$, where $se_i\in
I,\;\alpha(se_i)\neq a,\;a\in A$. Since $se_i\in I$ it follows
that $\alpha(se_i)=b\in A$, $b\neq a$, consequently,
$(sx=b)\in\Sigma(x)$ and $\Sigma(x)\models\neg\Phi(x)$. Thus for
any quantifier--free formula $\psi(x)$ either
$\Sigma(x)\models\psi(x)$ or $\Sigma(x)\models\neg\psi(x)$. Since
$T$ is a submodel complete theory, $T$ has elimination of
quantifiers by Theorem 3.3. Consequently, $\Sigma(x)\models p(x)$
for any type $p(x)\in S_1(A)$. We correspond the type $p$ to the
triple $\langle\theta,I,\alpha\rangle\in Tr^i(A)$. The mapping
$\varphi$ is built. Obviously, $\varphi|_{Tr^i(A)}$ is an
injective mapping for any $i,\;1\leqslant i\leqslant n$.

We claim that $\varphi$ is a surjective mapping. Let $p(x)\in
S_1(A)$; ${}_SB$ be a saturated model of $T$, realizing the type
$p(x)$ and being elementary extension of ${}_SA$; element $b_0\in
B$ realizes the type $p(x)$. Since ${\mathfrak R}$ is
axiomatizable class and ${}_SA\in{\mathfrak R}$,
${}_SB\in{\mathfrak R}$ and $_SSb_0\in{\mathfrak R}$. On Corollary
2.1 there exist an idempotent $f\in R$ and an isomorphism
$\eta:{}_SSb_0\longrightarrow{}_SSf$ such that $\eta(b_0)=f$.
Since $f\in R=\bigcup \{ e_iR\mid 1\le i\le n\}$ it follows that
$f\in eS$, that is, $f=ef$ for a certain $e\in\{e_1,\ldots,e_n\}$.
It is easy to understand that for any $s,r\in S$
\begin{equation}\tag{9.1}
  se=re\;\Longrightarrow\;sef=ref\;\Longleftrightarrow\;sf=rf\;\
  Longleftrightarrow\;sb_0=rb_0\;\Longleftrightarrow\;(sx=rx)\in
  p(x).
\end{equation}

Let $$\theta=\{\langle re,se\rangle\in(Se)^2\mid(rx=sx)\in
p(x),\;r,s\in S\},$$ $$I=\{re\in Se\mid(rx=a)\in p(x),\;r\in
S,\;a\in A\},$$
$$\alpha: I\longrightarrow A,$$ where $\alpha(re)=a$ for all $re\in I$
such that $(rx=a)\in p(x),\;a\in A$.

In view of the completeness of type $p(x)$ and (9.1) it is easy to
prove conditions (1)-(5) for the triple
$\langle\theta,I,\alpha\rangle$. Condition (6) follows from the
correlation ${}_SSe/\theta\cong{}_SSb_0\in{\mathfrak R}$. Clearly,
$\varphi(\langle\theta,I,\alpha\rangle)=p(x)$. Lemma is
proved.\hfill$\Box$

\smallskip {\bf Theorem 9.1} [Ste3]. Let the class ${\mathfrak R}$
for all regular polygons be axiomatizable and model complete. Then
$|R|\geqslant \omega$, $R=\bigcup\{e_iR|1\leqslant i\leqslant n\}$
for some $n\in\omega$, $e_1,\ldots,e_n\in R$, $e^2_i=e_i$
$(1\leqslant i\leqslant n)$, and

1) the class ${\mathfrak R}$ is stable;

2) the class ${\mathfrak R}$ is superstable if and only if for any
$i$, $1\leqslant i\leqslant n$, a semigroup $Se_i$ satisfies the
ascending chain condition for left ideals;

3) if a semigroup $R$ is  countable, then the class ${\mathfrak
R}$ is $\omega$--stable if and only if for any $i$, $1\leqslant
i\leqslant n$, a semigroup $Se_i$ satisfies the the ascending
chain condition for left ideals and has no more than countable
number of the left congruence $\theta$, such that
${}_S(Se_i/\theta)\in{\mathfrak R}$.

{\bf Proof}. Let the class ${\mathfrak R}$ for all regular
polygons be axiomatizable and model complete. Using Theorem 5.1
(3) we have $|R|\geqslant \omega$. The presentation of the
semigroup $R$ as the union of the finite number of right ideals is
established in Corollary 4.2.

Let ${}_SA\in{{\mathfrak R}}$. Note that the number of all
subpolygons ${}_SI$ of polygon ${}_SSe_i\;(1\leqslant i\leqslant
n)$ is not greater than $2^{|R|}$, the number of all left
congruences $\theta$ of semigroup $Se_i\;(1\leqslant i\leqslant
n)$ is not greater than $2^{|R|^2}$, the number of all
homomorphisms from the subpolygon ${}_SI$ of the polygon
${}_SSe_i\;(1\leqslant i\leqslant n)$ to the polygon ${}_SA$ is
not greater than $2|A|^{|R|}$. Thus,
$|Tr(A)|\leqslant2^{|R|}\cdot2^{|R|^2}\cdot|A|^{|R|}$.

If $|A|\leqslant2^{|R|}$ then $|Tr(A)|\leqslant2^{|R|}$ and, on
Lemma 9.1, $|S_1(A)|\leqslant2^{|R|}$. Consequently, the class
${\mathfrak R}$ is stable and 1) is proved.

Suppose for any $i$, $1\leqslant i\leqslant n$, the semigroup
$Se_i$ satisfies the ascending chain condition for left ideals.
Then each subpolygon ${}_SI$ of the polygon ${}_SSe_i\;(1\leqslant
i\leqslant n)$ is finite--generated. Consequently, the number of
such subpolygons is not greater than $|R|$ and the number of all
homomorphisms from the subpolygon ${}_SI$ of the polygon
${}_SSe_i\;(1\leqslant i\leqslant n)$ to the polygon ${}_SA$ is
not greater than ${\omega\cdot |A|\cdot |R|}$. So
$|Tr(A)|\leqslant{|A|\cdot |R|\cdot \omega\cdot|A|^{|R|}}$. If
$|A|\geqslant 2^{|R|}$ then $|Tr(A)|\leqslant{|A|}$ and on Lemma
9.1 $|S_1(A)|\leqslant{|A|}$. Consequently, the class ${\mathfrak
R}$ is superstable.

Assume $|R|=\omega$ and for any $i$, $1\leqslant i\leqslant n$,
the semigroup $Se_i$ satisfies the ascending chain condition for
left ideals and has no more than countable number of the left
congruence $\theta$ such that ${}_S(Se_i/\theta)\in{\mathfrak R}$.
Then, on proved above, $|Tr(A)|\leqslant{|A|\cdot |R|\cdot
\omega}$. If $|A|=\omega$ then a countability of the semigroup $R$
implies $|Tr(A)|\leqslant{\omega}$ and on Lemma 9.1
$|S_1(A)|\leqslant{\omega}$. Consequently, the class ${\mathfrak
R}$ is $\omega$--stable.

Let the class ${\mathfrak R}$ be superstable. Since the class
${\mathfrak R}$ is axiomatizable, the monoid $S$ is ${\mathfrak
R}$--superstable and, on Proposition 8.1, for any $i$, $1\leqslant
i\leqslant n$, the semigroup $Se_i$ satisfies the ascending chain
condition for left ideals. Thus 2) is proved.

Let the class ${\mathfrak R}$ be $\omega$--stable. Then on Theorem
3.4 the class ${\mathfrak R}$ is superstable. Let us show that for
any $i$, $1\leqslant i\leqslant n$, the semigroup $Se_i$ has no
more than countable number of the left congruences $\theta$ such
that ${}_S(Se_i/\theta)\in{\mathfrak R}$. Since
${}_SR\in{\mathfrak R}$ and $|R|=\omega$, it follows that
$|S_1(R)|\leq\omega$. On Lemma 9.1, $|Tr^i(R)|\leq|S_1(R)| $ for
any $i,\;1\leqslant i\leqslant n$, that is $|Tr^i(R)|\leq\omega$.
Let $\theta$ be a left congruence of the semigroup ${}_SSe_i$ such
that ${}_S(Se_i/\theta)\in{\mathfrak R}$. Suppose $I=Se_i$. Since
${}_S(Se_i/\theta)={}_SS(e_i/\theta)\in{\mathfrak R}$ it follows
that on Corollary 2.1 there exist an idempotent $f\in R$ and an
isomorphism $\alpha':{}_SS(e_i/\theta)\longrightarrow{}_SSf$. We
construct a homomorphism $\alpha:{}_SI\longrightarrow{}_SR$ as
follows: $\alpha(a)=\alpha'(a/\theta)$ for any $a\in Se_i$. Then
$\langle\theta,I,\alpha\rangle\in Tr^i(R)$ and the number of the
congruences of polygon ${}_SSe_i$ is not greater than $Tr^i(R)$.
Since $|Tr^i(R)|\leq\omega$, the semigroup $Se_i$ has no more than
countable number of the left congruences $\theta$ such that
${}_S(Se_i/\theta)\in{\mathfrak R}$ for any $i,\;1\leqslant
i\leqslant n$. Thus, 3) is proved. \hfill$\Box$

Note that the statements 1 and 2 of Theorem 9.1 we can derive as a
consequence of Theorems 7.1 and 8.1 accordingly.

\smallskip {\bf Corollary 9.1}. The class of regular polygons over a
countable semigroup is $\omega$--stable.

{\bf Proof.} On Corollary 5.1, the class of regular polygons over
countable group is axiomatizable and model complete. In a group
$S$ a unique idempotent is unit, a unique left ideal is a group
$S$ itself. Let $\theta$ be a left congruence of the group $S$
such that ${}_S(S/\theta)\in{\mathfrak R}$. Then $S/\theta=S\cdot
1/\theta$. On Corollary 2.1 there exists an isomorphism
$\varphi:S\cdot 1/\theta\longrightarrow S$ such that
$\varphi(1/\theta)=1$. If $t\in1/\theta$ then $t\cdot
1/\theta=1/\theta$, consequently, $t=1$ and $|1/\theta|=1$. Thus,
on Theorem 9.1, the class of regular polygons over a countable
semigroup is $\omega$--stable.\hfill$\Box$

We will use a formula $\exists^nx\Psi(x,x_1,\ldots,x_n)$ as an
abbreviation of a formula "there exist exactly $n$ elements $x$
such that $\Psi(x,x_1,\ldots,x_n)$".

\smallskip {\bf Lemma 9.1}. Let $a,b,c\in R$ be such that
$Sc\subset Sb\subset Sa$, $b=\alpha a,\;c=\beta b$ and there exist
a formula $\Phi(x,y,z)$ and $n\in\omega,\;n>0$, such that
$\Phi(S,S,c)=\Phi(Sa\setminus Sb,Sb\setminus Sc,c)$ and
${}_SSa\models\Phi(a,b,c)\wedge\forall y(\beta y= c\wedge\exists
x(\Phi(x,y,c)\wedge\alpha x=y)
\rightarrow\exists^nx(\Phi(x,y,c)\wedge \alpha x=y)).$ Then $S$ is
not an ${{\mathfrak R}}$--$\omega$--stabilizer.

\smallskip {\bf Proof.} Assume the conditions of lemma are
satisfied. Let $K\subseteq\omega$; ${}_SM_K=\bigcup\limits_{i\in
K}\bigcup\limits_{j\leqslant i}{}_SS\langle
a,j^i,K\rangle/\theta_K$; ${}_SM=\bigcup\limits_{K\subseteq
\omega}{}_SM_K$, where $S\langle a,j^i,K\rangle=\{\langle
x,j^i,K\rangle\mid x\in Sa\}$, $s\cdot \langle
x,j^i,K\rangle=\langle sx,j^i,K\rangle$ for all $x\in Sa,\;i\in
K,\;j\leqslant i$; $\theta_K$ be a congruence of the polygon
$\bigcup\limits_{i\in K}\bigcup\limits_{j\leqslant i}{}_SS\langle
a,j^i,K\rangle$, which is generated by a set $\{\langle\langle
b,j_1^i,K \rangle,\langle b,j_2^i,K\rangle\rangle\mid
j_1,j_2\leqslant i,i\in K \}$$\cup$$\{\langle\langle c,j_1^{i_1},K
\rangle,$$\langle c,j_2^{i_2},K\rangle\rangle\mid j_1\leqslant
i_1,j_2\leqslant i_2,i_1,i_2\in K \}$.  White $\Gamma_K(x)$ for a
set of the formulae $$ \{\exists y\exists^{n(i+1)}z(x=\beta
y\wedge y=\alpha z\wedge\Phi(z,y,x)\mid i\in K\}\cup$$
$$\cup \{{}^\neg\exists y\exists^{n(i+1)}z(x=\beta y\wedge
y=\alpha z\wedge\Phi(z,y,x)\mid i\not\in K\}.$$ It is not
difficult to understand that the set of the formulae $\Gamma_K(x)$
is realized by an element $\langle c,j^i,K\rangle/\theta_K$ and is
not realized by element $\langle c,j^i,K'\rangle/\theta_{K'}$ for
any $K'\subseteq \omega,\;K\neq K'$. Consequently,
$|S(\emptyset)|=2^\omega$, that is, $S$ is not ${{\mathfrak
R}}$--$\omega$--stabilizer.\hfill$\Box$

\smallskip {\bf Example 9.1.} We construct a non
${{\mathfrak R}}$--$\omega$--stabilizer and non stabilizer which
is ${\mathfrak R}$--su\-per\-stabilizer.

Let $S=\{a,b,c\} \cup\langle\alpha,\beta\rangle\cup\{1\}$, where
$\langle\alpha,\beta\rangle$ is a free two--generated commutative
semigroup with $\alpha$ and $\beta$ generators. On $S$ we define
the binary operation as follows: 1 unit in $S$, $a\cdot
x=a,\;b\cdot x=x$ for any $x\in \{a,b\}$, $c\cdot y=y\cdot c=y$
for any $y\in\{a,b,c\}$, $\gamma\cdot z=z\cdot \gamma=z$ for any
$z\in\{a,b,c\},\;\gamma\in\langle\alpha,\beta\rangle$. It is easy
to check that $S$ equipped with the operation is a monoid,
$\{a,b,c\}$ is a set of all idempotents of monoid $S$, $Sa\subset
Sb\subset Sc=S$ and $R=Sc=cS$. On Theorem 8.1 $S$ is ${\mathfrak
R}$--superstabilizer. For the elements $a,b,c$, a formula
$\Phi(x,y,z)\rightleftharpoons bx=y\wedge x\neq y$, and an $n=1$
the conditions of Lemma 9.4 are hold. Consequently, $S$ is non
${{\mathfrak R}}$--$\omega$--stabilizer. Since a monoid $S$ is not
linear ordered, $S$ is non stabilizer.\hfill$\Box$

\vskip 1cm

\medskip {\bf\centerline {References}}

\medskip

[ChK] Chang C.C., Keisler H.J. Model theory. North--Holland
Publishing Company, 1973.

[EP] Ershov Yu.L., Palyutin E.A. Mathematical Logic. M.: Nauka,
1987.

[FG] Fountain J.B., Gould V. Stability of S-sets. Preprint. 2001.

[KM] Kargapolov M.I., Merzlyakov Yu.I. Foundation of Group Theory.
-- M.: Nauka, 1982.

[KKM] Kilp M., Knauer U., Mikhalev A.V. Monoids, Acts and
Categories. Walter de Gruyet. Berlin. New York. 2000.

[Mus] Mustafin T.G. On Stable Theory of polygons
// Model Theory and its Application.-- Novosibirsk: Nauka. Sib.
otd--е, 1988. P.92-107.

[Ov1] Ovchinnikova E.V. Structure of the class of regular polygons
// Studies in theories of the algebraic systems. Mezhvuzovskiy
sbornik nauchnyh trudov. -- Karaganda: edition of Karaganda State
University, 1995. -- P. 84--86.

[Ov2] Ovchinnikova E.V. Complete classes of regular polygons with
finitely many idempotents // Sib Math. J. 1995. V.36, N. 2. P.
381--384.

[Ov3] Ovchinnikova E.V. Complete classes of regular polygons over
monoids of the depth 2 // Algebra and logic 2002. V. 41, N. 6. P.
745-753.

[Ov4] Ovchinnikova E.V. Monoid, over which the class of regular
polygons is complete, not model complete
// Siberian Math. J. 1997. V.38, N. 5. P. 110--114.

[Ov5] Ovchinnikova E.V. Not linearly ordered monoid, over which
the class of regular polygons is complete, but not model complete
// Algebra and model theory, 3. Collection of the papers.
Novosibirsk: edition of NSTU, 2001. P. 83-98.

[Sac] Sacks G.E. Saturated model theory. W.A Benjamin, Reading,
Mass., 1972.

[She] Shelah S. Classification theory and the number of
non-isomorphic models, North-Holland, Second edition: 1990.

[Ste1] Stepanova A.A. Axiomatizability and Model Completeness of
Monoids with Stable Theories for Regular Polygons // Sib. Math. J.
1994. V. 35. N. 1. P. 181-193.

[Ste2] Stepanova A.A. Monoids with Stable Theories for Regular
Polygons // Algebra and logic. 2001. V. 40. N. 4. P. 239-254.

[Ste3] Stepanova A.A. Stability of Class for Regular Polygons.
// Studies in theories of the algebraic systems. Mezhvuzovskiy
sbornic nauchnih trudov. -- Karaganda: edition of Karaganda State
University, 1995. -- P.95--102.

[Su] Sushkevich A.K. The theory of generalised groups. -- Harkov;
Kiev: Gos. nauch.-tehn. idition Ukr., 1937.

[Tra] Tram L.H. Characterization of monoid by regular acts //
Period. Math. Hungar. 1985. V.16. P.273-279.

[Zel] Zelmanowitz J. Regular modules // Trans. Amer. Math. Soc.
1972. V.163. P.341--355.

\end{document}